\documentclass[a4paper, 11pt]{amsart}

\usepackage[format=hang,font=small,labelfont=bf]{caption}

\usepackage[utf8]{inputenc}
\usepackage[english]{babel}

\usepackage[OT2,T1]{fontenc} 

\usepackage{etoolbox}

\usepackage{amsfonts}
\usepackage{amssymb}
\usepackage{amsthm}
\usepackage{hyperref}

\usepackage{mathrsfs}
\usepackage{mathtools}
\usepackage{enumerate}

\usepackage[shortlabels]{enumitem}

\usepackage[all]{xy}

\usepackage{fullpage}
\usepackage{xcolor}

\usepackage{wrapfig}

\usepackage{verbatim}
\usepackage{parskip}

\usepackage{pinlabel}
\usepackage[toc, page]{appendix}

\usepackage{tkz-graph}
\usetikzlibrary{arrows}
\usetikzlibrary{calc}
\usepackage{tikz-cd}

\usepackage{multicol}
\usepackage{cancel}

\usepackage{epigraph}
\usepackage{transparent}
\usepackage{extarrows}

\usepackage{stmaryrd}

\usepackage{calligra}

\usepackage{multirow}
\usepackage{booktabs}

\usepackage{thm-restate}

\usepackage[colorinlistoftodos]{todonotes}

\usepackage{float}
\usepackage{aliascnt}
\newaliascnt{eqfloat}{equation}
\newfloat{eqfloat}{h}{eqflts}
\floatname{eqfloat}{Equation}

\newcommand*{\ORGeqfloat}{}
\let\ORGeqfloat\eqfloat
\def\eqfloat{%
	\let\ORIGINALcaption\caption
	\def\caption{%
		\addtocounter{equation}{-1}%
		\ORIGINALcaption
	}%
	\ORGeqfloat
}

\newtheorem{thm}{Theorem}[section]
\newtheorem{theorem}{Theorem}
\newtheorem*{theorem*}{Theorem}
\newtheorem{mainTheorem}{Theorem}

\newtheorem{cor}[thm]{Corollary}

\newtheorem{lemma}[thm]{Lemma}
\newtheorem{prop}[thm]{Proposition}

\newtheorem{conj}[thm]{Conjecture}

\theoremstyle{rmk}
\newtheorem{remark_basic}[thm]{Remark}
\newtheorem{observation_basic}[thm]{Observation}
\newtheorem{notation_basic}[thm]{Notation}

\newtheorem{terminology_basic}[thm]{Terminology}

\theoremstyle{definition}
\newtheorem{definition_basic}[thm]{Definition}	
\newtheorem{example_basic}[thm]{Example}

\newenvironment{rmk}
{\pushQED{\qed}
	
	\begin{remark_basic}}
	{\popQED
\end{remark_basic}}

\newtheorem{asspt}{\textbf{Assumption}}

\newtheorem{quest}{Question}

\newenvironment{definition}
{\pushQED{\qed}
	
	\begin{definition_basic}}
	{\popQED
\end{definition_basic}}

\newenvironment{example}
{\pushQED{\qed}
	
	\begin{example_basic}}
	{\popQED
\end{example_basic}}

\DeclareMathAlphabet{\mathcal}{OMS}{cmsy}{m}{n}

\newcommand{\R}{\mathbb{R}}
\newcommand{\PP}{\mathbb{P}}
\newcommand{\C}{\mathbb{C}}

\newcommand{\A}{\mathbb{A}}

\newcommand{\Z}{\mathbb{Z}}
\newcommand{\N}{\mathbb{N}}

\newcommand{\Zmod}{\mathbb{Z}/2\mathbb{Z}}

\newcommand{\CP}{\mathbb{C}\mathbb{P}^n}
\newcommand{\RP}{\mathbb{R}\mathbb{P}^n}
\newcommand{\HP}{\mathbb{H}\mathbb{P}^n}
\newcommand{\CPtwo}{\mathbb{C}\mathbb{P}^2}
\newcommand{\CProduct}{\mathbb{C}\mathbb{P}^2\times \mathbb{C}\mathbb{P}^2}
\newcommand{\CPtimes}{\mathbb{C}\mathbb{P}^n\times \mathbb{C}\mathbb{P}^n}
\newcommand{\HPtwo}{\mathbb{H}\mathbb{P}^2}
\newcommand{\RPtwo}{\mathbb{R}\mathbb{P}^2}

\newcommand{\APn}{\mathbb{A}\mathbb{P}^n}
\newcommand{\APtwo}{\mathbb{A}\mathbb{P}^2}
\newcommand{\RPtre}{\mathbb{R}\mathbb{P}^3}
\newcommand{\CPtre}{\mathbb{C}\mathbb{P}^3}
\newcommand{\CPone}{\mathbb{C}\mathbb{P}^1}

\renewcommand{\L}{\mathcal{L}}

\newcommand{\J}{\mathcal{J}}

\newcommand{\Q}{\mathbf{Q}}

\renewcommand{\O}{\mathcal{O}}

\newcommand{\eps}{\varepsilon}

\newcommand{\OO}{\mathrm{O}}

\newcommand{\Symp}{\mathrm{Symp}}
\newcommand{\Diff}{\mathrm{Diff}}

\newcommand{\Crit}{\mathrm{Crit}}
\newcommand{\HF}{\operatorname{HF}}
\newcommand{\CF}{\operatorname{CF}}

\renewcommand{\leq}{\leqslant}
\renewcommand{\geq}{\geqslant}

\newcommand{\into}{\hookrightarrow}

\newcommand{\lra}{\longrightarrow}

\newcommand{\re}{\mathrm{Re}}
\newcommand{\im}{\mathrm{Im}}
\renewcommand{\ker}{\mathrm{ker}}

\newcommand{\Hom}{\mathrm{Hom}}

\renewcommand{\Im}{\operatorname{\mathrm{Im}}}

\makeatletter

\newcommand{\Rmnum}[1]{\expandafter\@slowromancap\romannumeral #1@}
\makeatother

\newcommand{\sheafHom}{\operatorname{\mathscr{H}\text{\kern -3pt {\calligra\large om}}\,}}

\newcommand{\sheafEnd}{\operatorname{\mathscr{E}\text{\kern -3pt {\calligra\large nd}}\,\,\,}}
\newcommand{\sheafEndmon}{\sheafEnd_{\text{\kern -3pt {\calligra mon}}\,}}

\renewcommand{\Im}{\operatorname{Im}} 
\newcommand{\ind}{\operatorname{ind}}

\newcommand*{\Jhat}[1]{#1\kern-0.37em\hat{\phantom{#1}}}

\title{Model projective twists and generalised lantern relations}
\author{Brunella Charlotte Torricelli}
\address{Brunella Torricelli, Centre for Mathematical Sciences, University of Cambridge, 
	CB3 0WB, UK}
\email{bct27@cam.ac.uk}
\begin{document}
	\maketitle

\begin{abstract}
	We use Picard-Lefschetz theory to introduce a new local model for the planar projective twists $\tau_{\APtwo} \in \Symp_{ct}(T^*\APtwo), \ \A \in \{ \R, \C \}$. In each case, we construct an exact Lefschetz fibration $\pi\colon T^*\APtwo\to \C$ with three singular fibres, and define a compactly supported symplectomorphism $\varphi \in \Symp_{ct}(T^*\APtwo)$ on the total space. Given two disjoint Lefschetz thimbles $\Delta_{\alpha},\Delta_{\beta} \subset T^*\APtwo$, we compute the Floer cohomology groups $\HF(\varphi^k(\Delta_{\alpha}), \Delta_{\beta};\Zmod)$, for $k\in \N$, and verify (partially for $\CPtwo$) that $\varphi$ is indeed isotopic to (a power of) the projective twist in its local model.
	
	The constructions we present are governed by  \emph{generalised lantern relations}, which provide an isotopy between the global monodromy of a Lefschetz fibration and a fibred twist along an $S^1$-fibred coisotropic submanifold of the smooth fibre. We also use these relations to generate non-exact fillings for the contact manifolds $(ST^*\CPtwo, \xi_{std}), (ST^*\RPtre,\xi_{std})$, and study two classes of monotone Lagrangian submanifolds of $(T^*\CPtwo, d\lambda_{\CPtwo})$.
\end{abstract}

\section{Introduction}

\subsection{Dehn and projective twists}

Given a Lagrangian sphere $L$ embedded in a symplectic manifold $(M^{2n},\omega)$, there is a compactly supported symplectomorphism of $M$ called the Dehn twist along $L$, and denoted by $\tau_L \in \Symp_{ct}(M)$. 
The construction and well-definedness of these symplectomorphisms rest on the fact that the tangent bundle to a sphere has a periodic geodesic flow, and all geodesics have a common period.

A celebrated feature that has made the study of Dehn twists particularly productive is the fact that Dehn twists occur as (local) monodromies of Lefschetz fibrations (\cite{arnold, seideles}).
This important property has opened up new avenues to the study of symplectic mapping class groups using tools from (symplectic) Picard-Lefschetz theory, a strategy pioneered by the early work of Seidel (\cite{seideles, seidelbook}). The monodromy of an exact Lefschetz fibration can never be isotopic to the identity in the symplectic mapping class group (\cite{bgzfilling}, \cite{bct}), so Dehn twists obtained that way are indeed an important source of non-trivial symplectomorphisms of exact symplectic manifolds. 

In \cite{seidelgraded}, Seidel introduced a wider class of symplectomorphisms defined from Lagrangian submanifolds with periodic geodesic flow; among these are notably spheres (in which case the resulting symplectomorphisms are squared Dehn twists) and projective spaces. This paper focuses on the latter class of symplectomorphisms, that we call \emph{projective twists} (the appellation \emph{Dehn} will be exclusively associated to twists along spheres). The definition of projective twist requires the existence of a Lagrangian embedding of a projective space in the ambient manifold, which might impose strong topological restrictions. Moreover, many results given by Picard-Lefschetz theory that do hold for Dehn twists do not apply to their projective counterparts. Nevertheless, there are a series of results that suggest that projective twists do have interesting properties worth of consideration; Evans (\cite{evans}), Harris (\cite{harristwist}), Mak-Wu (\cite{makwu1}), and the author (\cite{bct}). 

In this paper, we introduce a new local model for the real and complex projective twists in dimension two.

Although not formally discussed in the existing literature, there are ``standard'' local models for real and complex projective twists that can be deduced from the work of \cite{perutzmbl,chikoe,wwfibred}. The existence of such models hinges upon the non-trivial fact that projective twists can be modelled as local monodromies of (a specific type of) \emph{Morse--Bott--Lefschetz} fibrations (abbreviated MBL), a class of fibrations which admits singularities produced by Morse--Bott degenerations. This phenomenon, which echoes the symplectic Picard-Lefschetz formula, is due to the fact that projective twists concurrently belong to a class of symplectomorphisms called \emph{fibred twists} (introduced in full generality by Perutz), which are the symplectomorphisms that encode the symplectic monodromy maps of Morse--Bott degenerations (\cite{perutzmbl}).

Given a symplectic manifold $(M, \omega)$ and a coisotropic manifold $V\subset M$, endowed with the structure of a sphere bundle $V\to B$ over a symplectic manifold $(B,\omega_B)$, a fibred Dehn twist is defined as a compactly supported symplectomorphism of $M$ that acts as a Dehn twist on the fibres of $V$. Projective twists are a special type of $S^1$-fibred twists (we show this in Section \ref{stdmodel}), and as a result they admit a presentation as local monodromy of Morse--Bott--Lefschetz fibrations. From this point of view, the existing local model for projective twists emerges from the MBL fibration with fibres symplectomorphic to disc cotangent bundles of $\APtwo$ and whose monodromy is defined by the fibred twist along the unit cotangent bundle $ST^*\APtwo$, viewed as an $S^1$-bundle $ST^*\APtwo\to B$ over a symplectic manifold $(B, \omega_B)$.

\subsection{Picard-Lefschetz theory}

Interestingly, studying the presentation of projective twists as fibred Dehn twists and the local model resulting from this point of view also delivers the first ingredients to build our alternative models: Lefschetz fibrations on the cotangent bundles of $\RPtwo$ and $\CPtwo$ respectively.

Lefschetz fibrations were first thought of as fibration-like structures arising from families of hyperplane sections (with nodal singularities) on complex projective varieties (\cite{lamotke, lefschetz}), so called \emph{Lefschetz pencils} (see Definition \ref{defpencil}). 
A smooth projective variety $(X, \omega_X)$ with an ample line bundle $\L\to X$ can be viewed as a Kähler manifold with integral symplectic class $c_1(\L)=[\omega] \in H^2(X;\Z)$. Then, the zero locus of a section $s_{\infty}\in H^0(\L)$ vanishing transversely is a smooth hypersurface $\Sigma_{\infty}=s_{\infty}^{-1}(0)$ Poincaré dual to $[\omega]$. Given another generic section $s_0 \in H^0(\L)$, the family $\{ \lambda s_{0}+ \mu s_{\infty}=0 \}_{[\lambda:\mu]\in \CPone}$ parametrised by $\CPone$ determines a Lefschetz pencil; $\Sigma_{\infty} \subset X$ is a smooth fibre in this family. The above data determines the geometry of the pencil, in particular the following facts will be central to this paper.

\begin{enumerate}[label=(\roman*)]
	\item On one hand, the (closure of) the complement $X\setminus \Sigma_{\infty}$ is a Stein domain with plurisubharmonic function $f\colon X\setminus \Sigma_{\infty} \to \R$ given by $f:=-log\|s_{\infty}\|^2$, on which the ratio $s_0/s_{\infty}$ induces the structure of an exact Lefschetz fibration (see \cite[(19b)]{seidelbook}).  \label{observationi}
	
	\item On the other hand, if $f\colon X\setminus \Sigma_{\infty}\to \R$ is Morse, its unstable manifolds are isotropic. Let $\Xi\subset X$ be the union of these isotropic submanifolds.
	A theorem of Biran (\cite{biran1}) proves that $X \setminus \Xi$ is symplectomorphic to an open disc normal bundle to $\Sigma_{\infty}$ (see Theorem \ref{birandecompo}).
	In the special case that $\Xi$ embeds as a Lagrangian submanifold, endowed with Riemannian metric with periodic geodesic flow, then Biran's decomposition can be extended as in Proposition \ref{audinprop}. In that case, $X\setminus \Sigma_{\infty}$ is symplectomorphic to an open disc cotangent bundle $\mathring{D}_{\eps}T^*\Xi$, $\eps>0$ and the symplectic hypersurface $\Sigma_{\infty}$ is the space of geodesics in $\Xi$. \label{observationii}
\end{enumerate}

In Section \ref{examplescotangent}, we find, for $\A \in \{ \R, \C \}$, the appropriate projective variety $(X, \omega_X)$ admitting a line bundle $\L\to X$, and Lagrangian embedding of $\APn$ compatible with Biran's decomposition as in \ref{observationii}.

For $\A \in \{ \R, \C \}$, this yields the decomposition $X\cong \mathring{D}_{\eps}T^*\APn \cup \Sigma_{\infty}$, which, by \ref{observationi} delivers a Lefschetz fibration on the disc cotangent bundle of $\APn$ and by \ref{observationii} delivers information on the structure of the coisotropic $ST^*\APn$; an $S^1$-fibre bundle $ST^*\APn \to \Sigma_{\infty}$ over the space of geodesics $\Sigma_{\infty}$.

In particular, for $\A=\R$ (Example \ref{examplerp}), we have $X=(\CP, \omega_{FS})$ with the ample line bundle $\L:=\O_{\CP}(2) \to \CP$, a generic section of which defines a quadric $\Sigma_{\infty} \subset \CP$. This yields the decomposition

 \begin{align}\label{decompositions1}
\CP \cong \mathring{D}_{\eps}T^*\RP \cup \Sigma_{\infty}, \ \eps>0.
\end{align} 

And in the case $\A=\C$ (Example \ref{examplecp}), $X=(\CPtimes, \omega_{FS}\oplus \omega_{FS})$ with the ample line bundle $\mathcal{L}:=\mathcal{O}_{\CPtimes}(1,1):=pr_1^*(\mathcal{O}_{\CP}(1))\otimes pr_2^*(\mathcal{O}_{\CP}(1))$, where $pr_i\colon \CPtimes \to \CP$ is the projection to the $i$-th factor, $i=1,2$. A generic section of $\L$ defines a $(1,1)$-divisor $\Sigma_{\infty}\subset \CPtimes$ and a special version of Biran's decomposition gives

\begin{align}\label{decompositions2}
\CPtimes \cong \mathring{D}_{\eps} T^*\CP \cup \Sigma_{\infty}, \ \eps>0.
\end{align}

Combining the decompositions \eqref{decompositions1} and \eqref{decompositions2} with the study of a pencil of quadrics on $(\CPtwo, \omega_{FS})$ and a pencil of $(1,1)$-divisors on $(\CProduct, \omega_{FS}\oplus\omega_{FS})$ respectively, we obtain the following Lefschetz fibrations.

\begin{prop}\label{fibrationprop}
	\begin{enumerate}
		\item There is an exact Lefschetz fibration $\pi\colon E_{\RPtwo} \to \C$ with smooth fibre exact symplectomorphic to a sphere with four boundary components, three singular fibres and such that the completion of the total space $E_{\RPtwo}$ is exact symplectomorphic to $(T^*\RPtwo,d\lambda_{T^*\RPtwo})$ (Section \ref{pencil1}). 
		\item There is an exact Lefschetz fibration $\pi\colon E_{\CPtwo} \to \C$ with smooth fibre exact symplectomorphic to a plumbing of disc cotangent bundles $D_{\eps}(T^*S^3\#_{S^1}T^*S^3)$ along a circle and three singular fibres (see for example \cite{bct} for a definition of clean plumbing), such that the completion of the total space $E_{\CPtwo}$ is exact symplectomorphic to $(T^*\CPtwo,d\lambda_{T^*\CPtwo})$ (Section \ref{lfcp2}). 
	\end{enumerate}
\end{prop}

\subsection{Results}
As they arise from Lefschetz pencils of hypersurfaces on projective varieties, the two Lefschetz fibrations mentioned above have the property that for any vanishing cycle in the smooth fibre, $V\subset M$, the global monodromy of the fibration $\phi \in \Symp_{ct}(M)$ satisfies \begin{align}\label{preservecycles}
\phi(V)\simeq V. \end{align}
This is a consequence of a phenomenon we call \emph{general lantern relation}, a relation in the symplectic mapping class group of the smooth fibre $(M,\omega)$ (see Lemma \ref{generalantern}). This special feature of the monodromy enables us, in the two cases, to apply a construction of \cite{seideliterated}, which takes as input a symplectomorphism of the fibre of a Lefschetz fibration and outputs a compactly supported symplectomorphism of the total space (Section \ref{construction}). The construction involves lifting a Dehn twist along an annulus in the base $\C$ to a symplectomorphism $\Phi_{2\pi}$ of the total space, that we adjust to a compactly supported symplectomorphism $\varphi \in \Symp_{ct}(T^*\APtwo)$.

This construction represents a fundamentally new method of building a symplectomorphism. In Section \ref{chaptermodels} we illustrate this method, so far rather unexplored, for the fibrations above to obtain non-trivial elements of the symplectic mapping class groups $\pi_0(\Symp_{ct}(T^*\APtwo))$ that we compare with the standard planar projective twists (Sections \ref{rp2} and \ref{cp2twist} respectively).

To do that, we measure the Floer theoretical action of the symplectomorphism $\varphi$ on a Lefschetz thimble by computing the Floer cohomology groups $\HF(\varphi^k(\Delta_{\alpha}), \Delta_{\beta};\Zmod)$ for elements $\Delta_{\alpha}, \Delta_{\beta}$ in a distinguished basis of Lefschetz thimbles. In the real case, we can combine these computations with the knowledge of the mapping class group $\pi_0(\Symp_{ct}(T^*\RPtwo))$ (\cite{evans}) to obtain:

\begin{theorem}[{{Corollary \ref{corollaryreal}}}]
	The symplectomorphism $\varphi \in \Symp_{ct}(D^*\RPtwo)$ is isotopic to a power of the projective twist $\tau_{\RPtwo}^k$, $k\in \Z^*$.
\end{theorem}

In the complex case, very little is known about $\pi_0(\Symp_{ct}(T^*\CPtwo))$; we know for example about the existence of projective twists with non-standard framings (\cite{bct}), but only in higher dimensions. As a consequence, we can only deliver a partial result based on our Floer cohomological computations.

\begin{theorem}[{{Corollary \ref{corollaryorder}}}]
	The symplectomorphism $\varphi \in \Symp_{ct}(D^*\CPtwo)$ is of (symplectic) infinite order.
\end{theorem}

We can only conjecture that $\varphi$ is isotopic to the complex planar projective twist. However, there is strong evidence that suggests our constructions correspond to the projective twists in dimension two. Namely, in the two cases, a comparison with existing computations found in the literature confirms that the symplectomorphisms we build should indeed be projective twists.

Additional results are obtained by utilising generalised lantern relations to operate so called ``monodromy substitutions''. We show that these operations can produce strong fillings of unit cotangent bundles (Sections \ref{fillingrp3}, \ref{fillingcp2}).

As another byproduct of our results we obtain information on monotone Lagrangian submanifolds $L\cong S^1\times S^3 \subset (T^*\CPtwo, d\lambda_{\CPtwo})$.

\begin{theorem}[{{Lemma \ref{lemmaisotopyclass}}}]
	There are two distinguished Lagrangian isotopy classes of Lagrangian submanifolds $T_1, T_2 \cong S^1\times S^3 \subset T^*\CPtwo$.
\end{theorem}

\subsection{Organisation of the paper}

Section \ref{mblsection} reviews Lefschetz and Morse--Bott--Lefschetz fibrations in relation to Dehn and fibred twists respectively. In Section \ref{stdmodel} we outline the standard local model for projective twists that can be readily obtained when considering those symplectomorphisms as fibred $S^1$-twists.
In Section \ref{chapterfibrations} we look at the geometry and topology of exact Lefschetz fibrations that are obtained from Lefschetz pencils, providing the background for the two main examples of Section \ref{pencil1} (Lefschetz fibration on $T^*\RPtwo$) and Section \ref{lfcp2} (Lefschetz fibration on $T^*\CPtwo$). In Section \ref{chaptermodels} we give the general construction of a compactly supported symplectomorphism on the total space of an exact Lefschetz fibration satisfying certain properties, recipe that we apply to the two fibrations above in Section \ref{chaptermodels2} to obtain new local models for projective twists.
\subsection{Ackowledgements}

First and foremost, I would like to express my deepest gratitude to my supervisor Ivan Smith, whose support and guidance have been an invaluable component of the development and completion of this project. I am thankful to him for providing knowledgeable advice and helping me getting out of theoretical impasses along the way. Moreover, I am grateful to him for taking the time to comment an ample number of versions of this paper, and offering suggestions for improvement.

It was Ailsa Keating who first suggested that the model $\RPtwo$ twist could be constructed using a Lefschetz fibration as we present here. I am very thankful to her for sharing the ideas which provided the momentum to engender this project, and for helpful conversations in the initial stage of the research.

I am extremely thankful to Jack Smith for offering impeccable Floer theoretical advice, over numerous conversations, which helped me articulate the technicalities of Section \ref{hfcomputationsection}. I would like to thank Jonny Evans for sharing his insights about the Morse--Bott--Lefschetz fibration of Section \ref{mbldiscussion} and helpful discussions.

Many thanks to Denis Auroux, who helped clarify the structure of the bifibration of Section \ref{bifibration} via the map \eqref{net} and Agustín Moreno, who provided helpful comments on Section \ref{monodromysubst}, and mentioned \cite{oba} to me. I also want to thank Jeff Hicks and Navid Nabijou for helpful discussions.

I am enormously indebted to the keystone of DPMMS, Vivien Gruar, whose generous help and support (both administrative and personal) played an essential role in fostering the optimal work environment to pursue my research.

Last but not least, I wish to express my sincere gratitude to my family and friends, who accompanied me in these years of meandering explorations towards intellectual autonomy.\\

be water.

\section{(Mose-Bott-) Lefschetz fibrations}\label{mblsection}

In this section we review how Dehn and fibred twists arise as local monodromies of Lefschetz and Morse--Bott--Lefschetz fibrations respectively.

\subsection{Dehn twists and Lefschetz fibrations}\label{lfconventions}

\begin{definition}\label{defliouville}
	A Liouville manifold of finite type is an exact symplectic manifold $(W, \omega=d\lambda_W)$, where $\lambda_W \in \Omega^1(W)$ is called the Liouville form, such that there exists a proper function $h_W\colon W \rightarrow \lbrack 0, \infty )$ and $c_0>0$ with the following property. For all $ x \in ( c_0, \infty )$ and $y\in h^{-1}(x)$, the vector field $Z_W$ dual to $\lambda_W$, called the Liouville vector field, satisfies $dh_W(Z_W)(y)>0$. 
	
	For a regular value $c$ of $h_W$, a closed sublevel set $M:=h_W^{-1}([0,c])$ of a Liouville manifold $(W, d\lambda_W)$ is a compact symplectic manifold with contact type boundary $(\Sigma:=h_W^{-1}(c), \lambda_W|_{\Sigma})$, and it is called a Liouville domain. 
\end{definition}

\begin{definition}\label{exactsymplecto}
	An exact symplectomorphism between two Liouville manifolds $(W_1, d\lambda_1), (W_2, d\lambda_2)$ is a diffeomorphism $\psi\colon W_1 \to W_2$ satisfying $\psi^*\lambda_2-\lambda_1=df$, for a compactly supported function $f\colon W_1\to \R$. 
\end{definition}

\begin{definition}\label{defacs}
	Let now $(M,d\lambda)$ be a Liouville domain with contact boundary $(\Sigma=\partial M,\alpha=\lambda|_{\Sigma})$. The negative Liouville flow identifies a collar neighbourhood $C(\Sigma)$ of the boundary with $(-\eps, 0]\times \partial M$, such that $\lambda|_{C(\Sigma)}=e^t\alpha$. An almost complex structure $J$ of \emph{contact type near the boundary} is one that satisfies $de^t\circ J= -\lambda$.
\end{definition}

\begin{definition}
	Given a Liouville domain $(M,d\lambda)$ as above, we can use the identification of the collar neighbourhood $C(\Sigma)$ to glue an infinite cone and define the symplectic completion of $M$: \begin{align}
	(W, \omega_W):=(M \cup [ 0, \infty) \times \partial M, d(e^t \alpha)), 
	\end{align}
	where $t$ is the coordinate on $(0, \infty)$, such that the Liouville flow extends to $Z_W$ with $Z_W|_{[ 0, \infty) \times \partial M}=\partial_t$.
	
	An almost complex structure $J$ of contact type extends to an almost complex structure $J_W$ on the completion satisfying \begin{itemize}
		\item $J_W(\frac{\partial }{\partial t})=R_{\alpha}$, where $R_{\alpha}$ is the Reeb vector field associated to $\alpha$,
		\item $J_W$ is invariant under translations in the $t$-direction,
		\item $J_W|_{M}=J$.
	\end{itemize}
	This kind of almost complex structure will be called cylindrical. 
\end{definition}

We will only consider Liouville manifolds that are complete (i.e with complete Liouville vector field) and of finite type, which we can identify as the union of a Liouville domain with a cylindrical non-compact end, equipped with an almost complex structure cylindrical at infinity.

Let $(E^{2n+2}, \Omega_E,\lambda_E)$ be a Liouville manifold, with a compatible almost complex structure $J_E$, and consider the complex plane with its standard symplectic form and complex structure $j_{\C}$. Let $\pi\colon E \to \C$ be a map with finitely many critical points, which are all non-degenerate, and contained in a compact set of $E$.
Denote by $\Crit(\pi):=\{ x \in E, \ D_x\pi=0 \}$ the set of critical points, and by $\Crit v(\pi):=\pi(\Crit(\pi))$ the set of critical values.

Unless punctually specified in exceptional situations, by a Lefschetz fibration we will always mean the following.

\begin{definition}\label{deflf}
	A Lefschetz fibration on (the Liouville manifold) $E$ is a $(J_E,j_{\C})$-holomorphic map $\pi$, i.e $D\pi \circ J_E=j_{\C} \circ D\pi$, with the above properties and the following additional features.
	\begin{enumerate}
		
		\item For all $x\in E\setminus \Crit(\pi)$, $\ker(D_x\pi)\subset T_xE$ is symplectic.
		\item Every smooth fibre is symplectomorphic to the completion of a Liouville domain $(M, d\lambda_M)$.
		
		\item\label{trivialitycondition} There is an open neighbourhood $U^h \subset E$ such that $\pi\colon E \setminus U^h \lra \C$ is proper and $\pi|_{U^h}$ can be trivialised via an isomorphism $f\colon U^h \cong \C \times ([0, \infty) \times \partial M)$ such that
		\begin{align}\label{productlike}
		f^*(\lambda_E)=\lambda_{\C}+e^t\lambda_M.
		\end{align}\end{enumerate}
\end{definition}
For more details about how this fibration is modelled outside of a neighbourhood of the critical points, see \cite[(2.1)]{maysei}.

By the first point above, there is a symplectic splitting
\begin{align}\label{splitting} 
T_xE=\ker(D_x \pi)\oplus T_xE^h,
\end{align}

where $T_xE^h$ is the symplectic complement of $\ker(D_x \pi)$ with respect to $\Omega_E$. The decomposition in \eqref{splitting} defines a canonical connection over $\C\setminus \Crit v(\pi)$. By the triviality condition \ref{trivialitycondition}., for every path $\gamma\colon [0,1]\to D\setminus \Crit(\pi)$, there are well-defined parallel transport maps $h_{\gamma}\colon E_{\gamma(0)}\to E_{\gamma(1)}$ which yield symplectomorphisms between smooth fibres.

\begin{definition}\label{compatiblej}
	A pair $(J_E,j_{\C})$ is said to be \emph{compatible} with $\pi$ if the following holds.
	\begin{itemize}
		\item $D\pi\circ J_E=j_{\C}\circ D\pi$.
		\item There is a local Kähler structure $J_0$ such that $J_E=J_0$ in a neighbourhood of $\Crit(\pi)$.
		\item On the neighbourhood $U^h$, $J_E$ is a product, $f^*(J_E)=(j_{\C},J^{vv})$, where $J^{vv}$ is a cylindrical almost complex structure compatible with $d(e^t\lambda_M)$.
		\item $\Omega_E(\cdot, J_E \cdot)$ is symmetric and positive definite.
	\end{itemize}
\end{definition}

\begin{rmk}\label{nongenericrmk}
	This choice of almost complex structure is not generic. However, the space of compatible almost complex structures on the total space of an exact Lefschetz fibration is contractible (\cite[Section 2.1]{seideles}), and the moduli spaces we will consider still meet the usual regularity requirements (\cite[Section 2.2]{seideles}).
\end{rmk}

For a Lefschetz fibration on a Liouville manifold $(E,\Omega_E)$, the proper fibration obtained as $E\setminus U^h \lra \C$, for an open neighbourhood $U^h\subset E$ as above, carries the same symplectic information as $\pi$ with the difference that its fibres are Liouville domains, and as result the total space admits a non-trivial \emph{horizontal boundary}, given by the union of the boundaries of all fibres. 

In most of the paper we will employ this latter type of Lefschetz fibration (for notational simplicity), and unless specified, an \emph{exact} Lefschetz fibration will denote a fibration obtained in this way. 

Let now $\pi\colon E\to \C$ be an exact Lefschetz fibration, with smooth fibre given by the Liouville domain $(M,d\lambda)$. By the triviality assumption of Definition \ref{deflf}, there is a neighbourhood of $U^{\partial} \subset E$ of the horizontal boundary $\partial^h E$ that is isomorphic to an open neighbourhood of the trivial bundle $\C \times \partial M$: \begin{align}\label{horizontalboundary}
U^{\partial}\cong\C \times M^{out} \subset \C \times M
\end{align}
where $M^{out}\subset M$ is a collar neighbourhood of $\partial M$. The isomorphism is compatible with the Liouville forms and the almost complex structures.

Let $\pi\colon E \to \C$ be a Lefschetz fibration with exact compact fibre $(M, \omega)$ and distinct critical values $\Crit v(\pi)=\{ w_0, \dots , w_m \} \subset D_R$, where $D_R\subset \C$ is a disc of radius $R$. Fix a base point $z_{*} \in \R$, such that $z_*\gg R$, and an identification $\pi^{-1}(z_*)\cong M$. In what follows we will frequently use the fact that via parallel transport, any fibre $\pi^{-1}(z)$ for $z\in \C$ with $\re(z)>R$ can be symplectically identified with the smooth fixed fibre $M$ via parallel transport.

\begin{definition}
	\begin{enumerate}
		
		\item A vanishing path associated to a critical value $w_i\in \Crit v(\pi)$ is a properly embedded path $\gamma_i \colon \R^+\to \C$ with $\gamma_i^{-1}(\Crit v(\pi))=\{ 0 \}$, $\gamma_i(0)=w_i$ and $\lim_{t\to \infty} \re (\gamma(t))=\infty$ such that
		outside of a compact set containing the critical values, the image of $\gamma_i$ is a horizontal half ray at height $a_i \in \R$: \begin{align}\label{height0}
		\exists T>0 \text{ such that } \forall t>T, \  \re(\gamma_i(t))>R, \ \Im(\gamma_i(t))=a_i.
		\end{align}
		\item A distinguished basis of vanishing paths for $\pi$ is a collection of $m+1$ disjoint paths $(\gamma_0, \dots \gamma_m) \subset \C$ defined as above, with pairwise distinct heights satisfying $a_0<a_1 < \cdots < a_m$.

		\item The corresponding basis of Lefschetz thimbles is the unique set of Lagrangian discs $(\Delta_{\gamma_0}, \dots , \Delta_{\gamma_m}) \subset E$ where $\Delta_{\gamma_i}$ is defined as the set of points which under the limit $t \to 0$ of the parallel transport maps over $\gamma_i$ are mapped to the critical point in $\pi^{1}(w_i)$ (the proof of uniqueness can be found in \cite[(16b)]{seidelbook}).
		Given a general Lefschetz thimble $\L$, define its height $a(\L)$ as the value defined in \eqref{height0}. For a pair of thimbles $(\L_0, \L_1)$ set $\L_0 > \L_1$ if $a(\L_0)>a(\L_1)$.
		
		\item There is an associated basis of vanishing cycles $(V_0, \dots , V_m)$ where for all $i=0,\dots , m$, $$V_i=\partial \Delta_{\gamma_i}=\Delta_{\gamma_i}\cap M\subset M$$ (using the above identification for smooth fibres).
		Every vanishing cycle $V_i\subset M$ is an exact Lagrangian sphere which comes with an equivalence class in of diffeomorphisms $S^{n} \to V_i$ defined up to the action of O$(n+1)$ (called a \emph{framing}). This is induced by the restriction of a diffeomorphism $D^{n+1} \lra \Delta_i$ (which is canonical, see \cite[Lemma 1.14]{seideles}).
	\end{enumerate}
\end{definition}


\begin{figure}[h]
	\centering
\includegraphics[width=12cm]{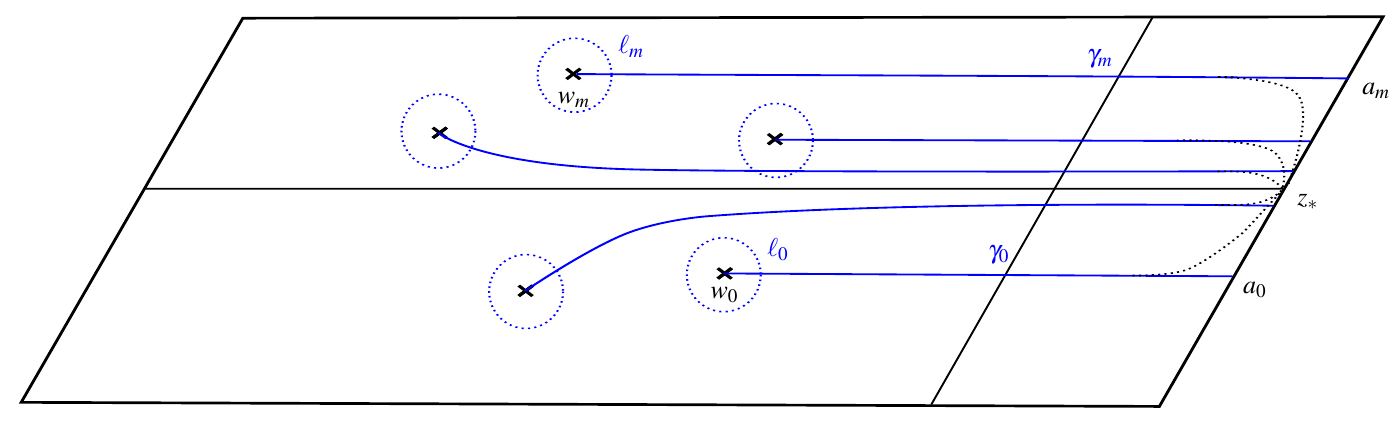}
	\caption{A distinguished basis of vanishing paths $(\gamma_0, \dots, \gamma_m)$.}\label{basis}
\end{figure}

\begin{definition}
	The global monodromy is the symplectomorphism $\phi\in \Symp_{ct}(M)$ whose Hamiltonian isotopy class is defined by the anticlockwise parallel transport map around a loop through the base point $z_*$ encircling all the critical values of the fibration. (Typically, this loop is defined as the smoothing of the concatenation of the loops centred at $z_*$ going around a single critical value as in Figure \ref{basis}).
\end{definition}
The symplectic Picard-Lefschetz theorem (\cite{arnold}) states that the global monodromy $\phi$ is isotopic to the product of the Dehn twists along the vanishing cycles $(V_0, \dots , V_m)$, \begin{align}
\phi \simeq \tau_{V_0}\cdots \tau_{V_m} \in \Symp_{ct}(M),
\end{align}
and the Hamiltonian isotopy class is independent of the choice of basis of vanishing paths.

On the other hand, given the data $( M, (V_0, \dots, V_m ) )$, there is an exact Lefschetz fibration $\pi\colon E\to \C$ with fibre $(M,\omega)$, and vanishing cycles $(V_0, \dots , V_m)\subset M$, unique up to exact symplectomorphism (\cite[(16e)]{seidelbook}).

\subsection{Fibred twists and MBL fibrations}\label{mbldiscussion}

Lefschetz fibrations can be viewed as a special case of \emph{Morse--Bott--Lefschetz} (abbreviated MBL) fibrations, a class of fibrations which allows non-isolated singularities. The monodromies of such fibrations are symplectomorphisms called \emph{fibred twists} (\cite{perutzmbl}), which naturally generalise Dehn twists.
The relevance of fibred twists in this discussion is Lemma \ref{projasfibred}, which shows that projective twists are just a special type of fibred twists, thereby providing a first local model for these symplectomorphisms (obtained as monodromy of a suitable MBL fibration).

In this section we briefly discuss MBL fibrations and ($S^1$-)fibred twists, to be able to state Lemma \ref{projasfibred}, which we'll prove in Section \ref{modelfibred}. Sections \ref{discbundlesection} and \ref{examplescotangent} are auxiliary to this proof.

Let $D$ be a disc with standard almost complex structure $j$, and $(E, \Omega_E)$ a $(2n+2)$-manifold with a closed $2$-form $\Omega_E$ and an almost complex structure $J_E$. A MBL fibration $\pi\colon E \to D$ is a proper $(J_E,j)$-holomorphic map such that $\Omega_E$ is non-degenerate on $\ker(D\pi)$, and whose singular fibres admit Morse--Bott singularities (the definition is more involved, see \cite[Definition 2.1]{perutzmbl}). This means that the critical locus $\Crit(\pi) \subset E$ can be a set of smooth closed symplectic submanifolds $(Q_i,\Omega_E|_{Q_i})$ of (real) dimension $2n-2k_i$, $k_i\geq 0$ (\cite[Definition 2.1]{perutzmbl}), and that the Hessian of $\pi$ is non-degenerate only along the normal bundle to $\Crit(\pi)$. 
For this reason, in general the two form on the total space cannot be made exact.

Every critical point in a component $Q_i \in E^{crit}$ of the singular locus has a neighbourhood with local complex coordinates $(x_0, \dots x_n)$ such that in that neighbourhood, \begin{align}
\pi(x_0, \dots x_n)=\sum_{j=0}^{k_i} x_j^2.
\end{align}

A MBL fibration is characterised by its smooth fibre $(M, \omega)$ and its ``fibred vanishing cycles'', a set of submanifolds of the smooth fibre playing the analogue role of vanishing cycles of a Lefschetz fibration. Fibred vanishing cycles are defined as the sets of points in the smooth fibre that are mapped into the critical locus via the parallel transport maps. For any component $Q_i$ of the critical locus, the associated fibred vanishing cycle is a coisotropic submanifold of $(M, \omega)$ which is an $S^{k_i}$-bundle $V_i \to Q_i$. The analogy with vanishing cycles is that the image of $V_i$ under the limit of the parallel transport into $Q_i$ collapses the spherical fibres of $V_i$ (\cite[Section 2.3]{perutzmbl}). 

The monodromy of a MBL fibration with smooth fibre $(M, \omega)$, one singularity with (compact) critical locus $Q$ and associated fibred coisotropic $p\colon V \to Q$ is given by the \emph{fibred twist} (see Definition \ref{fibretwistdef}) $\tau_V \in \Symp_{ct}(M)$ around $V$ (\cite[Theorem 2.16]{perutzmbl}). On the other hand, any fibred twist can be realised as (local) monodromy of such a MBL fibration (\cite[Section 2.4.1]{perutzmbl}). 

Instead of giving the general definition of fibred twists, we focus on a special case of interest, that of $S^1$-fibred twists along a coisotropic $V\to Q$ which has the structure of a \emph{prequantisation bundle}.\\

\begin{definition}\label{defbw}
	Let $p\colon V \to Q$ be a principal $S^1$-bundle over a closed symplectic manifold $(Q, \omega_Q)$.
	We call $p\colon V \to Q$ a \emph{prequantisation bundle} if $\omega_Q$ is an integral symplectic form (meaning $[\omega_Q]$ lies in the image of the map $H^2(Q;\Z)\to H^2(Q;\R)$), and the Euler class satisfies $e(p)=-k[\omega_Q]$, for some $k>0$. In that case, there is a connection $1$-form $\alpha\in \Omega^1(V;\R)$ (regarded as a regular real valued $1$-form) which is a contact form whose Reeb vector field generates the $S^1$-action on $V$, and $d\alpha=2k\pi p^*(\omega_Q)$ (see e.g \cite[Theorem 7.2.4]{geigesbook}). 
\end{definition}

Let $(W, \omega_W)$ be a symplectic manifold, and $V \subset W$ a coisotropic submanifold admitting the structure of a prequantisation bundle $p\colon V\to Q$ over a closed symplectic manifold $(Q,\omega_Q)$.

By the coisotropic neighbourhood theorem (see \cite[Section 2.4.1]{perutzmbl}, \cite[Section 2.3]{wwfibred}), a neighbourhood $U$ of $V\subset W$ is isomorphic to an associated symplectic bundle $(V \times_{S^1} T^*S^1, p^*(\omega_Q)+d\lambda_{T^*S^1}+d \langle \mu, \alpha \rangle)$, where the quotient is taken with respect to the $S^1$-action
\begin{align}\label{s1action}
(x, (\theta, t))\cdot \varphi= (x\cdot \varphi, (\theta- \varphi, t)), \ \varphi\in S^1, \ (\theta, t)\in T^*S^1.
\end{align}
The two-form given above is invariant under the $S^1$-action \eqref{s1action} so it descends to a well-defined two form on the quotient. 
Here $\lambda_{T^*S^1}=d(t\theta)$ is the standard Liouville form on $T^*S^1$ and $\mu\colon T^*S^1 \to \R$ is the moment map of the standard $S^1$-action on $T^*S^1$ (which is the distance function), and $\langle \cdot , \cdot \rangle$ is the Lie algebra pairing (see \eqref{nhoods} for the concrete description of this two-form).
\vspace{0.5cm}

\begin{definition}\label{fibretwistdef}
	\begin{enumerate}
		\item{\textbf{Local fibred twist.}} Let $\tau_{S^1}\colon T^*S^1\to T^*S^1$ the standard $S^1$-Dehn twist (see \cite[Section 4.b]{seidelgraded}). On the product space $V\times T^*S^1$, define a diffeomorphism 
		\begin{align}
		\bar{\tau} \colon	V \times T^*S^1 \lra V \times T^*S^1 \\
		(x, (\theta,t)) \mapsto (x, \tau_{S^1} (\theta,t))
		\end{align}
		Since $\bar{\tau}$ commutes with the diagonal $S^1$ action, it descends to a well-defined compactly supported diffeomorphism $\tau_{loc}$ of the quotient $V\times_{S^1}T^*S^1$ with the properties that
		\begin{enumerate}
			\item $\tau_{loc}$ is a symplectomorphism,
			\item $\tau_{loc}$ covers the identity on $Q$,
			\item $\tau_{loc}$ acts as a standard Dehn twist on each fibre of $V\times_{S^1}T^*S^1$
		\end{enumerate}
		see (\cite[Lemma 2.3]{perutzmbl}, \cite[Definition 2.7]{wwfibred})	
		
		\item{\textbf{Fibred twist along $V$.}} For $\eps>0$, let $\psi \colon V \times_{S^1} D_{\epsilon} T^*S^1 \to U$ be the symplectomorphism of neighbourhoods as above. Then \begin{align*}
		\tau_V:= \left\{
		\begin{array}{ll}
		\psi \circ \tau_{loc} \circ \psi^{-1}	 & \text{ on }  U\\
		Id & \text{ on } W\setminus U
		\end{array}
		\right.
		\end{align*}
		defines a compactly supported symplectomorphism in $\Symp_{ct}(W)$.
	\end{enumerate}
\end{definition}

We now describe another construction from the literature, formalised by Chiang--Ding--Koert, and we show it corresponds to the special instance of fibred twist above.

Let $p\colon (V, \alpha) \to (\Sigma, \omega_{\Sigma})$ be a prequantisation bundle. The flow $\sigma^{R_{\alpha}}_t$of the Reeb vector field $R_{\alpha}$, given its periodicity, gives rise to a loop of contactomorphisms on $V$. It is then possible to \emph{suspend} this loop of contactomorphisms to a symplectomorphism of the symplectisation of $V$ (see \cite[Proposition 3.1]{cks} for a general statement about such a construction).

\begin{definition}[{{\cite[Definition 2.7, Lemma 2.1]{chiangtwist}}}]\label{chiangfibred}
	Let $(W,\omega_W)$ be a Liouville domain such that $\partial W=V$, $\omega_W|_{V}=d\alpha$. Consider a neighbourhood of $V\subset W$, i.e a piece of symplectisation \begin{align}
	(V\times [0,1], d(e^t\alpha)).
	\end{align}
	Fix a smooth function $f\colon [0,1]\to \R$ with $f(0)=2\pi$ and $f(1)=0$, define
	\begin{equation}\label{fibredtwistbw}
	\begin{aligned}
	\tau_V' \colon V\times [0,1] \lra V \times [0,1] \\
	(x, t) \mapsto (\sigma_{f(t)}^{R_{\alpha}}(x), t)
	\end{aligned}
	\end{equation}
	and extend to $W\setminus (V\times [0,1])$ via the identity ($\tau_V'$ is the identity near the boundary of $V\times [0,1]$). Then $\tau_V'$ defines a symplectomorphism of $W$. 
\end{definition}

\begin{lemma}\label{equivalencedef}
	The construction of the fibred twist in Definition \ref{fibretwistdef} is equivalent to a suspension of a loop of contactomorphisms in the sense of Definition \ref{chiangfibred}.

\end{lemma}
\proof 
As a first step, we explain how the neighbourhoods of $V\subset W$ used to define $\tau_V$ and $\tau_V'$ are symplectomorphic (abstractly, this is clear by the coisotropic neighbourhood theorem).
In Definition \ref{fibretwistdef}, we used an associated symplectic $T^*S^1$-bundle
\begin{align}\label{nhoods}
(V \times_{S^1} D_{\eps}T^*S^1, p^*(\omega_Q)+d\lambda_{T^*S^1}+d \langle \mu, \alpha \rangle=\frac{d \alpha}{2\pi}+d(t\theta)+d(t\alpha)).
\end{align}
For $\eps=\frac{1}{2}$ we can use the identification $D_{\eps}T^*S^1\cong S^1\times [0,1]$ and get $V\times_{S^1} D_{\eps}T^*S^1\cong V\times [0,1]$, so that the Liouville form $\lambda_{T^*S^1}$ gets identified with $d(t\alpha)$ and \eqref{nhoods} is symplectomorphic to a collar neighbourhood $ (V \times [0,1] , d((2t+\frac{1}{2\pi})\alpha))$.
The correct symplectic form is obtained after a change of variable.

It is clear from the definition that $\tau_V'$ covers the identity on $\Sigma$ and after the above identification, it acts as a standard $S^1$-Dehn twist on the fibres of $V \times [0,1]$ (note that we can linearly interpolate between the choices of cut-off functions in the two definitions).

Therefore, the symplectomorphism $\tau_{V'}$ agrees with the fibred twist we have defined in Section \ref{mbldiscussion}.
\qed

\section{Models of projective twists via MBL degenerations}\label{stdmodel}

The projective twists $\tau_{\APn}, \A\in \{ \R, \C \}$ are examples of $S^1$-fibred twists along the unit cotangent bundle of $(T^*\APn, d\lambda_{T^*\APn})$.

\begin{lemma}\label{projasfibred}
	Let $\A \in \{ \R, \C \}$. The projective twist $\tau_{\APn} \in \Symp_{ct}(T^*\APn)$ is isotopic to the $S^1$-fibred Dehn twist (of Definitions \ref{fibretwistdef} and \ref{chiangfibred}) along the coisotropic 
	\begin{align}\label{fibredcoiso}
	(V, \alpha)=(ST^*\APn, \lambda_{T^*\APn}|_{ST^*\APn}) \to (\Sigma, \omega_{\Sigma})
	\end{align} 
	where \begin{enumerate}
		\item For $\A=\R$, $\Sigma:=\{ [z_0: \dots : z_n] \in \CP, \ \sum_{i=0}^{n}z_i^2=0\} \subset \CP$ is the $(n-1)$-complex quadric with the restriction of the Fubini--Study form $\omega_{FS}$;
		\item For $\A=\C$, $\Sigma:= \{ ([x_0: \dots :x_n], [y_0:\dots y_n])\in \CPtimes, \   \sum_{i=0}^{n} x_iy_i=0 \} \subset \CPtimes$ with the restriction of $\omega_{FS}\oplus\omega_{FS}$
		
	\end{enumerate}
	
	In particular, a projective twist is the local monodromy of a MBL fibration with fibred vanishing cycle $V\to \Sigma$.

\end{lemma}

The proof will be given in Section \ref{modelfibred}, after the excursus of Sections \ref{discbundlesection} and \ref{examplescotangent}, which are aimed at clarifying what the coisotropic submanifolds $(V, \alpha)\to (\Sigma, \omega_{\Sigma})$ are in the two examples. In that way, we will be able to visualise what the local model of projective twists is when considered as fibred twists.

\subsection{Symplectic decompositions of complex projective varieties}\label{discbundlesection}

To understand the structure of the fibred coisotropics \eqref{fibredcoiso}, it is necessary to provide an alternative description of the cotangent bundle of the real and complex projective spaces.

This section contains the necessary theoretical tools to do so, by viewing these (disc) cotangent bundles as the complement of a divisor $\Sigma \subset X$ in a complex projective variety $X$ (the specific examples are carried out in Section \ref{examplescotangent}). The main reference for this section is \cite[Section 2.1]{biran1}.

Let $X^{2n}$ be a smooth complex projective variety with an ample line bundle $\mathcal{L}\to X$; we view $X$ as a Kähler manifold with an integral symplectic form $\omega_X$ with the property that $c_1(\mathcal{L})=[\omega_X]$.

Assume there is a section $s\in H^0(\mathcal{L})$ such that $\Sigma:=s^{-1}(0) \subset X$ is a smooth projective subvariety (transversely cut out). Then $(\Sigma, \omega_{\Sigma}=\omega_X|_{\Sigma})$ is a Kähler submanifold of $X$, with fundamental class $[\Sigma]\in H_{2n-2}(X;\Z)$ Poincaré dual to $[\omega_X]$.

The (closure of the) complement $X \setminus \Sigma$ is a Stein domain, with the plurisubharmonic function $\varphi:=-\frac{1}{4\pi}log(\rho)= -\frac{1}{4\pi}log (\| s\|^2)$ and Liouville vector field $Z=grad(\varphi)$. The restriction $\omega:=\omega_X|_{X \setminus \Sigma}$ is exact and coincides with $-dd^{\C} \varphi$. By the symplectic neighbourhood theorem (\cite[Theorem 3.4.10]{mcduffsal}), a neighbourhood of $\Sigma \subset X$ is symplectomorphic to the \emph{symplectic normal bundle} to $\Sigma$, that we define below in \eqref{symplecticnormal}. On an algebro-geometric level, the normal bundle to $\Sigma \subset X$ is modelled on $\mathcal{L}|_{\Sigma}$. Define $P:= \{ v \in \mathcal{L}^*|_{\Sigma},  \| v\|=1  \}$, for a hermitian metric on the dual bundle induced by one on $\mathcal{L}$ (we make this choice of orientation in order for $P$ to be viewed as a convex boundary, positively oriented with respect to the Liouville vector field $Z$ above).

The bundle $p\colon P \to \Sigma$ is a principal $S^1$-bundle with Euler class $e(p)=-[\omega_X]|_{\Sigma}=-[\omega_{\Sigma}]$. Let $\alpha_P \in \Omega^1(P;\R)$ be a connection $1$-form as in Definition \ref{defbw}, with curvature $d\alpha_P=2\pi p^*(\omega_{\Sigma})$ and such that it is a contact form whose associated Reeb vector field $R_{\alpha_P}$ generates the $S^1$-action on $P$. It gives $P$ the structure of prequantisation bundle \begin{align}\label{prequantum}
p\colon (P, \alpha_P) \lra (\Sigma, \omega_{\Sigma}).
\end{align}

Consider the following disc bundle associated to $P$, \begin{align}
(P\times_{S^1} D(1), p^*(\omega_{\Sigma})+d(r^2d\theta)-d(r^2\alpha_P) )
\end{align}
where $D(1)\subset \C$ is the open unit disc with radial coordinate $r$, and $\theta \in S^1$ acts on $(x,z)\in P\times D(1)$ as \begin{align}
(x,z)\cdot \theta=(x\cdot \theta, e^{2\pi i \theta}z). 
\end{align}

The form $p^*(\omega_{\Sigma})+d(r^2d\theta)-d(r^2\alpha_P)= d((\frac{1}{2\pi}-r^2)(\alpha_P-d\theta))=:d\lambda_{\nu}  \in \Omega^1(P\times D(1))$ is invariant under the $S^1$-action 
and therefore induces a well-defined one form on the quotient $P\times_{S^1} D(1)$. Note that this bundle is positive, with first Chern class given by $+[\omega_{\Sigma}]=-e(p)$.

By the symplectic neighbourhood theorem, there is a symplectic embedding
\begin{align}\label{symplecticnormal}
\nu\colon (P\times_{S^1} D(1), d\lambda_{\nu}) \lra (X, \omega_X),
\end{align}
and the symplectic disc normal bundle is the image $N_{\Sigma/X}:=\nu(P\times_{S^1} D(1))$. 
The zero section of the associated bundle on the left hand side is mapped to $\Sigma$, and on its complement there is a vector field $X_{\nu}:=-\frac{1-r^2}{2r}\partial_r$, such that $\nu_*(X_{\nu})$ is a multiple of the Liouville vector field $Z$ on $X\setminus \Sigma$ (see \cite[2.2]{oba}).

\begin{thm}[{{\cite[Theorem 1.A]{biran1}}}]\label{birandecompo}
	Let $(X, \omega_X)$ be a smooth projective variety as above, with an ample line bundle $\mathcal{L}\to X$ defining a submanifold $\Sigma \subset X$ as the smooth zero locus of a non-trivial section $s\in H^0(\mathcal{L})$. Assume that the function $\rho\colon X\to \R$, $\rho(x)=\| s(x) \|^2$ (with respect to the same hermitian metric as above) is Morse on $X\setminus \Sigma$ and let \begin{align}\label{isotropiccw}
	\Xi := \bigcup_{p \in \Crit(\rho)\setminus \Sigma} W_{p}^u(\rho, g),
	\end{align}
	where $W_p^u(\rho, g)$ is the unstable manifold of a critical point $p\in \Crit(\rho)$ with respect to the (positive) gradient flow, defined via the Kähler metric $g$ on $X$ compatible with $\omega_X$. Then $\Xi \subset X$ is an isotropic CW-complex whose complement, the open dense subset $(X\setminus \Xi, \omega_X)$, is symplectomorphic to the symplectic open disc normal bundle $N_{\Sigma/X}$.
\end{thm}

In Secton \ref{chaptermodels}, we will use this decomposition in the situation in which the isotropic CW complex is a smooth Lagrangian submanifold of the given Kähler manifold.

\begin{lemma}\label{inversedecompolemma}
	Assume $X, \omega, \rho$ and $\varphi:=-\frac{1}{4\pi}log(\rho)$ as above. Suppose that the isotropic CW complex $\Xi\subset X$ defined in \eqref{isotropiccw} is a smooth Lagrangian submanifold which coincides with $\Crit(\varphi)=\Crit(\rho)\setminus \Sigma$.
	Let $U_{\Sigma/X} \supset \Sigma$ be an open disc subbundle of the normal bundle $N_{\Sigma/X}$ given in \eqref{symplecticnormal}. Then the complement $X \setminus U_{\Sigma/X}$ is symplectomorphic to a disc subbundle (of some radius $\eps >0$) of the cotangent bundle $(T^*\Xi, d\lambda_{T^*\Xi})$:
	\begin{align}\label{inversedecomposition}
	X \cong D_{\eps}T^*\Xi \cup U_{\Sigma/X}.
	\end{align}
\end{lemma}
\proof 

This lemma follows as a corollary to Biran's decomposition theorem, see \cite[Section 7]{biran1}. 

Let $U_{\Xi/X}\supset \Xi$ be an open neighbourhood of $\Xi \subset X$, and $U_{\Sigma/X}\supset \Sigma$ an open neighbourhood of $\Sigma$, the subbundle of $N_{\Sigma/X}$ as in the statement.

Consider the Liouville vector field $Z= grad(-\frac{1}{4\pi} log \|s \|^2)=grad(\varphi)$ on $X\setminus \Sigma$ with $\mathcal{L}_Z\omega_X=-dd^{\C}\varphi=\omega_X|_{X\setminus \Sigma}$ for the plurisubharmonic function $\varphi$. For any $t>0$, the flow $Z^t$ of $Z$ acts by conformal symplectomorphisms.
Then, there is a symplectic cobordism between the boundaries of the closures of $U_{\Sigma/X}$ and $U_{\Xi/X}$, and since $\Xi$ is the only critical set of $\varphi$ in $X\setminus \Sigma$, the flow lines of the vector field $Z$ define a trivialisation of the cobordism (so the two boundaries are contactomorphic). By the Weinstein's neighbourhood theorem, $U_{\Xi/X}$ is a disc cotangent bundle, and there is $C>0$ and $t\in [-C,C]$ such that $X\cong Z^t(U_{\Xi/X})\cup_{ST^*\Xi} U_{\Sigma/X}$. \qed

\subsection{Examples: disc cotangent bundles of projective spaces}\label{examplescotangent}

Consider a compact Riemannian manifold $(\Xi, g)$ with periodic geodesic flow, and a circle bundle of radius $r\in (0,1)$ inside $(T^*\Xi,d\lambda_{T^*\Xi})$, denoted by $ST^*\Xi= \{ v \in T^*\Xi; \ \|v \|_g =r  \}$. The latter is a contact hypersurface of the cotangent bundle, with contact form $\lambda_{T^*\Xi}|_{ST^*\Xi}$. The Hamiltonian vector field of the distance function $\mu\colon T^*\Xi \to \R$ $\mu(v)=\| v \|_g$ is (a multiple of) the Reeb vector field of the contact form on $ST^*\Xi$, whose flow coincides with the (co-) geodesic flow on $\Xi$.
The latter is periodic, so the Reeb vector field generates a free $S^1$-action on $(ST^*\Xi, \lambda_{T^*\Xi}|_{ST^*\Xi})$, which is a regular level set of $\mu$. So the unit cotangent bundle is a principal $S^1$-bundle over the space of geodesics.

In this case, a result stronger than Lemma \ref{inversedecompolemma} holds.

\begin{prop}[{{\cite[Proposition 4.1, 4.3]{audin2}}}]\label{audinprop}
	Let $\Xi \in \{ S^n, \RP, \CP, \HP \}$. Then there is projective variety $(X, \omega_X)$ where $\Xi \subset X$ is embedded as Lagrangian submanifold, and such that there is an open disc bundle $\mathring{D}_{\eps}T^*\Xi\subset T^*\Xi$ with a symplectic decomposition\begin{align}\label{specialbiran}
	X \cong \mathring{D}_{\eps}T^*\Xi \cup \Sigma
	\end{align}
	where $\Sigma:=ST^*\Xi / S^1$ is the space of (unparametrised, oriented) geodesics of $\Xi$.
\end{prop}

The decomposition \eqref{specialbiran} is a special instance of \eqref{inversedecomposition}; there is an ample line bundle $\L\to X$ and a holomorphic section $s\in H^0(\L)$ such that $s^{-1}(0)\cong \Sigma$. The restriction $\omega_X|_{X \setminus \Sigma}$ is a multiple of the standard symplectic form on the cotangent bundle $d\lambda_{T^*\Xi}$, and the coisotropic $(ST^*\Xi, \lambda_{T^*\Xi}|_{ST^*\Xi}) \to (\Sigma, \omega|_{\Sigma})$ is a prequantisation bundle contactomorphic to $(P,\alpha_P)$ in \eqref{prequantum}.

We illustrate two examples in which we can use a decomposition of the type \eqref{inversedecomposition} to present a disc subbundle of $T^*\APn$, $\A \in \{ \R, \C \}$ as the complement of an ample divisor in a projective variety.

\begin{example}[{{\cite[3.1.2]{biran1}}}]\label{examplerp}\label{rp2decompo}
	Let $X=(\CP, \omega_{FS})$ with the ample line bundle $\L:=\O_{\CP}(2) \to \CP$. A generic section $s\in H^0(\mathcal{L})$ is given by a homogeneous polynomial of degree two in the homogeneous coordinates $[z_0: \cdots : z_n]$ of $\CP$. Consider the section $s:=\sum_{i=0}^n z_i^2$. The divisor arising as the zero locus of $s$ is the smooth quadric $\Sigma:=s^{-1}(0)= \{ z_0^2 + \cdots + z_n^2=0 \}\cong Q^{n-1} \subset \CP$, which is a Kähler submanifold with the integral symplectic structure $\omega_{\Sigma}=\omega_{FS}|_{\Sigma}$.
	The smooth function \begin{align*}
	\rho\colon \CP \lra \R, \	\rho([z_0: \cdots : z_n])=\frac{|\sum_{i=0}^{n}z_i^2|^2}{(\sum_{i=0}^{n}|z_i|^2)^2},
	\end{align*}
	is a reparametrisation of $\| s \|^2$, and is Morse--Bott on $X\setminus \Sigma$ with unique critical locus given by\begin{align*}
	\Crit(\rho|_{X\setminus \Sigma})=\Crit(\rho)\setminus \Sigma= \{ [z_0: \cdots : z_n] \in \CP| z_j \in \R \text{ for all } 0\leq j \leq n \} \simeq \RP \subset \CP,
	\end{align*}
	which is embedded as a Lagrangian in $(\CP,\omega_{FS})$.
	The positive gradient vector field points out of $\Crit(\rho)\setminus \Sigma$, so the union of unstable manifolds coincides with the critical locus of $\rho|_{X\setminus \Sigma}$: $\Xi=\Crit(\rho)\setminus \Sigma\simeq \RP \subset \CP$.
	
	Therefore, by Lemma \ref{inversedecompolemma} and Proposition \ref{audinprop}, the complement $X \setminus \Sigma$ is symplectomorphic to an open disc cotangent bundle of $\Xi$: \begin{align}\label{decomporpn}
	\CP \cong \mathring{D}_{\eps}T^*\RP \cup \Sigma, \ \eps>0.
	\end{align}

	The unit cotangent bundle is a prequantisation bundle $(ST^*\RP, d\lambda_{T^*\RP}) \lra (\Sigma, \omega_{\Sigma})$ over the quadric $\Sigma\subset \CP$.		
\end{example}

\begin{example}[{{\cite[3.3]{biran1}}}]\label{examplecp}
	Let $X=(\CPtimes, \omega_{FS}\oplus \omega_{FS})$ with the ample line bundle $\mathcal{L}:=\mathcal{O}_{\CPtimes}(1,1):=pr_1^*(\mathcal{O}_{\CP}(1))\otimes pr_2^*(\mathcal{O}_{\CP}(1))$, where $pr_i\colon \CPtimes \to \CP$ is the projection to the $i$-th factor, $i=1,2$. A generic section $s\in H^0(\mathcal{L})$ is a homogeneous polynomial of degree $(1,1)$, in the homogeneous coordinates $(\underline{x}, \underline{y})=([x_0: \cdots : x_n], [y_0: \cdots : y_n])$ of $\CPtimes$. Let $s:=\sum_{i=0}^n x_iy_i$; the divisor obtained as the zero set of $s$ is given as \begin{align}\label{smoothfibrecpn}
	\Sigma:= \left\lbrace \sum_{i=0}^{n}x_iy_i=0 \right\rbrace \subset \CPtimes.
	\end{align}
	
	The function \begin{align*}
	\rho\colon \CPtimes \to \R, \ \rho(\underline{x},\underline{y})=\frac{|\sum_{i=0}^{n}x_iy_i|^2}{(\sum_{i=0}^{n}|x_i|^2)^2(\sum_{i=0}^{n}|y_i|^2)^2}
	\end{align*}
	is a smooth reparametrisation of $\| s \| ^2$, and is Morse--Bott on $X\setminus \Sigma$ with critical locus given by 
	\begin{align*}
	\Crit(\rho|_{X\setminus \Sigma})=\Crit(\rho)\setminus \Sigma= \{ ([x_0: \cdots : x_n], [\overline{x_0}:\cdots : \overline{x_n}]), \ [x_0 \cdots : x_n] \in \CP \} \simeq \CP \subset \CPtimes,
	\end{align*}
	which is embedded as a Lagrangian in $(\CPtimes, \omega_{FS}\oplus \omega_{FS})$.
	Following the same reasoning as in the previous example, the union of unstable manifolds $\Xi$ is given by $\Xi = \Crit(\rho)\setminus \Sigma\simeq \CP$. 
	
	Using Lemma \ref{inversedecompolemma} and Proposition \ref{audinprop}, we obtain a decomposition 
	\begin{align}\label{decompocpn}
	\CPtimes \cong \mathring{D}_{\eps} T^*\CP \cup \Sigma, \ \eps>0.
	\end{align}
	
	As we explain below in Remark \ref{rmkflag}, the unit cotangent bundle is a prequantisation bundle $(ST^*\CP, d\lambda_{T^*\CP})\lra (\Sigma, \omega_{\Sigma})$ over the projectivisation of the holomorphic cotangent bundle; $\Sigma\cong \PP(T_{hol}^*\CP)$, which, for $n=2$, is isomorphic to the Flag $3$-fold.
\end{example}
\begin{rmk}\label{rmkflag}
	
	Consider
	\begin{align*}
	\Sigma ' := \{ (\ell, H), \ \ell \text{ is a line in } \C^{n+1}, \ H \text{ is a hyperplane in } \C^{n+1} \text{with } \ell \subset H \} \subset \CP \times (\CP)^*.
	\end{align*}
	This admits a projection \begin{equation}\label{projectivisedprojection}
	\begin{aligned}
	\Sigma '& \lra \CP \\
	(\ell,H)&\longmapsto \ell.
	\end{aligned} \end{equation} 
	that realises $\Sigma '$ as the projectivised cotangent bundle $\Sigma '\cong \PP(T_{hol}^*\CP)$ (see for example \cite[Proposition 10.22]{eisenbudharris}). [The fibre over $\ell\in \PP(\C^{n+1})$ of \eqref{projectivisedprojection} is given by $\{ w \in \PP((\C^{n+1})^*), w(\ell)=0 \}=\PP(\C^{n+1}/ \ell)$.]

	Consider the diffeomorphism
	\begin{equation}
	\begin{aligned}
	f\colon &\CP\to (\CP)^* \\
	&\ell \mapsto \ell^{\perp}
	\end{aligned}
	\end{equation}
	sending a line $\ell \in \CP$ to its orthogonal complement (with respect to a fixed hermitian metric), i.e the hyperplane whose normal vector is $\ell$. Then $\CPtimes$ can be identified with $\CP\times (\CP)^*$ using $f$. Under this identification, we can identify $\Sigma$ from \eqref{smoothfibrecpn} with $\Sigma'$.
\end{rmk}

In the two dimensional case, $\Sigma\subset \CProduct$ is a Flag $3$-fold embedded in $\CProduct$. Write the Flag 3-fold (the manifold of complete flags in $\C^3$) as
\begin{align}\label{flag3fold}
Fl_3:=\{ (\ell,H), \ \ell \subset \C^3 \text{ is line through 0 }, H\subset \C^3 \text{ is a plane such that } \ell \subset H \},
\end{align} 
then lines $\ell$ can be parametrised by the first $\CPtwo$ factor, and planes $H$ by the second $\CPtwo$ factor so that in coordinates we obtain the identification $Fl_3\cong \left \{ \sum_{i=0}^{2}x_iy_i=0 \right \} = \Sigma$.

Consider the Flag $3$-fold as the set of pairs $(\ell, H)$ as in \eqref{flag3fold}. Then the projection \eqref{projectivisedprojection} gives $Fl_3 \cong \PP(T_{hol}^*\CPtwo)$.

\subsection{The projective twist as a fibred twist: proof of Lemma \ref{projasfibred}}\label{modelfibred}

Let $\A \in \{ \R, \C \}$. We use the decompositions of the previous section to prove that the projective twists $\tau_{\APn} \in \Symp_{ct}(T^*\APn)$, defined via the geodesic flow (\cite{seidelgraded}), can also be viewed as $S^1$-fibred twists along the unit cotangent bundle $ST^*\APn \subset T^*\APn$.

\proof[Proof of Lemma \ref{projasfibred}] 
The second part of the statement about the shape of the coisotropic manifolds is clear by Examples \ref{rp2decompo} and \ref{examplecp}.

As we have shown (in Lemma \ref{equivalencedef}) the equivalence of Definitions \ref{fibretwistdef} and \ref{chiangfibred}, it suffices to prove that the projective twist coincides with the fibred twist of Definition \ref{chiangfibred}.

Let $j \colon ST^*\APn \times [0,1] \into T^*\APn$ be the inclusion with $j^*(\omega)=d(e^t\alpha)$. Then the projective twist, defined as in \cite[4.b]{seidelgraded}, pulls back under the inclusion $j$ to a symplectomorphism of the shape of \eqref{fibredtwistbw} (see \cite[4.a.]{seidelgraded}), which acts on the symplectisation as the Reeb flow in one component and the identity in the other.

The result is a symplectomorphism that is obtained as suspension of a loop of contactomorphisms induced by the Reeb flow on level-set hypersurfaces; this is exactly as $\tau_V'$ was defined in \eqref{fibredtwistbw}. Since we proved that $\tau_V'$ coincided with the fibred twist $\tau_V$ of Definition \ref{fibretwistdef}, this completes the proof.

\endproof

\section{Lefschetz fibrations from Lefschetz pencils}\label{chapterfibrations}

In Section \ref{examplescotangent}, we have used Biran's symplectic decomposition (\cite{biran1}) to understand the complement of a smooth hypersurface of two complex projective varieties $(\CPtwo,\omega_{FS})$ and $(\CProduct, \omega_{FS}\oplus\omega_{FS})$. In Sections \ref{pencil1} and \ref{lfcp2}, we will look at these projective varieties as carrying the structure of a Lefschetz pencil inducing Lefschetz fibrations $\pi \colon T^*\APtwo \to \C$, $\A \in \{ \R, \C \}$ on cotangent bundles of projective spaces. These fibrations will be a key component in the construction of the model planar projective twists of Section \ref{chaptermodels}. 

Before looking at the specific examples we have in mind, in the current section we review some general facts of symplectic Picard--Lefschetz theory.

\subsection{Geometry of Lefschetz pencils}\label{lffrompencil}\label{pencilbasic}
The Lefschetz fibrations we build in this section are obtained from pencils of hyperplane sections. For this reason, we include a short review of classical material about Lefschetz pencils (\cite{lamotke}) and the properties, on a symplectic level, of Lefschetz fibrations arising from this context (\cite{aurouxmon,seidelbook,oba2}).

Let $(X,\omega_X)$ be an $(n+1)$-dimensional smooth complex projective variety, and let $\mathcal{L}\to X$ be a holomorphic ample line bundle. Let $s_0, s_{\infty} \in H^0(\L)$ be linearly independent holomorphic sections of $\L$. Then \begin{align}\label{pencil0}
\Sigma_{[\lambda:\mu]}:=	\{ \lambda s_{0}+ \mu s_{\infty}=0 \}_{[\lambda:\mu]\in \CPone} \subset X
\end{align} 
defines a family of (projective) hypersurfaces. There is a rational function \begin{align}\label{pencil}
X \dashrightarrow \CPone, \ z\mapsto [s_0(z):s_{\infty}(z)]
\end{align}
defined away from the base locus \begin{align}
B:=\bigcap_{[\lambda:\mu] \in \mathbb{C}\mathbb{P}^1} \Sigma_{[\lambda:\mu]}=\{s_0=s_{\infty}=0 \}.
\end{align}
\begin{definition}\label{defpencil}
	Let $\Sigma_{\infty}:=\Sigma_{[0:1]}=s_{\infty}^{-1}(0)$ be smooth.
	The family \eqref{pencil0} is called a (algebraic) Lefschetz pencil if \begin{enumerate}
		\item The base locus $B$ is a smooth submanifold of $X$ (of real codimension $4$).
		\item The map $p_X \colon X\setminus B \to \CPone$ induced by \eqref{pencil} admits a (finite) set of non-degenerate critical points $\Crit(p_X)$, and is a submersion away from $\Crit(p_X)$. Denote the set of critical values by $\Crit v(p_X)$.
	\end{enumerate}
\end{definition}
The two above properties imply that:
\begin{enumerate}[resume]
	\item In a neighbourhood of $p \in B$ there are local coordinates $(z_0, \dots , z_n)$ such that $B=\{ z_0=z_1=0 \}$ and $p_X(z_0, \dots , z_n)=\frac{z_1}{z_0}$.
	\item In a neighbourhood of $p \in \Crit(p_X)$ there are local coordinates $(z_0, \dots , z_n)$ such that $p_X(z_0, \dots , z_n)=z_0^2 + \cdots +z_n^2$.
\end{enumerate}

\begin{rmk}
	Pencils exist on any projective variety (see for example \cite[Lemma 2.10]{voisin}), and Donaldson (\cite{donaldson1, donaldson2}) proved the existence of a Lefschetz pencil on any closed integral symplectic manifold $(X, \omega)$.
\end{rmk}

\subsection{From pencil to fibration}\label{monodromy}
A natural way of turning a Lefschetz pencil into a Lefschetz fibration is to perform a blow-up of $X$ at the base locus $B$. Let $\widetilde{X}:=Bl_{B} X=\{ (z,y) \in \CPone \times X : \ s_0(y)=zs_{\infty}(y) \}.$
The projection map $\widetilde{p}\colon \widetilde{X}\to \CPone$ defines a Lefschetz fibration over $\CPone$ with closed fibres, each of which contains a copy of the base locus embedded as a smooth hypersurface. 

Note that the geometry of the blown-up space yields a well known formula relating the Euler characteristics of the various components of a Lefschetz pencil, which will be helpful for later computations.

\begin{lemma}\label{lemmauler}
	Let $X \subset \C\PP^N$ be a smooth projective variety of complex dimension $n+1$. Let $(\Sigma_{[\lambda:\mu]})_{[\lambda:\mu] \in \CPone}$ be a Lefschetz pencil on $X$ with base locus $B$, generic fibre $\Sigma$ and $r+1$ critical points, then the Euler characteristics of $X$, $\Sigma$ and $B$ are related as follows: \begin{align}\label{eulerpencil}
	\chi(X)=2 \chi(\Sigma)- \chi(B)+ (-1)^{n+1} (r+1). 
	\end{align}\qed
\end{lemma}

Another strategy to produce a Lefschetz fibration from the projection map of a Lefschetz pencil is to remove the smooth generic fibre $\Sigma_{\infty}=s_{\infty}^{-1}(0)$ from the family \eqref{pencil0}. 
In the process, the base locus is removed, and the restriction $p_X|_ {X \setminus \Sigma_{\infty}}$ gives rise to a well-defined map
\begin{align}\label{prelf}
p \colon  X \setminus \Sigma_{\infty} \lra \C.
\end{align}

The total space $X\setminus \Sigma_{\infty}$ is affine, and as noted in Section \ref{discbundlesection}, $\varphi:=-\frac{1}{4\pi}log(\| s \|^2)$ is a plurisubharmonic function on it. Similarly, a generic smooth fibre of \eqref{prelf} is given by the complement of the base locus $\Sigma \setminus B$, which is also affine.

The ``affine'' fibration $p$ can be turned into a an Lefschetz fibration with exact compact fibres as follows. Consider the space $E:=X\setminus (\Sigma_{\infty}\cup U_{B/X})$, where $U_{B/X}$ is a disc subbundle of the normal bundle $N_{B/X}$. 
The latter can be obtained as in \eqref{symplecticnormal}, with the difference that in this case, the normal bundle to $B\subset X$ is modelled on $\L\oplus \L$. Equivalently, $E$ is obtained as the complement $X\setminus U_{\Sigma_{\infty}/X}$ of an open subbundle $U_{\Sigma_{\infty}/X}\subset N_{\Sigma_{\infty}/X}$ of the normal bundle to the smooth fibre in $X$.
The smooth fibres of the restriction $p|_{E}$ are Stein domains $M:=\Sigma \setminus U_{B/\Sigma}$, where $U_{B/\Sigma}\subset N_{B/\Sigma}$ is an open subbundle of the normal bundle to $B\subset \Sigma$. 

\begin{prop}[{{\cite[Proposition 2.6]{oba2}}}]\label{frompenciltolf}
	Let $(X, \omega_X)$ be a smooth projective variety with a Lefschetz pencil induced by a holomorphic line bundle $\L\to X$, with base locus $B$ and generic smooth fibre $\Sigma_{\infty}.$ Then, for a disc bundle $U_{B/X}\supset B$ as above there is a closed two form $\Omega$ on $E:=X\setminus  U_{\Sigma_{\infty}/X}$ such that 
	\begin{enumerate}
		\item $\pi=p|_{E} \colon (E, \Omega)\lra \C$ is an exact Lefschetz fibration with smooth fibre $(M, \omega)$,
		\item For every $z\in \C$, $\omega_X|_{\pi^{-1}(z)}=\Omega|_{\pi^{-1}(z)}$,
		\item $\Omega$ coincides with $\omega_X$ outside a collar neighbourhood of the horizontal boundary of $E$.
	\end{enumerate}
	The total space $(E,\Omega)$ can be (symplectically) completed to a Stein manifold.
\end{prop}

\begin{rmk}\label{coisotropicboundary}
	By construction, the boundary of the smooth fibre $(M,\omega)$ has the structure of a prequantisation bundle $P:=\{ \nu \in \L^*|_B, \ \| \nu \| =1 \}$ over the base locus $B$ (see Section \ref{discbundlesection}).
\end{rmk}

The symplectic Picard--Lefschetz theorem expresses the global monodromy of a Lefschetz fibration as the product of Dehn twists along the Lagrangian vanishing cycles of the singularities. The global monodromy of an exact Lefschetz fibration induced by a Lefschetz pencil, like that of Theorem \ref{frompenciltolf}, has an additional property, that we call the \emph{generalised lantern relation} (the isotopy \eqref{generalisedlantern}).

\begin{lemma}[{{\cite[Section 2.4]{aurouxmon}, \cite[Section 3]{gompf} \cite[Theorem 2.7]{oba2}}}]\label{lantern1}\label{generalantern}
	Consider a Lefschetz pencil on the projective variety $(X, \omega_X)$, with generic smooth fibre $\Sigma$ and base locus $B$. Let $\pi\colon E\to \C$ be the exact Lefschetz fibration induced by the pencil as in Theorem \ref{frompenciltolf}, with smooth fibre $(M,\omega)$, and global monodromy $\phi \in \Symp_{ct}(M)$. Then there is a symplectic isotopy \begin{align}\label{generalisedlantern}
	\phi \simeq \tau_V,
	\end{align}
	where the map $\tau_{V} \in \Symp_{ct}(M)$ is a \emph{fibred} twist along the $S^1$-fibred coisotropic submanifold $V=\partial M\to B$, obtained as in Remark \ref{coisotropicboundary}.
	
\end{lemma}

\subsection{The monodromy as a graded symplectomorphism}\label{gradingshift}\label{gradingmonod}\label{lthimbles}\label{vcycles}

In some cases, it is possibe to define a notion of $\Z$ grading, for Lagrangian submanifolds, Floer cohomology groups and symplectomorphisms. In this section we study how the global monodromy of a Lefschetz fibration that is built as in Section \ref{monodromy} acts (as a graded symplectomorphism) on vanishing cycles (as graded Lagrangians). To illustrate this, we follow \cite[2]{seidelgraded} and \cite[11]{seidelbook}.

For now, let $(M^{2n}, \omega)$ be a general symplectic manifold satifying $2c_1(M)=0$.
Since $2c_1(M)=2c_1(TM)=c_1((\bigwedge^n_{\C}TM)^{\otimes 2})$, $2c_1(M)=0$ implies the bundle $\mathcal{K}_M^2:=(\bigwedge^n TM)^{\otimes -2}$ is trivial so it admits a nowhere vanishing section, i.e a quadratic complex volume form $\eta_M$. 

Let $LGr(n)$ be the Grassmannian of Lagrangian planes in $(\R^{2n}, \omega_{std})$, and $LGr(TM)$ the bundle of Lagrangian Grassmannians in the tangent bundle of $M$, which is a $LGr(n)$-bundle over $M$.

\begin{definition}
	The squared phase map $\alpha_M\colon LGr(TM)\to S^1$, is defined as
	\begin{align*}
	\alpha_M \colon LGr(TM) \to S^1, \	\alpha_M(V)=\frac{\eta_M((v_1 \wedge \cdots \wedge v_n)^2)}{|\eta_M((v_1 \wedge \cdots \wedge v_n)^2)|}=\arg(\eta_M|_V),
	\end{align*}
	where $v_1, \dots , v_n$ is any basis of $V$.

\end{definition}

Let $L\subset M$ be a closed connected Lagrangian submanifold of $M$. There is a section $s_L\colon L \to LGr(TM)$ sending $p\in L$ to the corresponding tangent plane $T_pL \in LGr(T_pM)$. Consider the composition $\alpha_L=\alpha_M \circ s_L\colon L\to S^1$.

\begin{definition}\label{maslovnumbers}
	The Maslov class of $L$, denoted by $\mu_L\in \Hom(\pi_1(L), \Z)=H^1(L;\Z)$ is the homotopy class of the map $[\alpha_L]\in [L,S^1]=H^1(L;\Z)$. In other words, $\mu_L$ is the pullback of the angle form $d\theta \in H^1(S^1;\Z)$ under the map $\alpha_L$.
	
	Given a pseudoholomorphic (with respect to some given compatible almost complex structure) disc $u\colon (D, \partial D) \to (M,L)$, its Maslov index is given by $\mu(u)=\int_{\gamma} \mu_L$ where $\gamma$ is the loop in $LGr(n)$ given by the (symplectic) trivialisation of $(u|_{\partial D})^*(TL)$.
	
	The minimal Maslov index of $L$ is defined as $$N_L:=\min \left\lbrace \mu([u])>0, \ u\colon (D, \partial D)\to (M,L) \text{ pseudoholomorphic} \right\rbrace.$$
\end{definition}

\begin{definition}\label{gradingdef}
	A graded Lagrangian is a pair $(L, \tilde{\alpha}_L)$ consisting of a Lagrangian submanifold $L \subset M$ with a lift of $\alpha_L$, i.e a map $\tilde{\alpha}_L \colon L \to \R $	such that $\exp(2\pi i \tilde{\alpha}_L(x))=\alpha_L(TL_x)$. If $\mu_L=0$, such a lift exists.
\end{definition}

We shall focus on grading shifts that graded exact Lagrangian submanifolds undergo under a symplectomorphism.

\begin{definition}\label{gradedsymp}
	A graded symplectomorphism is a pair $(\phi, \alpha_{\phi})$ consisting of a symplectomorphism $\phi \in \Symp(M)$ and a map $\alpha_{\phi} \colon LGr(TM)\to \R$ satisfying $\exp(2\pi i \alpha_{\phi}(V))=\frac{\alpha_M(D\phi(V))}{\alpha_M(V)}$.
\end{definition}

An important example of graded symplectomorphism is the \emph{shift} operator (\cite[(11k)]{seidelbook}), which can be considered as the pair $(id, -k)$, $k\in \Z$. 

On a graded Lagrangian $(L, \tilde{\alpha}_L)$, the shift sends $\tilde{\alpha}_L$ to $\tilde{\alpha}_L-k$. We will simply denote the action of shift on the Lagrangian as $L[k]$.

\begin{definition}
	Given a pair $(L_0,L_1)$ of graded exact closed Lagrangian submanifolds of $M$, the degree of an intersection point $x \in L_0\cap L_1$ is defined as \begin{align}
	\deg(x)=\tilde{L}_1(x)- \tilde{L}_2(x),
	\end{align} where $\tilde{L}_i(x):=\tilde{\alpha}_{L_i}(x)$.
\end{definition}

The Floer cohomology group $\HF(L_0,L_1)$ can be endowed with a grading (\cite[2.f]{seidelgraded}), and we have the following formula for shifts (\cite[Section 2f]{seidelgraded}) \begin{align*}
\HF^*(L_0[j],L_1)=\HF^{*-j}(L_0,L_1)=\HF^*(L_0,L_1[-j]), \ j \in \Z. \\
\end{align*}

For the remainder of the section, assume that $\pi\colon E^{2n+2} \to \C$ is an exact Lefschetz fibration with smooth fibre $(M,\omega)$, that is induced by a Lefschetz pencil on a projective variety $(X,\omega_X)$ as in Proposition \ref{frompenciltolf}. Let $\L\to X$ be the ample line bundle that generates the pencil. The generalised lantern relation (Lemma \ref{lantern1}) ensures that the (ungraded) global monodromy of $\pi$ preserves the vanishing cycles. We introduce an extra assumption on $X$ that yields a graded version of this statement (Lemma \ref{lemmamonodromy}).

\begin{asspt}\label{anticanonical1}
	The canonical bundle of $X$, $\mathcal{K}_{X}$, satisfies:
	\begin{align}\label{anticanonicalpole}
	\mathcal{K}_{X}^2 \cong \mathcal{L}^{\otimes -d}, \ d\in \Z^{>0}
	\end{align} 
\end{asspt}

Below we explain how, under Assumption \eqref{anticanonical1}, the symplectic parallel transport maps and in particular the global monodromy can be treated as graded symplectomorphisms. 

Let $s_0, s_{\infty} \in H^0(\L)$ be linearly independent sections of $\L$. The preimage of $s_{\infty}^{\otimes -d}$ (which is a section of $\mathcal{L}^{\otimes -d}$) under \eqref{anticanonicalpole} is the square of a $(n+1)$-form, which has a zero ($d<0$) or a pole ($d>0$) along $\Sigma_{\infty}=s_{\infty}^{-1}(0)$. That gives a trivialisation of $\mathcal{K}_X^2$, the bundle of quadratic forms, over the complement $W=X\setminus \Sigma_{\infty}$, where it restricts to a quadratic form $\eta_X$, which induces a quadratic volume form on the total space $E$ of the Lefschetz fibration.
There is an associated relative quadratic form $\eta_X /dz^2$ (here $dz$ is the volume form on $\C$), which is a section of $\mathcal{K}^2_{X/\C}=\pi^*(\mathcal{K}_{\C}^{-2}) \otimes \mathcal{K}_X^2$, on the complement of $\Sigma_{\infty}$. This defines a quadratic complex volume form on the smooth fibre $(M, \omega)$ of $E$.

\begin{lemma}[{{\cite[(19b)]{seidelbook}}}]\label{lemmamonodromy}
	The global monodromy $\phi \in \Symp_{ct}(M)$ of the exact Lefschetz fibration $\pi\colon E\to \C$ above (obtained from a Lefschetz pencil on a projective variety $X$ satisfying Assumption \eqref{anticanonical1}) acts as a shift on every graded closed exact Lagrangian vanishing cycle $V \subset M$. In other words, it acts as the identity, equipped with a constant grading, and in the notation of the discussion above we have \begin{align}
	\phi(V)= V[4-d].
	\end{align}
\end{lemma}
\proof 
By Lemma \ref{lantern1}, the global monodromy $\phi\in \Symp_{ct}(M)$ of a fibration as in the statement is isotopic to a fibred twist whose support is disjoint from the vanishing cycles. Therefore, $\phi$ preserves setwise the vanishing cycles, and will act on each of them as a shift. This shift depends on the pole order of $\eta_X/dz^2$. Under the assumption \eqref{anticanonical1}, the relative quadratic form $\eta_X/dz^2$ does not extend smoothly over infinity---where the section $s_{\infty}$ has a pole.

By \eqref{anticanonicalpole}, $\eta_X$ has a pole of order $d$ at infinity, while a standard argument shows that $dz^2$ has a pole of order $4$. Therefore, $\eta':=z^{4-d}(\eta_X/dz^2)$ is holomorphic, so that the phase function of $\eta'$ is preserved under parallel transport around a large circle through $z_*$ (that can be homotoped to be away at infinity). This implies that, on the other hand, the phase function of $\eta_X/dz^2=z^{d-4}\eta'$ must shift by $d-4$ when doing parallel transport 'around infinity' (since the phase function of $\eta'$ must remain constant).
Then, by the convention of \cite[(11k)]{seidelbook}, the grading of the graded Lagrangian vanishing cycles $V_i$ is shifted as $\phi(V_i)= V_i[-(d-4)]=V_i[4-d]$. 
\endproof

\section{Lefschetz fibration on $T^*\RPtwo$}\label{pencil1}
From a pencil of quadrics on $(\CPtwo, \omega_{FS})$, we define an exact Lefschetz fibration $\pi\colon E_{\RPtwo} \to \C$, with smooth fibre a $2$-sphere with four punctures, three singular fibres, and whose total space admits a symplectic completion to a Stein manifold exact symplectomorphic to $(T^*\RPtwo, d\lambda_{T^*\RPtwo})$.

This is a well-known case study in symplectic toplogy (see for example \cite[Section 3.1]{aurouxmon}).

Recall Example \ref{examplerp}. Consider a Lefschetz pencil $C_{[\lambda:\mu]}:=\{ \lambda s_{0}+\mu s_{\infty}=0 \}_{[\lambda:\mu] \in \CPone}$ on $\CPtwo$, generated by sections of $\mathcal{O}_{\CPtwo}(2)$ denoted by $s_0$ and $s_{\infty}$, which are homogeneous polynomials in degree two in the coordinates $[z_0:z_1:z_2]$ on $\CPtwo$.
In a generic family of conics, the (complex) curves intersect in a common set of four points (two degree two complex curves generically intersect in four points), the base locus $B=\{ s_0=0=s_{\infty} \}$. Every fibre contains these four points, and there are three degenerate cases, in which the fibre is made up of two curves (each containing two of the base points) intersecting at a point. These are the three singular fibres, any of which has one singularity (at the aforementioned intersection point) which is at most an ordinary double point. 
This might be best illustrated with an example; take $s_0:=z_0^2-z_2^2$, $s_{\infty}:=z_0^2+z_1^2+z_2^2$. Then $B=\{ [1:i\sqrt{2}:1], [1:-i\sqrt{2}:1], [1:i\sqrt{2}:-1], [1:-i\sqrt{2}:-1] \}$, and there are three singular fibres at $[\lambda:\mu]=\{ [0:1], [1:1], [1:-1] \}$.

One can check that each singular conic is composed by two lines in $\CPtwo$ containing two of the base points each; three such configurations can arise (the number of critical fibres can also be checked with the formula \eqref{eulerpencil}). For a basis of vanishing paths, the three associated vanishing cycles are the classes of circles which collapse to the intersection of the two lines in each singular conic. Hence, each vanishing cycle in a smooth fibre is represented by a loop encircling two points of the base locus (see Figure \ref{fibre}).

\begin{figure}[htb]
	\centering
	\def\svgwidth{160pt}
	\input{cycles11.pdf_tex}
	\caption{A smooth fibre of the pencil, and the three configurations of vanishing cycles.}	\label{fibre}
\end{figure}

Let $[z_0:z_1:z_2]$ be homogeneous coordinates on $\CPtwo$. There is a rational map \begin{align*}
\CPtwo \dashrightarrow  \CPone, \ \underline{z}=[z_0:z_1:z_2] \mapsto [s_{0}(\underline{z}):s_{\infty}(\underline{z})]
\end{align*}sending each point to the hyperplane in $\CPtwo$ containing it. By removing one of the smooth fibres, for example $\Sigma_{\infty}:=s_{\infty}^{-1}(0)$, we obtain a well defined map  \begin{align}\label{originalfibration}
p: =\frac{s_0}{s_{\infty}}: \CPtwo \setminus \Sigma_{\infty} \longrightarrow \mathbb{C}\mathbb{P}^1 \setminus \{ \infty \} \cong \mathbb{C}.
\end{align} This defines a Lefschetz fibration whose fibres are $2$-spheres with four punctures (the base locus $B$ has been removed in the process) and base $\C$.

\begin{lemma}\label{sympdecom}
	The total space $\CPtwo \setminus \Sigma_{\infty}$ of $p$ is exact symplectomorphic to an open disc subbundle of $(T^*\RPtwo, d\lambda_{T^*\RPtwo})$.
\end{lemma}
\proof Follows from Example \ref{examplerp}. \qed

The fibration \eqref{originalfibration} can be adjusted (Proposition \ref{frompenciltolf}) to become an exact Lefschetz fibration \begin{align}\label{lfadjusted1}
\pi\colon E_{\RPtwo}\to \C
\end{align} with smooth fibre a $2$-sphere with four boundary components, and whose total space admit a symplectic completion to $(T^*\RPtwo, d\lambda_{T^*\RPtwo})$.

The vanishing cycles $V_0,V_1,V_2 \subset M$ are exact Lagrangian circles which partition the four boundary components in three possible configurations of pairs. 
There is only one Hamiltonian class for each vanishing cycle (each homotopy class has only one exact Lagrangian representative, by Stokes's theorem) and since $\Diff(S^1)\simeq O(2)$, there is a unique choice of framing.

\subsection{Monodromy}\label{monodromy1}
In this section we study the properties of the monodromy of the Lefschetz fibration $\pi\colon E_{\RPtwo} \to \C$. Let $\Crit v(\pi)=\{ w_0, w_1, w_2 \}$ be the set of critical values of $\pi$. Fix a base-point $z_* \in \C$, smooth fibre $\pi^{-1}(z_*)\cong M$ (a four punctured sphere) and basis of vanishing paths $(\gamma_0, \gamma_1, \gamma_2)$. Let $(\Delta_0, \Delta_1, \Delta_2)$ and $(V_0, V_1, V_2) \subset M$ be the associated bases of Lagrangian thimbles and vanishing cycles respectively.

\begin{lemma}\label{condition1}
	Let $(d_1, d_2, d_3 ,d_4) \subset M$ be simple closed curves with $d_i \cap V_j=\emptyset$ ($i=1,2,3,4$ and $j=0,1,2$), such that each $d_i$ encircles a distinct boundary component of the fibre.
	The global monodromy $\phi \in \Symp_{ct}(T^*\RPtwo)$ of $\pi$ satisfies \begin{align}\label{lantern}
	\phi:=\tau_{V_0} \tau_{V_1} \tau_{V_2} \simeq \prod_{i=1}^4 \tau_{d_i}.
	\end{align} 
	In particular, $\phi$ commutes with each individual twist, and preserves each vanishing cycle.
\end{lemma}
\proof
The expression \eqref{lantern} is \emph{the} lantern relation (see for example \cite[Proposition 5.1]{bmprimer}), a special instance of the generalised lantern relation of Lemma \ref{lantern1}, according to which the global monodromy is isotopic to the fibred twist in the circle bundle of the normal bundle to the base locus $B\subset \CPone$. In this low-dimensional case, the unit normal bundle is simply the union of the four boundary circles. \qed
\begin{rmk}
	The isotopy \eqref{lantern} also implies that any cyclic permutation of $\tau_{V_0} \tau_{V_1} \tau_{V_2}$ defines the same element of the mapping class group. 
\end{rmk}

\begin{figure}[htb]
	\centering
	\def\svgwidth{200pt}
	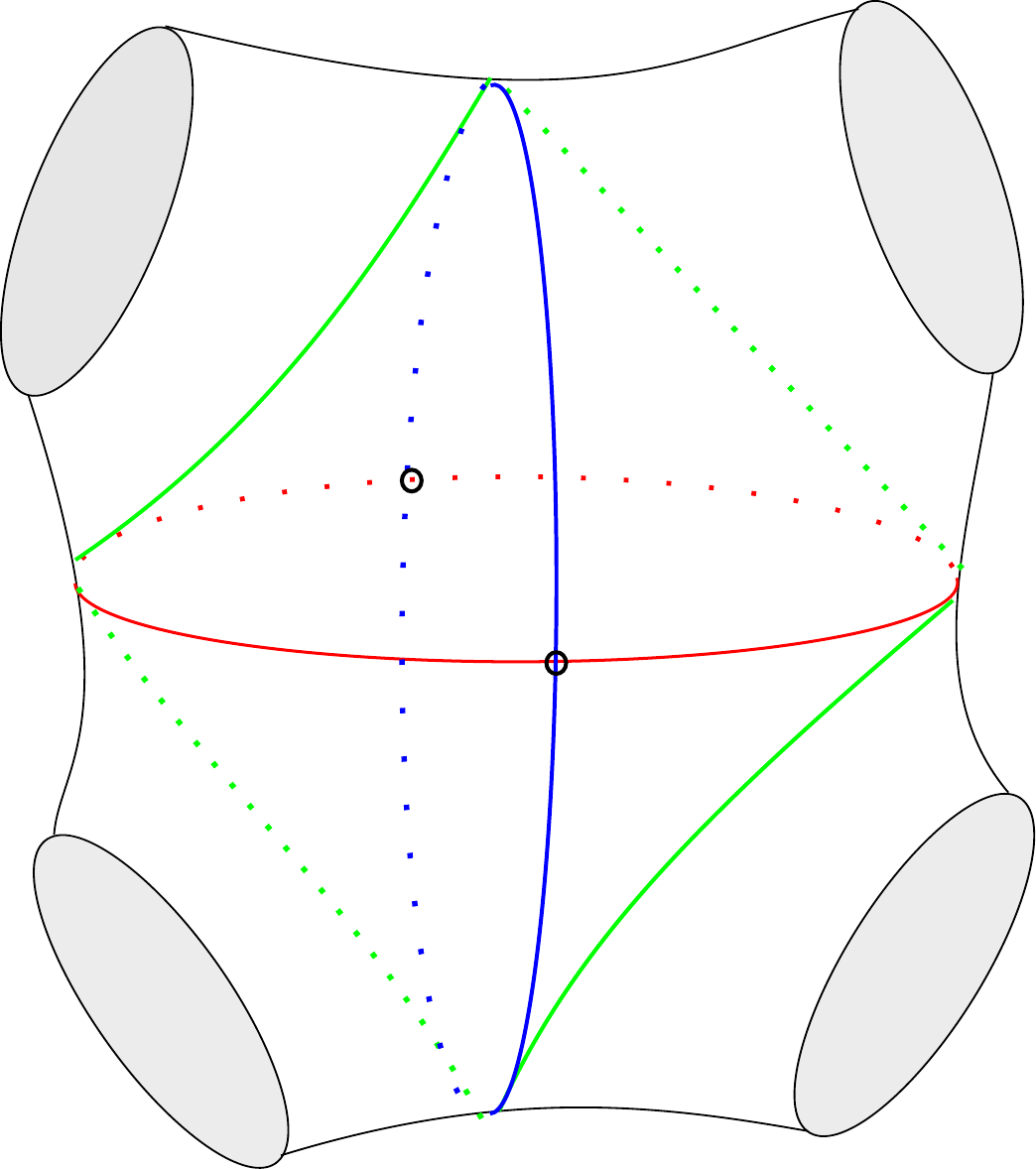
	\caption{ The boundary circles $d_i$ around the boundary components of the smooth fibre are disjoint and therefore the twists in the composition $\Pi_{i=1}^4 \tau_{d_i}$ commute.}\label{fibrecircles}
\end{figure}

\noindent The lantern relation forces the (ungraded) monodromy to act trivially on vanishing cycles. However, if we consider $\phi$ as a graded symplectomorphism, it acts on vanishing cycles by a shift determined as in Lemma \ref{lemmamonodromy}.
The line bundle $\mathcal{L}=\mathcal{O}_{\CPtwo}(2)$ defining the pencil, and the canonical bundle $\mathcal{K}_{\CPtwo}=\mathcal{O}(-3)$ satisfy Assumption \eqref{anticanonical1}, with $d=3$: \begin{align}\label{pencilrp2}
\mathcal{O}(-3)^{\otimes2}\cong \mathcal{O}(2)^{\otimes -3}.
\end{align}
On graded Lagrangian vanishing cycles $V_i \subset M$, we therefore have the shift $\phi(V_i)= V_i[4-3]=V_i[1]$.

\section{Lefschetz fibration on $T^*\CPtwo$}\label{lfcp2}

In this section we discuss an exact Lefschetz fibration $\pi\colon E_{\CPtwo} \to \C$ with smooth fibre a clean plumbing of disc cotangent bundles of $3$-spheres $D_{\eps}(T^*S^3\#_{S^1}T^*S^3)$ along a circle (see \cite{bct} for a reference on the construction of clean Lagrangian plumbing), three singular fibres and whose total space admits a symplectic completion to a Stein manifold exact symplectomorphic to $(T^*\CPtwo, d\lambda_{T^*\CPtwo})$. 

At the end of the section, we add a short digression in which we use the results we have gathered to study Lagrangian submanifolds of $(T^*\CPtwo, d\lambda_{T^*\CPtwo})$ that are diffeomorphic to $S^1\times S^3$. By viewing such manifolds as Lagrangians of the total space of the fibration $\pi\colon E_{\CPtwo}\to \C$, we show that they can be classified in (at least) two Lagrangian isotopy classes (Lemma \ref{lemmaisotopyclass}).

Recall Example \ref{examplecp}. Let $(\CProduct, \omega_{FS}\oplus \omega_{FS})$ with homogeneous coordinates $(\underline{x}, \underline{y}):=([x_0:x_1:x_2], [y_0:y_1:y_2])$ on the first and second factor respectively. Fix a hypersurface \begin{align}
\Sigma:= \left \{ \sum_{i=0}^{2}x_iy_i=0 \right \} \subset \CPtwo \times \CPtwo,
\end{align}
obtained as the image of the embedding of the Flag variety $Fl_3\hookrightarrow \CProduct$. Let $ s_0, s_{\infty}$ be sections of $\mathcal{O}_{\CProduct}(1,1)$. Consider a pencil \begin{align}\label{pencilfirst}
\Sigma_{[\lambda:\mu]}:=	\{ \lambda s_0(\underline{x},\underline{y}) + \mu s_{\infty}(\underline{x},\underline{y})=0 \}_{[\lambda:\mu]\in \CPone} \subset \CProduct.
\end{align} A generic fibre of such a pencil is isomorphic to the three-fold Flag variety $Fl_3 \subset \C^3$.

\begin{lemma}\label{lemmabaselocus}
	The base locus $B$ is diffeomorphic to the 3-point blow up of $\CPtwo$, $B\cong \CPtwo\#3\overline{\CPtwo}$, equipped with its monotone symplectic form.
\end{lemma}
\proof The base locus is a symplectic manifold of dimension $4$, and since it is obtained as the intersection of hyperplane sections of $\CProduct$ of bidegree $(1,1)$, its Chern class is a positive class by the adjunction formula. Hence $B$ is a monotone Fano (or del Pezzo) surface, i.e a non-singular projective algebraic surface with ample anticanonical divisor. Thus $B$ is either isomorphic to $\CPone \times \CPone$ or to the blow up $\CPtwo\#r\overline{\CPtwo}$ at $r$ points in general position, for $0 \leq r \leq 8$ (see for example \cite{delpezzo}).

Embed $\CPtwo \times \CPtwo$ in $\C\PP^{8}$ via the Segre embedding, and by Lemma \ref{lemmauler}, we know \begin{align}
\chi(B)= 2\chi(Fl_3)-\chi(\CPtwo \times \CPtwo)+3.
\end{align}
Now $\chi(\CPtwo \times \CPtwo)=\chi(\CPtwo)^2=9$ and the flag manifold $Fl_3$ is known to have Euler characteristic $\chi(Fl_3)=3!=6$,	so that $\chi(B)=6$, and hence $B\cong \CPtwo\#3\overline{\CPtwo}$.
\endproof

Given a generic pencil of $(1,1)$ divisors \eqref{pencilfirst}, whose smooth fibre ``at infinity'' is denoted by $\Sigma_{\infty}:=s_{\infty}^{-1}(0) \cong Fl_3$ 
consider the map 
\begin{align}\label{lf1}
p= \frac{s_0(x,y)}{s_{\infty}(x,y)}\colon	\CProduct \setminus \Sigma_{\infty} \longrightarrow \C.
\end{align}

\begin{lemma}[{{see Lemma \ref{jonnylemma}}}]\label{affinefibre}
	The fibre of \eqref{lf1}, which is the affine $3$-fold $Fl_3 \setminus \CPtwo\#3\overline{\CPtwo}$ is the interior of a Stein domain whose symplectic completion is exact symplectomorphic to a clean plumbing $T^*S^3\#_{S^1}T^*S^3$.
\end{lemma}

\begin{lemma}\label{sympdecom2}
	The total space $\CProduct \setminus \Sigma_{\infty}$ of $p$ is exact symplectomorphic to an open disc subbundle of $(T^*\CPtwo, d\lambda_{T^*\CPtwo})$.
\end{lemma}
\proof Follows from Example \ref{examplecp}. \qed

The fibration \eqref{originalfibration} can be adjusted (Proposition \ref{frompenciltolf}) to become an exact Lefschetz fibration \begin{align}\label{lfadjusted2}
\pi\colon E_{\RPtwo}\lra \C
\end{align} whose total space admit a symplectic completion to $(T^*\CPtwo, d\lambda_{T^*\CPtwo})$. In the next section, we study the fibres of $\pi$.

\subsection{Topology of the fibre: A MBL fibration on the affine Flag 3-fold}\label{mblfib}
In this section we enhance our understanding of the topology of the fibres of the Lefschetz fibration $\pi \colon E_{\CPtwo} \to \C$ via a construction due to Jonny Evans \cite{evansnotes}. We will show why Lemma \ref{affinefibre} is true.

Let $Fl_3=\{ x_0y_0+ x_1y_1+x_2y_2=0 \} \subset \CProduct$ be the Flag $3$-fold as in Example \ref{examplecp} (where we use the same choice of coordinates on $\CProduct$ as before). Consider a pencil of divisors $Y_{[\lambda :\mu]}:=\{ \lambda(x_0y_0-x_2y_2)+\mu x_1y_1=0 \}_{[\lambda:\mu]\in \C\PP^1} \subset Fl_3$. By the definition of the pencil, generic fibre of $Y_{[\lambda:\mu]}$ can be understood as the intersection of two hyperplane sections of $\CProduct$ of bidegree $(1,1)$. Therefore, by Lemma \ref{lemmabaselocus}, the generic fibre is a copy of the del Pezzo surface $\CPtwo\#\overline{3\CPtwo}$. There are three critical fibres at $[\lambda: \mu]\in \{ [1:0], [1:1], [-1:1] \}$. The base locus of this pencil is the set $\Theta:=\{ x_0y_0=x_1y_1=x_2y_2 =0 \}$, composed of six lines $L_{ijk}'=\{ x_i=x_j=y_k=0 \}$, $L_{ijk}'=\{ y_i=y_j=x_k \}$ where $ijk$ is a permutation of $012$. 
For all $i$, let $z_i=x_iy_i$. Remove a smooth fibre $Y_{\infty}$ and define a Morse--Bott Lefschetz fibration (see the proof of \cite[Lemma 4.5]{smithwemyss})\begin{align}\label{mbl}
p_1\colon Fl_3 \setminus Y_{\infty} \to \C, \ ([x_0:x_1:x_2], [y_0:y_1:y_2]) \mapsto \frac{x_1y_1}{x_2y_2-x_0y_0}.
\end{align}
The general fibre is isomorphic to $(\C^*)^2$ and there are three singular fibres over $\{ 0 , \pm 1 \}$ isomorphic to $\C^*\times (\C\vee_0\C)$, so that in every singular fibre, there is a copy of $\C^*$ that is the singular locus of that fibre.

\begin{lemma}[{{\cite{evansnotes}}}]\label{matchings}
	The fibration $p_1$ admits three $3$-dimensional matching spheres of the total space, which pairwise intersect cleanly along a circle.
\end{lemma} 

\proof 

The symplectic structure on $Y\setminus Y_{\infty}$ yields a preferred choice of symplectic connection (the symplectic orthogonal to each fibre of $p_1$) which gives local parallel transport maps.

We first explain why these parallel transport maps are globally defined on the total space of the fibration $p_1$. There is a Hamiltonian $T^2$-action on $\CProduct$ given by \begin{align}\label{torica}
(\theta, \phi)\cdot([x_0:x_1:x_2],[y_0:y_1:y_2]) = ([e^{i\theta}x_0:x_1:e^{i \phi}x_2], [e^{-i\theta}y_0:y_1:e^{-i\phi}y_2]), \end{align}

which preserves the fibres of $p_1$. The (restricted) moment map $H\colon Fl_3\setminus Y_{\infty} \to \R^2$ of \eqref{torica} has the shape 
\begin{align}
H([x_0:x_1:x_2], [y_0:y_1:y_2])=\frac{1}{2}\left( \frac{|x_0|^2}{|a|^2}-\frac{|y_0|^2}{|b|^2}, \frac{|x_2|^2}{|a|^2}-\frac{|y_2|^2}{|b|^2} \right)
\end{align}
and its level sets are compact. These level sets are preserved by the parallel transport maps (see \cite[Lemma 4.1]{smithwemyss} for a proof). Therfore, parallel transport maps are defined globally via the action \eqref{torica}.

Let $\gamma_{\pm}(t)=\pm t$, $t\in [0,1]$ be the path in $\C$ connecting the critical value $0$ to $\pm 1$. Over $z\notin \{ -1, 0, 1\}$, the restriction $T_z:=H^{-1}(0) \cap p_1^{-1}(z)$ is a smooth torus in the smooth fibre $(\C^*)^2$. Since parallel transport preserves the level sets of $H$, set $S_{\pm}:=H^{-1}(0)\cap p_1^{-1}(\gamma_{\pm}([0,1]))$ obtained by parallel transporting such a torus along the path $\gamma_{\pm}$, defines a matching cycle. To understand $S_{\pm}$, we look at how the torus $T_z$ degenerates in the critical fibres over $z\in \{-1,0,1 \}$.

The action \eqref{torica} exhibits three circles given by $\{ \phi=0 \}, \{ \theta=0\}$ and $\{\phi-\theta=0 \}$. We can set the generators of $H_1(T^2;\Z)\cong \Z^2$ to be the $\theta$- and $\phi$-circles. We explain below that each circle collapses in one of the singular fibres.


\begin{figure}[h]
	\centering
\includegraphics[width=9cm]{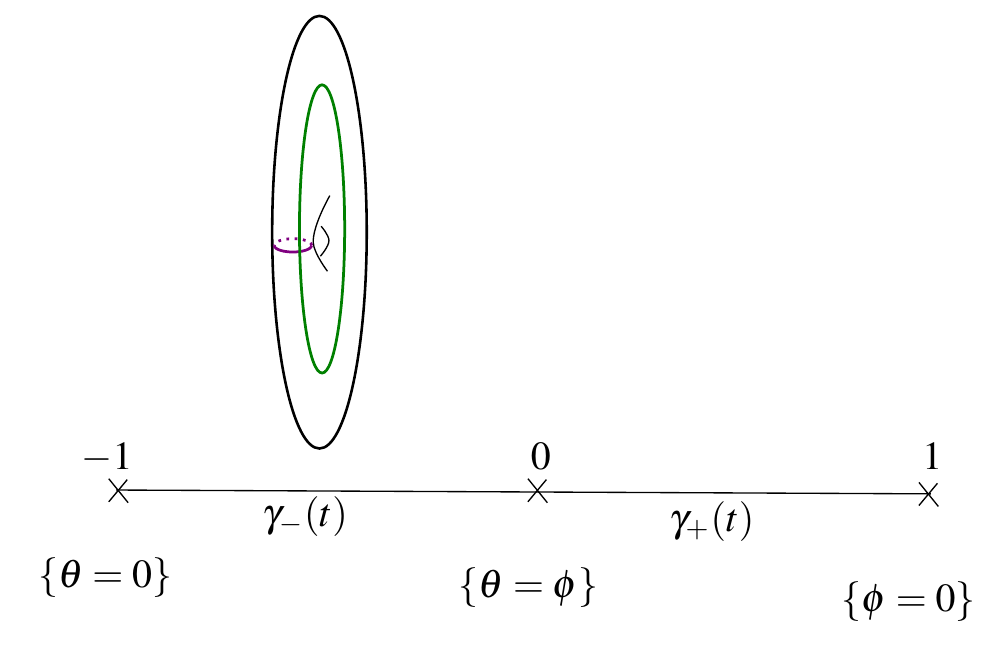}
	\caption{Depiction of the restriction of the smooth fibre $T_z=H^{-1}(0)\cap p_1^{-1}(z)$, such that over each singular value one of the circles $\{ \theta=0\}$, $\{ \theta=\phi \}$, $\{\phi=0 \}$ collapses.}
	\label{mbltorus}
\end{figure}

The critical locus over $z=-1$ is given by the equations $\{ x_0=y_0=0, \ x_1y_1 +x_2y_2 =0 \}$ so it is made of pairs $([0:x_1:x_2], [0:-x_2:x_1])$. Therefore, as $t\to -1$ on $\gamma_-$, the $\{  \theta=0 \}$ circle must collapse to a point.

The critical locus over $0$ is given by the equations $\{ x_1=y_1=0, \ x_0y_0+x_2y_2=0$, so it is made of pairs $([x_0:0:x_2], [-x_2:0:x_0])$. The $T^2$ action becomes trivial if $\theta=\phi$, which means that at $t\to 0$ on $\gamma_{\pm}$, the circle $\{ \theta=\phi \}$ collapses to a point. 

The critical locus over $1$ is given by $\{ x_2=y_2=0, \ x_0y_0=-x_1y_1\}$ so it is made of pairs $([x_0:x_1:0], [x_1:x_0:0])$. In this case, as $t\to 1$ on $\gamma_+$, the $\{ \phi=0\}$ circle collapses to a point.

From these observations, one can see that the matching cycles $S_{\pm}$ are unions of solid tori $S^1\times D^2\cup_{S^1\times S^1} D^2\times S^1$ glued along the torus $S^1\times S^1\subset \C^*$ in the smooth fibre; see also \cite[Example 4.2]{smithwemyss}. This union has to be a three sphere, because the circles $\{ \theta=0\}$ and $\{ \phi=0 \}$ become trivial homology classes in $S_{\pm}$.

By this description we also understand that the matching spheres $S_{\pm}$ intersect in $p_1^{-1}(0)$ cleanly in a circle (see also discussion of \cite[Lemma 4.5]{smithwemyss}, these intersection loci are the real part of the three components of the critical locus).
A third matching sphere can be obtained as the Morse--Bott surgery of the other two (there are two possible such surgeries).

\endproof

\begin{lemma}[\cite{evansnotes}]\label{jonnylemma}
	The union of $S_{\pm}$ form a Lagrangian skeleton for the clean Lagrangian plumbing $T^*S^3 \#_{S^1}T^*S^3$.
	In particular, the symplectic completion of $M:=Fl_3 \setminus Y_{\infty}$ is symplectomorphic to this plumbing.
\end{lemma}
The proof of this lemma is beyond the scope of this section. The idea in \cite{evansnotes} is to define a pluri-subharmonic function $Fl_3 \setminus Y \to \R$ with skeleton formed by the pair of Lagrangians $3$-spheres over $[-1,0]$ and $[0,1]$, meeting cleanly in a circle over $0$. 
A clean plumbing of spheres $T^*S^3\#_{S^1}T^*S^3$ is determined by the topology of the surgery of the Lagrangian cores (\cite{smithwemyss}), and the plumbing of the statement is such that the surgery of the core spheres is another sphere.

\subsection{Topology of the fibre: a bifibration}\label{bifibration}

Consider $(1,1)$-divisors on $(\CProduct, \omega_{FS}\oplus \omega_{FS})$ given by a linear combination of monomials $ x_0y_0, x_1y_1, x_2y_2$. This gives a rational map \begin{align}\label{net}
\CProduct & \dashrightarrow\CPtwo \\
([x_0:x_1:x_2], [y_0:y_1:y_2]) &\longmapsto [z_0=x_0y_0:z_1=x_1y_1:z_2=x_2y_2].
\end{align} 
The base locus of this map is the set $\Theta$ as in Section \ref{mblfib}, so consider the well-defined restriction $$q\colon \CProduct \setminus \Theta \to \CPtwo.$$ 

A generic fibre over $z=[z_0:z_1:z_2]$ is determined by the equations $ \{ x_iy_i=z_i \}$. The smooth fibre is $(\C^*)^3$. However, it reduces to $(\C^*)^2$ after quotienting by the $\C^*$ action given by $(x_i,y_i)\mapsto (\zeta x_i, \zeta^{-1}y_i)$ for all $i \in \{ 0,1,2 \}$, $\zeta \in \C^*$.

Over the coordinate lines $\{ z_i=x_iy_i=0 \} \setminus \{ [1:0:0], [0:1:0], [0:0:1] \}$ of $\CPtwo$, the fibration is singular with singular fibres $\C^*\times (\C \vee_0 \C)$. There are three such degenerations, one for each factor of $(\C^*)^3$.

\begin{lemma}
	There is a Lefschetz fibration over $\C$, with total space exact symplectomorphic to a disc cotangent bundle of $T^*\CP$ and fibres exact symplectomorphic to the affine Flag $3$-fold (see Lemma \ref{affinefibre}) with the following property. The three vanishing cycles (associated to a basis of vanishing paths) of this fibration are three-spheres which pairwise intersect cleanly in a circle.
	
\end{lemma}
\proof To show the claim, we consider a $1$-parameter family of lines in \eqref{net}.
For every generic line $\ell:=\{ b_0z_0+b_1z_1+b_2z_2=0, \ b_i\neq 0 \}$ in the base $\CPtwo$, the restriction $q|_{\ell}$ gives a $(\C^*)^2$-fibration of an open subset of the Flag 3-fold in $\CProduct$. There are three singular fibres, over the points in which $\ell$ intersects the coordinate axes.

A deformation of $\ell$ in a family $\ell_{[\alpha:\beta]}$ of lines in $\CPtwo$ corresponds to a pencil of flag 3-folds in $\CProduct$, and the restriction of $q$ over $\ell_{[\alpha:\beta]}$, $[\alpha:\beta] \in \C\PP^1$ defines a Lefschetz fibration $\pi$ as in Section \ref{lfcp2}. 

In this family, consider a path $\delta(t)$ joining a generic line $\delta(0):=\ell_{gen}$ to a line $\delta(1)=\ell_{crit}$ passing through either of $[1:0:0], [0:1:0], [0:0:1]$; in the total space this traces a family of Morse--Bott--Lefschetz fibrations, each with a configuration of three singular fibres. We can choose $\ell_{gen}$ such that $q|_{\ell_{gen}}=p_1$ from the previous section.
As the path reaches $\ell_{crit}=\delta(1)$, two of the critical loci of this configuration will come together over the singular restriction $q|_{\ell_{crit}}$ and this produces a matching sphere of the fibration $p_1$, as described in Lemma \ref{matchings}.
Namely, for any $w \in \{ [1:0:0], [0:1:0], [0:0:1]\}$, $q^{-1}(w)=(\C\vee \C)\times  (\C\vee \C)$: as two of the singular points of $q$ come together, \emph{two} of the original factors in $(\C^*)^3$ collapse, and after quotienting by the $\C^*$-action, the degenerations are recognisable as ordinary double points. 
As the restriction $q|_{\ell_{crit}}$ defines a singular fibre for a fibration of the type of $\pi$, the matching spheres of $p_1$ (Lemma \ref{matchings}) coincide with the vanishing cycles of $\pi$.

\qed

\begin{rmk}
	Any smooth fibre $q^{-1}(z)$, $z=[z_0:z_1:z_2]$ with $z_i \neq 0$ compactifies to $q^{-1}(z) \cup \Theta =\CPtwo\#3\overline{\CPtwo}$. Note that this coincides with the base locus of the pencil in Lemma \ref{lemmabaselocus}. 
	For the closure of the singular fibres, see \cite[Section 2]{tyurin}.
	
\end{rmk}

\subsection{Monodromy}\label{lanterncp2}
Let $\pi \colon E_{\CPtwo} \to \C$ be the exact Lefschetz fibration discussed in the previous subsections, with $\Crit v=\{ w_0, w_2, w_4 \}$. Let $z_*\in \C$ be a basepoint, which fixes the smooth fibre $\pi^{-1}(z_*) \cong M$, which admits a symplectic completion to the plumbing $T^*S^3\#_{S^1}T^*S^3$. Choose a distinguished basis of vanishing paths $(\gamma_0, \gamma_2, \gamma_4)$ and associated vanishing thimbles $(\Delta_0, \Delta_2, \Delta_4)$ (so any two of them, $\Delta_j, \Delta_k$ are disjoint and satisfy $h(\Delta_j)>h(\Delta_k)$ if $j>k$). Let $(V_0, V_2, V_4) \subset M$ be the basis of vanishing cycles associated to the three singular points, from which we define the Dehn twists $\tau_{V_0}, \tau_{V_2}, \tau_{V_4} \in \Symp_{ct}(M)$. The global monodromy of the fibration $\pi$ is isotopic to $\tau_{V_0}\tau_{V_2}\tau_{V_4}.$


\begin{lemma}\label{monodromycp2}
	The global monodromy satisfies the generalised lantern relation \begin{align}\label{lantern2}
	\tau_{V_0}\tau_{V_2}\tau_{V_4}\simeq \tau_{V} \in \pi_0(\Symp_{ct}(M))
	\end{align}
	where $\tau_{V}$ is the fibred twist in the unit normal bundle $V\to B$ to base locus $B=\CPtwo\#3\overline{\CPtwo} \subset \CProduct$.
\end{lemma}
\proof 
This follows from Lemma \ref{lantern1}, for the pencil defined in Section \ref{lfcp2}. \endproof 

\begin{cor}\label{commutingrel}
	The global monodromy commutes with every single twist $\tau_{V_i}$ for $i=0,2,4$.
\end{cor}
\proof
Follows from Lemma \ref{monodromycp2}. 
\qed

Recall that the monodromy can be made into a graded symplectomorphism, and that for the pencil we consider, the canonical bundle $\mathcal{K}_{\CProduct} \cong \mathcal{O}_{\CProduct}(-3,-3)$ satisfies Assumption \eqref{anticanonical1} with $d=6$: 
\begin{align}\label{pencilcp}
\mathcal{K}_{\CProduct}^2=\mathcal{L}^{\otimes - d}
\end{align}for the line bundle $\mathcal{L}:=\mathcal{O}_{\CProduct}(1,1)$ generating the pencil. So if $V_i\subset M$ is a graded vanishing cycle, we have $$\phi(V_i)= V_i[-2].$$

\subsection{Monotone Lagrangian submanifolds of $T^*\CPtwo$}\label{sectiontori}
This final section is a short digression in which we include an application derived from the construction of Section \ref{lfcp2} to discuss monotone Lagrangian submanifolds of $(T^*\CPtwo, \lambda_{T^*\CPtwo})$ diffeomorphic to $S^1 \times S^3$. We show that there are at least two distinct Lagrangian isotopy classes of such Lagrangians (\ref{lemmaisotopyclass}).

The Lagrangians we study in this section can be obtained by parallel transport of the vanishing cycles of the exact Lefschetz fibration $\pi\colon E_{\CPtwo} \to \C$ defined in Section \ref{lfcp2}. The construction of such Lagrangians imitates that of two well known families of Lagrangians in $T^*S^2$; the Chekanov and Clifford tori. These tori can be defined in various ways (see below), but in particular they do admit a presentation as Lagrangian submanifolds of the Lefschetz fibration on $(T^*S^2, d\lambda_{T^*S^2})$, a fibration obtained from a Lefschetz pencil of conics on $(\CPone\times \CPone,\omega_{FS}\oplus\omega_{FS})$. 
These two tori are not Hamiltonian isotopic to each other, and the motivation behind this section is to disclose whether the Lagrangians $S^1\times S^3 \subset T^*\CPtwo$ obtained in the analogous way (in the total space of $E_{\CPtwo}\to \C$) can be divided into two (or more) distinct Lagrangian/Hamiltonian isotopy classes. 

The Chekanov torus is an important type of monotone torus originally studied in its Lagrangian embedding in $(\C^n,\omega_{std})$ (\cite{chekanov} and later (\cite{chesch}) in complex projective spaces $\CP$ and products of spheres $\times_n\CPone$ (in particular in $(\CPtwo, \omega_{FS})$ and $(\CPone\times \CPone,\omega_{FS}\oplus\omega_{FS})$). In all such examples, this kind of torus stands out as being not Hamiltonian isotopic to the Clifford torus, which is the ``standard'' product of circles $\times_nS^1(r)$.

Consider the decompositions of the Kähler manifolds 
\begin{align}
\CPtwo \cong D_{\eps}T^*\RPtwo \cup C_{\infty}\\
\CPone\times \CPone \cong D_{\eps}T^*S^2 \cup \Delta
\end{align}
(see Section \ref{modelfibred}, the second is a special case of Example \ref{examplecp} that can also be found in \cite{biran1,audin2} as referred above) where $C_{\infty}\subset \CPtwo$ is the quadric at infinity and $\Delta \subset \CPone\times \CPone$ is the holomorphic diagonal. In \cite{gadbled}, Gadbled proved that the Chekanov tori in $(\CPtwo, \omega_{FS})$ and $(\CPone\times \CPone,\omega_{FS}\oplus\omega_{FS})$ can be constructed as circle bundles over $C_{\infty}$ and $\Delta$ respectively (this was proved using the \emph{circle bundle construction} of \cite{biran2}). Therefore, these tori must be preserved in the complement of the divisors in the decompositions above, giving rise to Lagrangian tori in the cotangent bundles $(T^*\RPtwo, d\lambda_{T^*\RPtwo})$ and $(T^*S^2,d\lambda_{T^*S^2})$.

The well-studied Lefschetz fibration on $(T^*S^2, d\lambda_{T^*S^2})$ with smooth fibre $(T^*S^1, d\lambda_{T^*S^1})$, has two critical fibres, one vanishing cycle $V_0$ (the zero section of the smooth fibre) and a matching sphere corresponding to the zero section.
This fibration arises from a pencil of conics on $(\CPone\times \CPone, \omega_{FS}\oplus\omega_{FS})$ (whose base locus is a pair of points), by removing the conic at infinity, the holomorphic diagonal $\CPone \cong \Delta \subset \CPone \times \CPone$ (see \cite[Example 3.2.1]{biran1} and \cite[Section 4.3]{audin2}).

It is known that the Chekanov torus in $T^*S^2$ admits a presentation as Lagrangian submanifold of the total space of the Lefschetz fibration $T^*S^2\to \C$ above (see \cite[5.1]{aurouxanti}). In this picture, the torus emerges from flowing the vanishing cycle $V_0$ by parallel transport over a loop that does not encircle any of the critical values. Then, the Chekanov torus is fibred by $V_0$ over this loop.

Recall the decomposition \begin{align}\label{decompocp2}
\CProduct\cong D_{\eps}T^*\CPtwo \cup \Sigma_{\infty}
\end{align}
where the divisor at infinity is the Flag $3$-fold.

\begin{quest}\label{questpc2}
	Can we use \ref{decompocp2} to get an interesting Lagrangian $S^1\times S^3\subset T^*\CPtwo$, comparable (in its construction) to the Chekanov torus in $T^*S^2$?

\end{quest}
We could attempt to answer Question \ref{questpc2} by applying Biran's circle bundle construction (\cite{biran2}) over a three-sphere contained in $\Sigma_{\infty}$, to generalise the methods used for the Chekanov tori. However, there is no distinguished choice of sphere in the divisor $\Sigma_{\infty}$, contrasting the Chekanov torus construction in which the circle bundle is taken over the equator in $S^2\cong\Delta \subset \CPone\times \CPone$.

We discard this approach, and privilege another, which is to consider the products $S^1\times S^3$ as Lagrangians in the total space of the Lefschetz fibration $\pi\colon E_{\CPtwo} \to \C$.

\subsubsection{Lagrangian in the total space of the standard fibration on $\C^4$}
We begin by studying monotone Lagrangian submanifolds of $\C^4$ obtained by parallel transport in the total space of the ``standard fibration'' \begin{equation}
\begin{split}\label{stdmap}
q\colon\C^4& \lra \C \\
z&\longmapsto  z_1^2+z_2^2+z_3^2+z_4^2
\end{split}
\end{equation}

with smooth fibres $\ q^{-1}(z) = \{ (x,y)\in \C^4, \ |x|^2-|y|^2=1, \ \langle x, y \rangle =1 \} \cong T^*S^3$ for $z\neq 0$ and a singular fibre over $0$, in which the zero section of $T^*S^3$ collapses to a point. The monodromy around a loop in the base $t \mapsto e^{2\pi i t}$ is then given by the Dehn twist along this vanishing cycle, i.e $\tau_{S^3} \in \Symp_{ct}(T^*S^3)$. For $z\in \C^*$, let $V_z \cong S^3\subset q^{-1}(z)$ be the zero section of the fibre over $z$ (a representative of the vanishing cycle).

We construct two (distinct) families of Lagrangian $S^1\times S^3 \subset \C^4$ as follows. Let $\sigma \colon [0,1]\to \C$ be a loop in the base and let
\begin{align}
T_{\sigma}:=\bigcup_{z\in \im(\sigma)}  V_z\subset \C^4. \end{align}
This is a Lagrangian submanifold of the total space that can be obtained by flowing the vanishing cycle by parallel transport over the loop $\sigma$.
If $0 \in \im(\sigma)$ then $T_{\sigma}$ is an immersed four-sphere (with a nodal singularity at $0$). If, on the other hand, $0 \notin \im(\sigma)$, we call $T_{\sigma} \cong S^3 \times S^1$ of type $1$ if $\sigma$ encloses the origin, and of type $2$ otherwise.


\begin{figure}[h]
	\centering
\includegraphics[width=14cm]{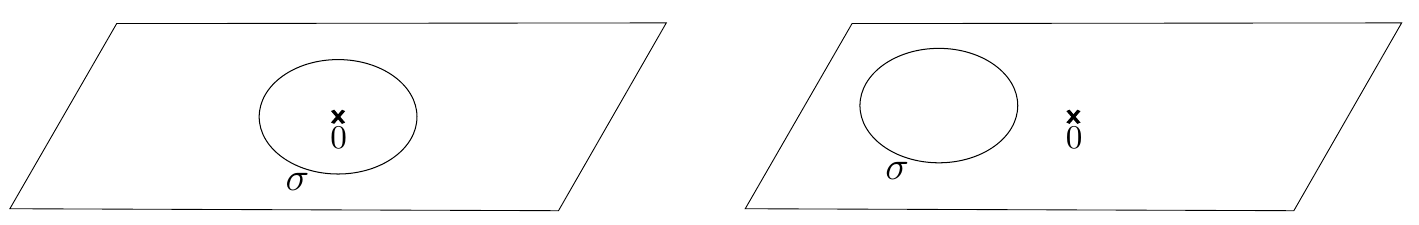}
	\caption{Type 1 (left) and Type 2 (right) Lagrangians.}
\end{figure}

For simplicity, consider the exact Lefschetz fibration $f\colon E\to \C$ obtained by cutting the the fibres into Liouville domains exact symplectomorphic to a disc cotangent bundles $D_rT^*S^3$, $r>0$, and call the new fibration $f\colon E\to \C$ (see \cite[Section 1.2]{seideles} for the precise description of the exact local model).

\begin{lemma}
	Let $0 \notin \im(\sigma)$. The Lagrangian $T_{\sigma}$ is monotone.
\end{lemma}
\proof 
Recall that a Lagrangian $L\subset (X, \omega)$ is monotone if the the area homomorphism $\omega\colon \pi_2(X,L)\to \R$ and the Maslov (index) homomorphism $\mu\colon \pi_2(X,L)\to \Z$ are proportional, i.e \begin{align}
\forall u\in \pi_2(X,L), \ \omega(u)=c \mu(u).
\end{align}
By the homotopy long exact sequence for pairs, we have $\pi_2(\C^4,T_{\sigma})\cong \pi_1(T_{\sigma}) \cong \Z$, for $T_{\sigma}\cong S^3 \times S^1\subset \C^4$.
Then, the two homomorphisms are forced to be proportional.
\endproof

The characterisation of pseudoholomorphic discs we present below revolves around the use of a section-count invariant associated to Lefschetz fibrations as defined in \cite[Section 2.1]{seideles}.
Let $(J,j)$ be a pair of almost complex structures on $\C^4$ and $\C$ compatible with $\pi$ (as in Definition \ref{compatiblej}).

Let $D\subset \C$ be the disc in the base bounded by the loop $\sigma$, such that $\partial D=\im(\sigma)$, and consider pseudoholomorphic sections of the restriction $f|_D$, the curves $s\colon D \to E$ satisfying \begin{enumerate}
	\item $\forall z\in D$: $f (s(z))=z$,
	\item $s(\partial D)=s(\sigma) \subset T_{\sigma}$ (i.e the section $s$ has a boundary condition defined by $T_{\sigma}$),
	\item $Ds(z)+J(s)\circ Ds(z)\circ j =0$.
\end{enumerate}

Let $\mathcal{M}_{E/D}:=\{ s\colon D \to E \text{ satisfies } (1), (2), (3) \}$ be the moduli space of $(J,j)$-holomorphic sections described above. A Fredholm analysis can be adapted to this situation to show that $\mathcal{M}_{E/D}$ is a smooth manifold (see \cite[Lemma 2.5]{seideles}, \cite[p.237]{seidelbook}).
The choice of almost complex structures forces $J$-holomorphic discs in the total space to project to $j$-holomorphic discs in the base, so by the open mapping theorem, the moduli space $\mathcal{M}_{E/D}$ corresponds to the moduli space of $J$-holomorphic discs $u\colon (D, \partial D)\to (E, T_{\sigma})$.
Moreover, if $0 \in int(D)$, the moduli space $\mathcal{M}_{E/D}$ can be identified with the unit cotangent bundle $ST^*S^3$ (see \cite[p.1033]{seideles}).

\begin{lemma}\label{maslovnumberdiscs1}
	Assume the loop $\sigma\colon S^1\to \C$ bounding $D\subset \C$ is such that $0\in int(D)$. Then the Lagrangian $T_{\sigma}$ has minimal Maslov index $4$.
\end{lemma}

\proof Let $N_{T_{\sigma}}$ be the minimal Maslov index associated to $T_{\sigma}\subset \C^4$ (see Definition \ref{maslovnumbers}). The moduli space $\mathcal{M}_{E/D}$ of $J$-holomorphic discs $s \colon (D, \partial D)	\to	 (X,T_{\sigma})$ with boundary in the homotopy class of the generating loop of $\pi_1(L)$
has expected (virtual) dimension $\dim (\mathcal{M}_{E/D})=4+N_{T_{\sigma}}-3$ (this is a standard result that derives from the study of the Fredholm operator involved in the equations defining pseudo-holomorphic discs, see for example \cite{oh}). By the observation above, we know $\dim(\mathcal{M}_{E/D})=\dim(ST^*S^3)=5$, so that $N_{T_{\sigma}}=4$. \endproof

\begin{lemma}\label{chektype}
	Assume the loop $\sigma\colon S^1\to \C$ bounding $D\subset \C$ is chosen such that $0\notin D$. Then the Lagrangian $T_{\sigma}$ has minimal Maslov index $2$.
\end{lemma}
\proof Since $0 \notin D$, the restriction $f|_D\colon E|_D \to D$ can be trivialised so that $E|_D=f^{-1}(D)\cong D_rT^*S^3 \times D$.
The only non-trivial discs are generated by the base loop $\sigma$ onto which $T_{\sigma}$ projects to (there is no contribution from the smooth fibre since $T_{\sigma}|_{z}:=T_{\sigma}\cap f^{-1}(z) \subset f^{-1}(z) \cong T^*S^3$ are both simply connected). Any disc generated this way has Maslov index 2.

\endproof

\begin{rmk}
	Lemma \ref{maslovnumberdiscs1} can be somewhat surprising at first. Namely, one could be mislead by the two-dimensional situation and expect the index of those discs to be two. 
	The minimal Maslov index of any non-trivial disc with boundary conditions on any Lagrangian torus in $(\C^2,\omega_{std})$ is two (\cite{viterbo}), so in particular that is the case for the Clifford torus (which would be a type-one torus in the total space of the model Lefschetz fibration on $\C^2$). This is also the result obtained from the expected dimension of the moduli space of such discs, which is $\dim(ST^*S^2)=3$ (so that the minimal Maslov index is indeed $3-2+1=2$).
	But this only holds in dimension two; it is incorrect to think that this would generalise to the discs of Lemma \ref{maslovnumberdiscs1}.

\end{rmk}

\begin{cor}
	The two Lagrangians are distinguished by their Maslov index, hence they are neither Lagrangian isotopic, nor Hamiltonian isotopic.
\end{cor}
\qed 



\subsubsection{Monotone Lagrangian $S^3\times S^1\subset T^*\CPtwo$}

Let $\pi\colon E_{\CPtwo} \to \C$ the exact Lefschetz fibration \eqref{lfadjusted2}, with critical values $\Crit v(\pi)=\{ w_0,w_2,w_4 \}$ and vanishing cycles $V_0,V_2,V_4$. For $i=0,2,4$, set $V_{z,i}$ to be the representative of $V_i$ in the fibre $\pi^{-1}(z)$.

For each critical value $w_i\in \Crit v(\pi)$ there is a pair of Lagrangians $S^3\times S^1 \subset T^*\CPtwo$ as follows.
\begin{definition}\label{lagincp2}
	For $i\in \{ 0,2,4\}$, fix $w_i\in \Crit v(\pi)$. Let $\sigma\colon S^1\to \C$ be a loop in the base of the fibration $\pi\colon E_{\CPtwo}\to \C$, bounding a disc $D$ with $\partial D=\im(\sigma)$, $\Crit v(\pi)\cap \im(\sigma)=\emptyset$ and $w_j \notin D$ if $j\neq i$.
	
	There is a monotone Lagrangian \begin{align}\label{lagrangianwi}
	T_{\sigma, w_i}:=\bigcup_{z\in \im(\sigma)}V_{z,i} \subset E_{\CPtwo},
	\end{align}
	that we can view as Lagrangian of $(T^*\CPtwo, d\lambda_{T^*\CPtwo})$, obtained by flowing the vanishing cycle $V_i$ under parallel transport around the loop $\sigma$. 
	
	We say $T_{\sigma, w_i}$ is
	\begin{enumerate}
		\item of Type $1$ if $w_i \in int(D)$, 
		\item of Type $2$ if $w_i\notin D$.
	\end{enumerate}
\end{definition}

\begin{lemma}\label{lemmaisotopyclass}
	Fix $w_i \in \Crit v(\pi)$. Let $\sigma_1, \sigma_2\colon S^1\to \C$ be two loops as in Definition \ref{lagincp2}, and assume that the Lagrangian $T_{\sigma_1, w_i}$ is of type $1$, and $T_{\sigma_2,w_i}$ is of type $2$, both associated to the same critical value. Then the Lagrangians $T_{\sigma_1,w_i}$ and $T_{\sigma_2,w_i}$ are not Lagrangian isotopic.
\end{lemma}

\proof 
Consider the homotopy long exact sequence for the pair $(X,L)=(T^*\CPtwo, S^3\times S^1)$. Exactness of $X$ and monotonicity of $L$ imply that the Maslov homomorphism $\mu\colon \pi_2(X,L)\to \Z$ descends to a map $\mu\colon \pi_1(L)\to \Z$, which represents an element of $H^1(L;\Z)$. Therefore, we can use the considerations of the previous section to prove the claim. \endproof

\section{A symplectomorphism on the total space of a Lefschetz fibration}\label{chaptermodels}

 Given a Lefschetz fibration $\pi\colon E\to \C$ with certain properties, we use Picard-Lefschetz theory together with a construction from \cite{seideliterated}, to build a compactly supported symplectomorphism $\varphi \in \Symp_{ct}(E)$. The construction (presented in Section \ref{construction}) involves lifting a Dehn twist along an annulus in the base $\C$ to a symplectomorphism of the total space, that we adjust to a compactly supported symplectomorphism.

In Section \ref{hfcomputationsection} we measure the Floer theoretical action of $\varphi$ on a Lefschetz thimble by computing the Floer cohomology groups $\HF(\varphi^k(\Delta_{\alpha}), \Delta_{\beta};\Zmod)$ for elements $\Delta_{\alpha}, \Delta_{\beta}$ in a distinguished basis of Lefschetz thimbles. Section \ref{degreeshifts} consider the special case of a Lefschetz fibration obtained from a Lefschetz pencil, where useful grading-related information can be obtained.

This general construction will be applied in Section \ref{chaptermodels2} to the fibrations of Sections \ref{pencil1} and \ref{lfcp2} to produce local models for the real and complex planar projective twists.

\subsection{Building a compactly supported symplectomorphism}\label{construction} 

This section outlines the construction (after \cite{seideliterated}) of a compactly supported symplectomorphism on the total space of a class of Lefschetz fibrations (satisfying the Assumption \eqref{seidelcondition} below).
This is the general recipe that we use in building the models of the real and complex projective twists (in Sections \ref{rp2}, \ref{cp2twist} respectively) starting from the exact Lefschetz fibrations $\pi\colon E_{\APtwo} \to \C$, $\A \in \{ \R, \C \}$.

Let $\pi\colon E^{2n+2} \to \C$ be an exact Lefschetz fibration with smooth fibre a Liouville domain $(M^{2n}, \omega)$, and denote the set of critical values by $\Crit v (\pi)=\{ w_0, \dots , w_r \} \subset D_R$, where $D_R\subset \C$ is a disc of radius $R>0$. Fix a base point $z_{*} \in \C$, with $\re(z_*)\gg R$, and use it to fix a representative of the symplectomorphism class of the smooth fibre $\pi^{-1}(z_*)\cong M$.

Fix a basis of vanishing paths $(\gamma_0, \dots ,\gamma_r)$, the corresponding set of framed Lagrangian vanishing cycles $(V_0, \dots , V_r)$ in $M$ (following the conventions of Section \ref{lfconventions}) and denote by $\phi \in \Symp_{ct}(M)$ the total monodromy.

\begin{asspt} \label{seidelcondition}
	There is a symplectomorphism $\phi' \in \Symp_{ct}(M)$ with the following property. There is a symplectic isotopy $\phi\simeq \phi'$ such that
	
	\begin{align}\label{isotopy}
	\forall i=0, \dots , r, \ \     \phi' (V_i)=V_i \ \  \text{ and $\phi'|_{V_i}$ is homotopic to an element of $\OO(n+1)$.}
	\end{align}
\end{asspt}

To simplify the notation, in what follows the perturbed map will also be denoted by $\phi$.

Note that the second condition in \eqref{isotopy} ensures that for all $i=0, \dots , r$, the homotopy class of the framing of the vanishing cycle $V_i$ and that of the framing of its image under $\phi$ coincide.

\begin{lemma}\label{symptotspace}
	Under Assumption \eqref{seidelcondition}, $\phi \in \Symp_{ct}(M)$ can be extended to an element of $\Symp(E)$ which acts fibrewise as $\phi$.
\end{lemma}
\proof 
For $i=0, \dots, r$, let $f_i \colon S^n\to V_i$ be the framings of the vanishing cycles of the Lefschetz fibration $\pi\colon E\to \C$.
Then the data $(M, \phi \circ f\colon S^n\to V_i)$ define a new Lefschetz fibration $\pi' \colon E' \to \C$, which is obtained by pulling back all data in the construction of $E$ (as in \cite[(16e)]{seidelbook}) by $\phi$. This construction yields a symplectomorphism $q\colon E \to E'$. Since the framings of the vanishing cycles of $E$ and $E'$ are in the same homotopy class, there is a $1$-parameter family of exact Lefschetz fibrations $(E_t, \Omega_{E_t})_{t\in [0,1]}$ that interpolates between $E$ ($t=0$) and $E'$ ($t=1$) (\cite[(16e)]{seidelbook}).
By exactness, Moser's argument yields another symplectomorphism $\tilde{q} \colon E\to E'$. The map $\tilde{q}^{-1}\circ q \colon E\to E$ is then a symplectomorphism of the total space which by construction acts (up to isotopy) fibrewise as $\phi$. \qed

Lemma \ref{symptotspace} applies to the total monodromy of $E$ but also to its inverse, $\phi^{-1}\in \Symp_{ct}(M)$.

Let $D_{R-2\eps} \subset D_{R-\eps} \subset \C$ be two discs centered at the origin, of radii $R-2\eps$ and $R- \eps$ respectively, such that $\Crit v(\pi)\subset D_{R-2\eps}.$
Consider an anticlockwise rotational vector field on the base, supported on the annulus $A_{2\eps}:=D_R \setminus D_{R-2\eps}$, defined as follows. Let $\psi \colon \R^+ \to \R^+$ be a smooth function such that 
\begin{equation}\label{vfbase}
\begin{array}{ll}
\psi (r)=0	&r \leq R-2\eps \\
\psi'(r)=1& r > R-\eps.
\end{array}
\end{equation}

Define a Hamiltonian function $H\colon \C \to [0, \infty)$, $H(z)=\psi(|z|)$ and let $(b_t)$ be its associated Hamiltonian flow. The $2\pi$-flow $(b_{2\pi})$ defines a Dehn twist in $A_{2\eps}$, which acts as the identity on $D_{R-2\eps}$, as an anticlockwise $2\pi$ rotation on $\C \setminus int(D_R)$ and interpolates between the two on $A_{2\eps}$. 

For every $t\in [0, 2\pi]$, $(b_t)$ can be lifted via parallel transport to a family of symplectomorphisms $(\Phi_t)$ of the total space (see also \cite[Observation 6.4]{maydanskiy}). The element $\Phi_{2\pi}$ covers the base twist $b_{2\pi}$, and therefore is fibre preserving over $\C \setminus int(A_{2\eps})$, namely over $D_{R-2\eps}$ and over $\C \setminus int(D_R)$. In particular:

\begin{itemize}
	\item For $z\in D_{R-2\eps}$ the map acts fibrewise as $\Phi_{2\pi}|_{\pi^{-1}(z)}=Id$,
	\item By the Picard--Lefschetz theorem, $\Phi_{2\pi}|_{\pi^{-1}(z_*)}\simeq \phi$ in $\Symp_{ct}(M)$, since $z_*\in \C \setminus int(D_R)$.
\end{itemize} 


\begin{figure}[h]
	\centering
\def\svgwidth{160pt}\includegraphics[width=6cm]{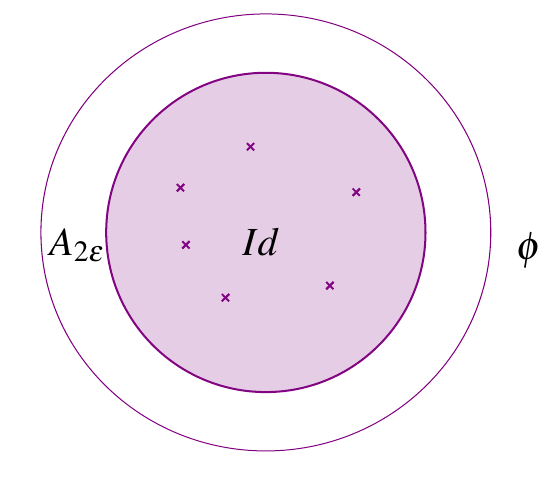}
	\qquad
\def\svgwidth{160pt}\includegraphics[width=6cm]{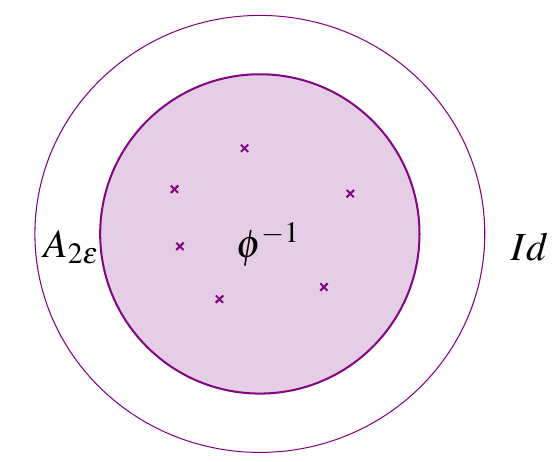}
	\caption{The fibrewise action of $\Phi_{2\pi}$ (left) and that of $\widetilde{\varphi}$ (right) over the disc $D_{R-2\eps}$ and over the complement $\C\setminus D_R$.}%
\end{figure}

By Lemma \ref{symptotspace}, $\phi^{-1} \in \Symp_{ct}(M)$ can be extended to a symplectomorphism of the total space $\tilde{\phi}^{-1} \in \Symp(E)$. Let $\tilde{\varphi}:=\tilde{\phi}^{-1} \circ \Phi_{2\pi}$. This is a symplectomorphism $\tilde{\varphi} \in \Symp(E)$ which defines another lift of the base twist $b_{2\pi}$, but is supported over the compact region $D_{R}$. 
So $\tilde{\varphi}$ is compactly supported in the horizontal direction, however not so in the vertical direction. Namely, $\tilde{\varphi}$ acts non-trivially on $\pi^{-1}(A_{2\eps}) \cap \partial^hE$.

In what follows we adjust the map by finding an isotopy $\tilde{\varphi} \simeq \varphi$ to a compactly supported symplectomorphim $\varphi \in \Symp_{ct}(E)$.

The horizontal boundary $\partial ^h E$ can be trivialised as in \eqref{horizontalboundary}: an open neighbourhoood of $\partial^h E$ is isomorphic to $U^{\partial}\cong\C \times M^{out} \subset \C \times M$, where $M^{out}\subset M$ is an open neighbourhood of $\partial M$.

The above trivialisation induces a decomposition $\widetilde{\varphi}|_{U^{\partial}} \cong \tau_{\C}\times Id_M$, where $\tau_{\C}$ is the Dehn twist in a circle in $\C$. The latter admits a symplectic isotopy with the identity, $(\alpha_t)_{t\in [0,1]}$, $\alpha_0=\tau_{\C}, \alpha_1=Id_{\C}$.
Since the symplectic form is product-like in the region $U^h$, there is a symplectic isotopy $(\alpha_t \times Id)_{t\in [0,1]}$ connecting $\tau_{\C} \times Id_M$ to $Id_{\C}\times Id_M$, supported in $U^{\partial}$ and fixing $\partial U^{\partial}$.
This isotopy does not interfere with the support of $\tilde{\varphi}$ in the compact region, so after applying the isotopy, we obtain a compactly supported symplectomorphism $\varphi\in \Symp_{ct}(E)$ whose support is $U^c:=D_R \times (M\setminus M^{out})$, and $\varphi|_{U^c}=\tilde{\varphi}|_{U^c}$.

\subsection{Floer cohomology computations}\label{hfcomputationsection}

For a given Lefschetz fibration $\pi\colon E\to \C$ admitting a compactly supported symplectomorphism $\varphi \in \Symp_{ct}(E)$ as in the previous section, the question arises on whether $\varphi$ is a non-trivial element in the mapping class group of $E$, and if so, how to study its properties. 

One possibility is to compute the Floer theoretical action of $\varphi$. This section focuses on the computation of Floer cohomology groups $\HF^*(\varphi^k(\Delta_{\alpha}), \Delta_{\beta};\Zmod)$, where $\Delta_{\alpha}, \Delta_{\beta}\subset E$ are (disjoint) Lefschetz thimbles of $\pi$. The image $\varphi^k(\Delta_{\alpha})$ is another thimble, which intersects $\Delta_{\beta}$ over $k$ regular values $z_i \in \C\setminus \Crit v(\pi)$, $i=0, \dots , k-1$, and the components of the intersection are given by $\varphi^k(\Delta_{\alpha}) \cap \Delta_{\beta}=\bigcup_{i=0, \dots , k-1} \{ V_{z_i,\alpha}\cap V_{z_i,\beta} \}$, where $V_{z_i,\alpha}:= \varphi^k(\Delta_{\alpha})\cap \pi^{-1}(z_i)$, $V_{z_i,\beta}:= \Delta_{\beta} \cap \pi^{-1}(z_i)$.

The main computation tool is a spectral sequence adapted from \cite{maysei} (Section \ref{hfthimbles}), whose first page is given by $\bigoplus_{i=0}^{k-1}\HF^*(V_{z_i,\alpha}, V_{z_i,\beta}; \Zmod)$, and which converges to $\HF(\varphi^k(\Delta_{\alpha}), \Delta_{\beta}; \Zmod)$.

The computations take place in Sections \ref{flcomplex}, \ref{hfthimbles} and \ref{degreeshifts}. Before that, we set up, in Section \ref{hfbasics}, the conventions for ungraded Floer cohomology of Lagrangians in the total space and in the smooth fibre of a Lefschetz fibration.

\subsubsection{Floer cohomology conventions}\label{hfbasics}
In what follows the ground field for all the Floer cohomology groups will always be assumed to be $\Zmod$; this avoids the use of spin structures. 

Let $\pi\colon E\to \C$ be an exact Lefschetz fibration. Let $\J(E, \pi, j)$ be the set of almost complex structures compatible with $\pi$ in the sense of Definition \ref{compatiblej}. Such a choice of almost complex structure is not generic, but one can always choose a generic element $J_E\in \J(E, \pi, j)$ (as in \cite[Section 2.1]{seideles}) making the elements of the moduli spaces below regular. Fix such a $J_E \in \J(E,\pi,j)$.

Let $\mathcal{L}_0, \mathcal{L}_1 \subset E$ be two Lagrangian thimbles with $h(\mathcal{L}_0)>h(\L_1)$, and define $\mathscr{C}(\L_0,\L_1)$ to be the set of intersection points of the pair $(\L_0, \L_1)$ after a suitable Hamiltonian perturbation to make the intersection transverse.

Given $\xi_{\pm} \in \mathscr{C}(\L_0,\L_1)$, let $u\colon \R \times [0,1]\to E$ be a solution to the (possibly perturbed) Floer equation with the properties \begin{equation}\label{floersolution}
\begin{split}
\partial_s u(s,t)+J_E(u) \partial_t u(s,t)=0, \\
\forall s\in \R: \ \ u(s,0)\in \L_0, \ u(s,1)\in \L_1\\
\lim\limits_{s \to - \infty}u(s, t)=\xi_{-}, \   \lim\limits_{s \to + \infty}u(s, t)=\xi_{+}. \\
\end{split}
\end{equation}

Let $\mathcal{M}(\xi_{-}, \xi_{+}, [u]; J_E)$ be the moduli space of unparametrised $J_E$-holomorphic curves in the class $[u]$, satisfying \eqref{floersolution}. Then $\mathcal{M}(\xi_{-}, \xi_{+}, [u]; J_E)$ is a smooth manifold of dimension $\mu([u])-1$ (\cite{oh}).

The Floer complex $\CF(\L_0, \L_1)$ is generated, as a $\Zmod$-vector space, by the elements of $\mathscr{C}(\L_0,\L_1)$, and the Floer differential of an element $\xi_{+} \in \mathscr{C}(\L_0,\L_1)$ is given by \begin{align}\label{differential}
\partial \xi_{+}= \sum_{\substack{	\xi_{-} \in \mathscr{C}(\L_0, \L_1) \\ \mu([u])=1}} (\sharp_{\Zmod}\mathcal{M}(\xi_{-}, \xi_{+}, [u]; J_E) )\langle \xi_{-} \rangle.
\end{align}

The moduli spaces of curves $\mathcal{M}(\xi_{-}, \xi_{+}, [u]; J_E)$ are compact (we explain this below), so that the above expression is well-defined.

We show that the elements of these moduli spaces satisfy compactness properties in both ``vertical'' and ``horizontal'' directions (see also \cite[Section 6]{maydanskiy}).

For compactness in the fibre (vertical) direction, let $U^{\partial} \subset E$ be an open neighbourhoood of $\partial^h E$. As before, there is an isomorphism $U^{\partial}\cong\C \times M^{out} \subset \C \times M$ for an open neighbourhood $M^{out}\subset M$ of $\partial M$, under which both the Liouville form and the almost complex structure split in a product-like fashion.
Assume a pseudoholomorphic curve $u$ as in \eqref{floersolution} was to enter the neighbourhood $\partial^hU$. Then the projection of $u$ to the second factor would be a pseudoholomorphic curve in $M$, with interior points in $M^{out}$ and boundary conditions in the compact part away from $M^{out}$, a contradiction to the maximum principle.

On the other hand, since $\pi$ is $(J_E,j_{\C})$-holomorphic, the projection of a non-trivial $J_E$-holomorphic curve $u$ to the base $\C$ is a holomorphic strip with boundary condition on $(\gamma_0(\R^+), \gamma_1(\R^+))$, and by the maximum principle the strip cannot escape a compact neighbourhood in $\C$ containing the intersection points of the paths. Therefore, $u$ cannot go arbitrarily far in the horizontal direction.

By exactness, the action functional only depends on the endpoints $\xi_{\pm}$ and therefore gives a common upper bound for the energy of the curves in these moduli spaces. So by the Gromov compactness theorem (for a statement of this theorem see \cite{mcduffsal2}) together with $\dim(\mathcal{M}(\xi_{-}, \xi_{+}, [u]; J_E))=0$, the only remaining issue are bubble phenomena. But these cannot occur by exactness.

Now consider two Lagrangian thimbles $\mathcal{L}_0, \mathcal{L}_1 \subset E$ with associated vanishing paths $\gamma_0, \gamma_1\colon \R^+\to \C$ and assume $h(\L_0)\leq h(\L_1)$. For a point in the base $z \in \C$, let $x:=\re(z)$, $y:=\Im(z)$.

\begin{definition}
	For $\epsilon > h(\L_1)-h(\L_0)$, define a map $H_{\epsilon} \in C^{\infty}(\C,\R)$ satisfying
	
	\begin{equation}\label{admissible}
	\begin{array}{ll}
	H_{\epsilon} (x)=0	&x<c \\
	H_{\epsilon}'(x)=\epsilon & x>c+K.
	\end{array}
	\end{equation}
	for $c>R$ and $K >0$. 	
	Define the Hamiltonian vector field $Y_{\epsilon}=H'_{\epsilon}(x)\partial_y$, call its time-1 flow $\chi_{\epsilon}$ and let $\widetilde{\chi_{\epsilon}}$ be a lift of $\chi_{\epsilon}$ to the total space. If $\mathcal{L}_0$ is a Lefschetz thimble, then $\widetilde{\chi_{\epsilon}}(\mathcal{L}_0)$ is a another Lefschetz thimble isotopic to $\mathcal{L}_0$ and $\pi(\widetilde{\chi_{\epsilon}}(\mathcal{L}_0))=\chi_{\epsilon}(\gamma_0)$.
\end{definition}
In the case $h(\L_0)\leq h(\L_1)$, the Floer complex $\CF(\L_0, \L_1)$ is defined as above, but considering the modified pair $(\chi_{\eps}(\L_0), \L_1)$.

For Lagrangian vanishing cycles, the situation is fairly standard. For a pair $(V_0, V_1) \subset (M,\omega)$ of closed exact Lagrangian submanifolds of the fibre, the Floer cohomology groups $\HF(V_0, V_1;\Zmod)$ are well defined (see for example \cite[8,9]{seidelbook}).

\subsubsection{The Floer complex $\CF(\varphi^k(\Delta_a), \Delta _c)$}\label{flcomplex}

Let $\pi\colon E \to \C$ be an exact Lefschetz fibration with smooth fibre a Liouville domain $(M,d\lambda)$, $r+1$ critical fibres and associated critical values $\Crit v=\{ w_0, \dots , w_{r} \} \subset \C$. Let $z_* \in \C$ be a base-point with $\re(z_*)\gg 0$, fixing the smooth fibre $ M:=\pi^{-1}(z_*)$. Fix a distinguished basis of vanishing paths $(\gamma_0, \dots , \gamma_{r})$ for the critical values, and the associated basis of Lagrangian vanishing thimbles $(\Delta_0, \dots , \Delta_{r})$. Call $(V_0,\dots , V_{r}) \subset M$ the resulting basis of vanishing cycles in the fixed fibre. Choose two elements of the basis $\Delta_{\alpha}, \Delta_{\beta}$, $\alpha, \beta \in \{ 0, \dots , r\}$ with the property that $h(\Delta_{\alpha})>h( \Delta_{\beta})$. 

For the entire section, assume that $\phi:= \tau_{V_0}\cdots \tau_{V_r}$ satisfies Assumption \eqref{seidelcondition}, which ensures the existence of a symplectomorphism $\varphi \in \Symp_{ct}(E)$ as constructed in Section \ref{construction}. Note that we don't use Assumption \eqref{anticanonical1} until Section \ref{degreeshifts}.

Let $R>0$, and $D_{R-2\eps} \subset D_R \subset \C$ two discs in the base (of radii $R-2\eps$ and $R$ respectively, for $\eps>0$) containing the critical values. Let $b_{2\pi}$ be a Dehn twist in the annulus $A_{2\eps}=D_R \setminus D_{R-2\eps}$ as defined in Section \ref{construction}. Note that the image $ \varphi(\Delta_{\gamma})$ of the vanishing thimble associated to a vanishing path $\gamma$ is Lagrangian isotopic to the vanishing thimble $\Delta_{b_{2\pi}( \gamma)}$ associated to the twisted path $b_{2\pi}( \gamma)$.

For $k\in \Z$, set $\gamma^k_{\alpha}:=b_{2\pi}^{k}(\gamma_{\alpha})$ and $\Delta^k_{\alpha}:=\Delta_{\gamma^k_{\alpha}}=\varphi^k(\Delta_{\alpha})$.

Call the set of intersection points $I:=\{\gamma_{\alpha}^k\cap \gamma_{\beta} \} = \{  z_{0}, \dots ,  z_{k-1} \} \subset A_{2\eps}$, where $z_0$ is the innermost and $z_{k-1}$ the outermost intersection point.


\begin{figure}[h]
	\centering
\def\svgwidth{400pt}\includegraphics[width=9cm]{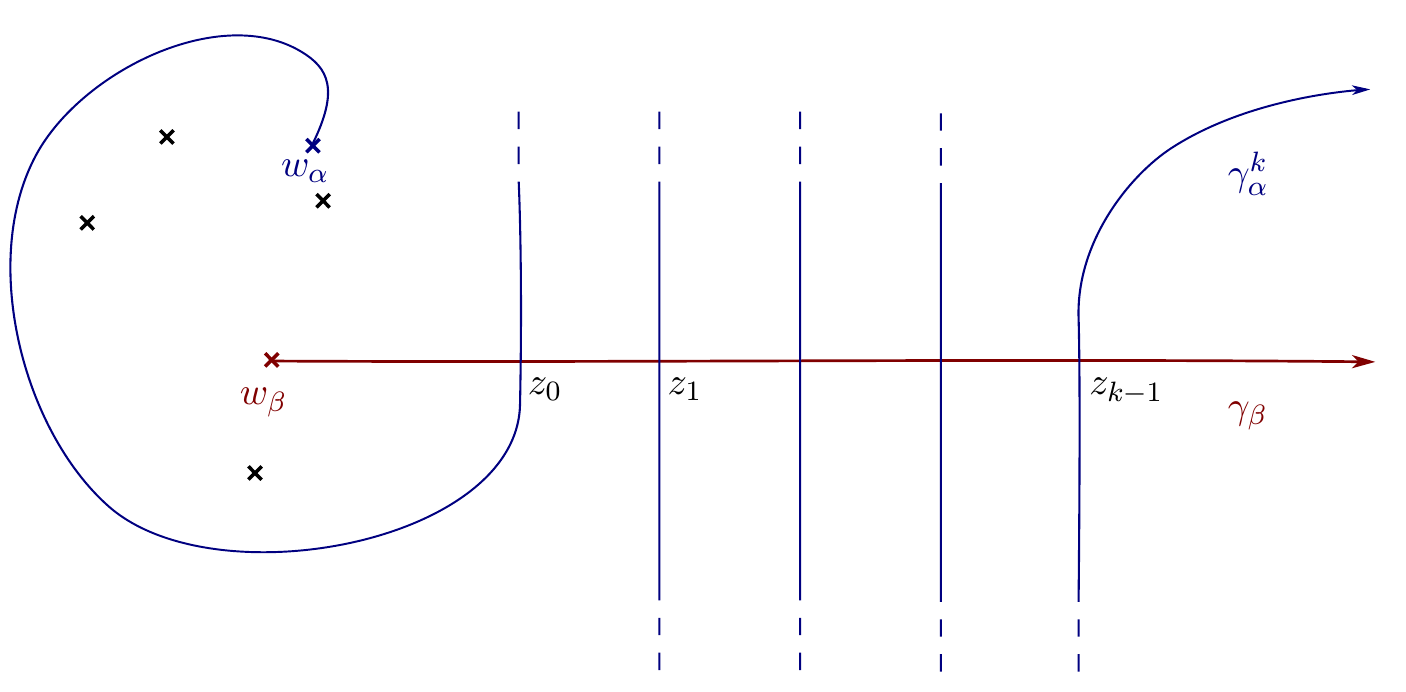}
	\caption{The intersection pattern of the pair of paths $(\gamma_{\alpha}^k, \gamma_{\beta})$}%
\end{figure}

Each component of the intersection locus $\{ \Delta_{\alpha}^k \cap \Delta_{\beta} \}$ lying over $z_j \in I$ is determined by the intersection of the pair $(V_{z_{j},\alpha}, V_{z_j,\beta}) \subset E_{z_j}:=\pi^{-1}(z_j)$, where $V_{z_j,{\alpha}}:=E_{z_j} \cap \Delta_{\alpha}^k, \ V_{z_j,{\beta}}:=E_{z_j} \cap \Delta_{\beta}.$

Moreover, by the discussion of Section \ref{gradingshift}, the global monodromy $\phi\in \Symp_{ct}(M)$ of $\pi\colon E\to \C$ preserves the vanishing cycles up to a shift in their gradings. As a consequence, the intersection $\Delta_{\alpha}^k \cap \Delta_{\beta} \cap E_{z_j}$ is the same at any $z_j\in I$. Therefore, we can apply identical (local, compactly supported) Hamiltonian perturbations in each fibre $E_{z_j}$ in order to turn all intersections into a configuration of points $\{\xi_{z_j,1}, \dots , \xi_{z_j,\ell} \}$ for a fixed $\ell \in \N$, and $j\in \{ 0, \dots , k-1\}$. Fix a generic element $J_E \in \J(E,\pi, j)$ as in Section \ref{hfbasics}. For all $j \in \{ 0, \dots , k-1 \}$, the (ungraded) Floer cohomologies $\HF(\Delta_{\alpha}^k, \Delta_{\beta}; \Zmod)$ and $\HF(V_{z_j,{\alpha}}, V_{z_j,{\beta}};\Zmod)$ are well defined by the previous section (\ref{hfbasics}). 

Recall that $\Delta_{\alpha}$ and $\Delta_{\beta}$ are disjoint and $h(\Delta_{\alpha})>h(\Delta_{\beta})$ so $\HF(\Delta_{\alpha}^0, \Delta_{\beta};\Zmod)=\HF(\Delta_{\alpha},\Delta_{\beta};\Zmod)=0$ (by the conventions of Section \ref{hfbasics}).

\begin{lemma}\label{firstpowerhf}
	The first power of $\varphi$ satisfies
	\begin{align}
	\HF(\Delta_{\alpha}^1, \Delta_{\beta})=	\HF(\varphi(\Delta_{\alpha}), \Delta_{\beta})\cong \HF(V_{z_0,\alpha}, V_{z_0,\beta}).
	\end{align}
\end{lemma}
\proof If we consider the first power of the map $\varphi$, the vanishing paths associated to the two thimbles intersect in a single point $z_0$. Then there can be no ``horizontal'' Floer differential. Namely, all generators are contained in $\pi^{-1}(z_0)$, so any pseudoholomorphic curve $u\in \mathcal{M}(\Delta_{\alpha}^1, \Delta_{\beta})$ connecting these generators must be confined in that fibre, since $\pi$ is $(J_E,j_{\C})$-holomorphic and the open mapping theorem applies to the image $\pi(u)$.
Therefore, we have $\HF(\varphi(\Delta_{\alpha}),\Delta_{\beta}) \cong \HF(V_{z_0,\alpha},V_{z_0,\beta})$.  \endproof



\begin{figure}[h]
	\centering
\def\svgwidth{200pt}\includegraphics[width=9cm]{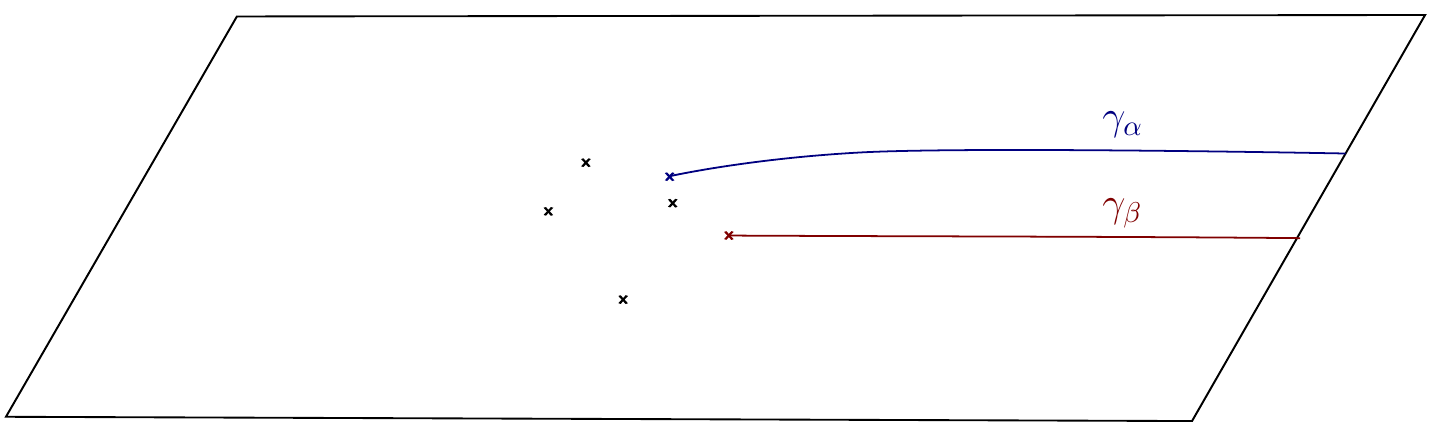}
	\qquad
\def\svgwidth{200pt}\includegraphics[width=9cm]{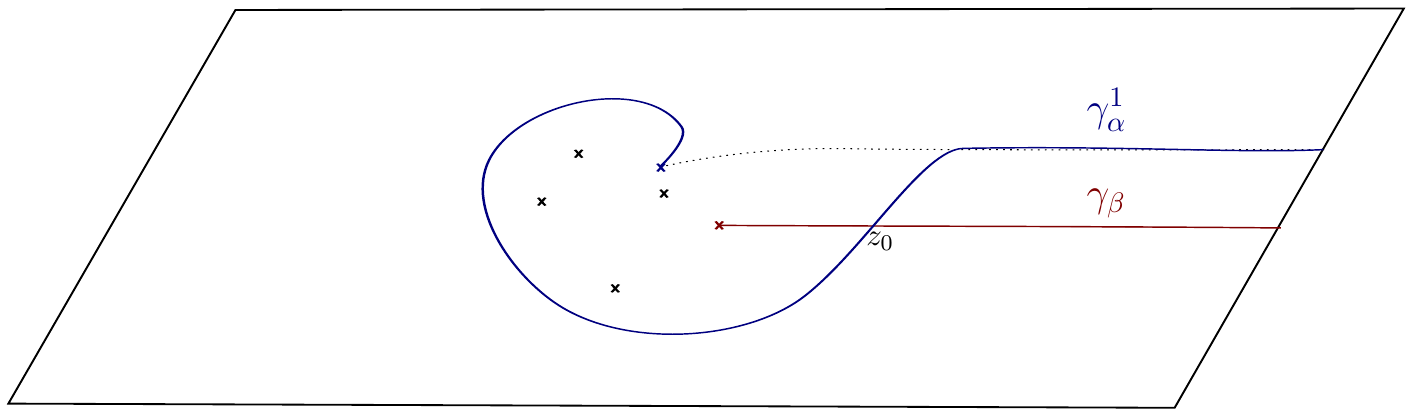}
	\qquad
\def\svgwidth{200pt}\includegraphics[width=9cm]{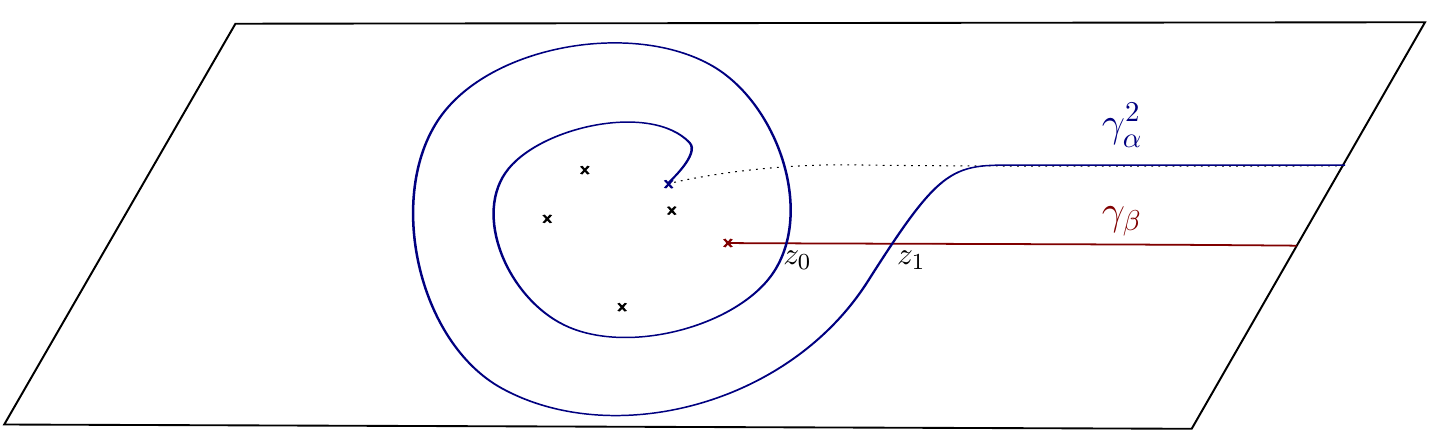}
	\vspace{1cm}
	\caption{Top: the initial vanishing paths $(\gamma_{\alpha},\gamma_{\beta})$, middle: the first power $(\gamma^1_{\alpha}, \gamma_{\beta})$, bottom: the second power $(\gamma_{\alpha}^2, \gamma_{\beta})$.}
\end{figure}

To determine the groups $\HF(\varphi^k(\Delta_{\alpha}), \Delta_{\beta})$ for $k>1$, we will need to introduce additional computational tools. In the next section, we opt for a spectral sequence that was expressly developed to compute the Floer cohomology of Lefschetz thimbles in \cite{maysei}.

\subsubsection{A spectral sequence}\label{hfthimbles}\label{arcs}
In this section we complete the computation of the Floer cohomology groups $\HF^*(\varphi^k(\Delta_{\alpha}), \Delta_{\beta};\Zmod)$ of the Lefschetz thimbles $\varphi^k(\Delta_{\alpha})$ and $\Delta_{\beta}$ with the help of a spectral sequence originally defined in \cite{maysei}.

Let $\pi\colon E \to \C$ be the Lefschetz fibration we have considered in the previous subsection, and assign gradings for the the Lefschetz thimbles and the vanishing cycles.

Consider the holomorphic strips in the base formed by the pair of arcs $(\gamma_{\alpha}^k, \gamma_{\beta})$, and connecting the intersection points of $I=\{\gamma_{\alpha}^k\cap \gamma_{\beta} \}=\{ z_0, \dots , z_{k-1} \}$. We assign indices to these points as follows. Fix $\ind(z_0)=0$ and for $j=1, \dots , k-1$, set  \begin{align}\label{basegr}
\ind(z_{j})=j \cdot (-2).
\end{align}
For any $z \in I$, $\ind(z)$ is an even negative value, $-2(k-1)\leq\ind(z)\leq 0$. 
Geometrically, $\ind(z_j)$ corresponds to the Maslov index of the holomorphic strip with endpoints $z_{j-1}, z_j$ and boundary condition on $(\gamma_{\alpha}^k, \gamma_{\beta})$.

We define an order on the set $I$ as follows: $z_-, z_+ \in I$ satisfy $z_+>z_-$ if there is a sequence $v_1,\dots , v_m \colon \R\times [0,1] \to \C$ of holomorphic curves with bounded energy, of the form 
\begin{align*}
v_i(\R\times \{ 0 \} )\subset \gamma_{\alpha}^k(\R^+), \ v_i(\R\times \{ 1 \} )\subset \gamma_{\beta}(\R^+),  \\
\lim_{s\to +\infty}v_i(s, \cdot)=\lim\limits_{s \to - \infty}v_{i+1}(s,\cdot), \\
\   \lim_{s\to -\infty}v_1(s, \cdot)=z_-,      \lim_{s\to +\infty}v_m(s, \cdot)=z_+, 
\end{align*}

By the open mapping theorem, $z_{0}>  \cdots >z_{k-1}$ and note that for $z_-, z_+ \in I$, $\ind(z_+)<\ind(z_-) \iff z_+>z_-$.

\begin{prop}\label{spectral}
	There is a cohomological spectral sequence with bigraded differentials $d_{c} \colon E_{c}^{p,q}\to E_c^{p+c,q-c+1}$, converging to $HF^*(\Delta_{\alpha}^k,\Delta_{\beta};\Zmod)$. The starting page is generated by 
	
	\begin{equation}\label{spectralsequence}
	E_1^{pq} = \left\{
	\begin{array}{ll}
	\HF^{p+q}(V_{z_{(-p)},\alpha}, V_{z_{(-p)},\beta};\Zmod)	 & -(k-1) \leq p \leq 0 \\
	0 & \emph{otherwise.}
	\end{array}
	\right.
	\end{equation}

\end{prop}

\proof From now onwards, we will omit the coefficient field $\Zmod$ in the notation. The spectral sequence is obtained as a special instance of \cite[Proposition 4.1]{maysei}, whose proof is adapted below. Let $(J_t)_{t\in [0,1]}$ be a family of almost complex structures on $E$ such that for each $t \in [0,1]$, the tuple $(J_t, j)$ is compatible with the fibration $\pi$ in the sense of Definition \ref{compatiblej}. Then, for a pair of intersection points $(\xi_-,\xi_+) \in \Delta_{\alpha}^k\cap \Delta_{\beta}$ we can define a relation $\xi_+>\xi_-$ if there is a sequence $u_1, \dots , u_m \colon \R \times [0,1]\to E$ of pseudoholomorphic curves of bounded energy such that for all $i=1, \dots , m$
\begin{align*}
\partial_su_i(s,t)+J_t\partial_tu_i(s,t)=0, \\
u_i(\R\times \{ 0 \} )\subset \Delta_{\alpha}^k, \ u_i(\R\times \{ 1 \} )\subset \Delta_{\beta}, \\
\lim_{s\to +\infty}u_i(s, \cdot)=\lim\limits_{s \to - \infty}u_{i+1}(s,\cdot), \\
\ \lim_{s\to -\infty}u_1(s, \cdot)=\xi_-,  \lim_{s\to +\infty}u_m(s, \cdot)=\xi_+.
\end{align*}

\noindent If $\xi_+>\xi_-$, then $\xi_{\pm}$ must either lie in the same fibre or their projection satisfy $\pi(\xi_+)>\pi(\xi_-)$ in the ordering defined for intersection points on the base. This means that $\xi_+>\xi_-$ implies $\ind(\pi(\xi_+)) \leq \ind(\pi(\xi_-))$.

As noted in Section \ref{lfconventions}, the choice of almost complex structures is not generic, but by applying a small perturbation to the family $(J_t)_{t\in [0,1]}$, these curves meet the usual regularity and compactness requirements.

Let $I_p:= \{ z \in I, \ \ind(z)\geq 2p \} \subset I$, so that $I_p$ is non-empty for $-(k-1)\leq p \leq 0$.
The complex $\CF(\Delta_{\alpha}^k, \Delta_{\beta})$ admits a filtration $F^*$, where each term $F^{p}$ is generated by the intersection points in $\pi^{-1}(I_{p})$. By definition of the order, the Floer differential preserves $F^{p}$ so that there is an induced differential on $F^{p}/F^{p+1}$. The latter is therefore a cochain complex, generated by the intersection points in $\pi^{-1}(z)$, for the unique element $z=z_{-p} \in I_p\setminus I_{p+1}$
This gives rise to a spectral sequence, whose first page is given by the cohomology $H^*(F^{p}/F^{{p}+1})=\HF(V_{z,\alpha}, V_{z,\beta})$. The last equality holds because regularity for $J_t$-holomorphic strips in the total space is equivalent to regularity of $J_t|_{\pi^{-1}(z)}$-holomorphic strip in the fibre obtained by restriction (see the original proof of \cite{maysei} for details). \endproof

\subsection{The special case of a Lefschetz fibration obtained by a Lefschetz pencil}\label{degreeshifts}

The spectral sequence above becomes computationally relevant when there is a concrete understanding of the grading shift that the monodromy causes on vanishing cycles. This is the case for Lefschetz fibrations that are induced by Lefschetz pencils satisfying Assumption \eqref{anticanonical1}; we now restrict to this specific situation to obtain a more enlightening version of \eqref{spectralsequence}.

We have seen in Lemma \ref{generalantern} that in the case of a Lefschetz fibration $\pi \colon E\to \C$ induced by a Lefschetz pencil, the total monodromy is isotopic to a symplectomorphism (a fibred twist) whose restriction to the vanishing cycles is the identity (because the fibred twist is supported away from the cycles). In particular, (a stronger version of) Assumption \eqref{seidelcondition} holds, and one can build a symplectomorphism $\varphi \in \Symp_{ct}(E)$ as in Section \ref{construction}. When moreover we assume the pencil satisfies Assumption \eqref{anticanonical1}, we can observe interesting patterns in the gradings of vanishing cycles and their images under iterations of $\varphi$.

We keep all the notation as in the previous subsection, and again assume the Lagrangians $\Delta_{\alpha}^k$, $\Delta_{\beta}$ have been perturbed by a Hamiltonian isotopy to make their intersection transverse (a set of $\ell$ points) over each $z_j \in I$, $j\in \{ 0, \dots, k-1\}$. This section studies the grading of intersection points in $\{\xi_{j,1}, \dots , \xi_{j,\ell} \}$ obtained after the perturbations applied to $\{\Delta^k_{\alpha} \cap \Delta_{\beta} \}$.

For $j\in \{ 0, \dots , k-1\}$ and $i \in \{ 1, \dots , \ell \}$, choose $\xi_{z_j,i} \in \Delta_{\alpha}^k \cap \Delta_{\beta} \cap E_{z_j}$ and let $\xi_{z_{j+1},i}\in \Delta_{\alpha}^k \cap \Delta_{\beta} \cap E_{z_{j+1}}$, the corresponding point in the next fibre.


\begin{figure}[h]	
	\centering
\def\svgwidth{200pt}\includegraphics[width=6cm]{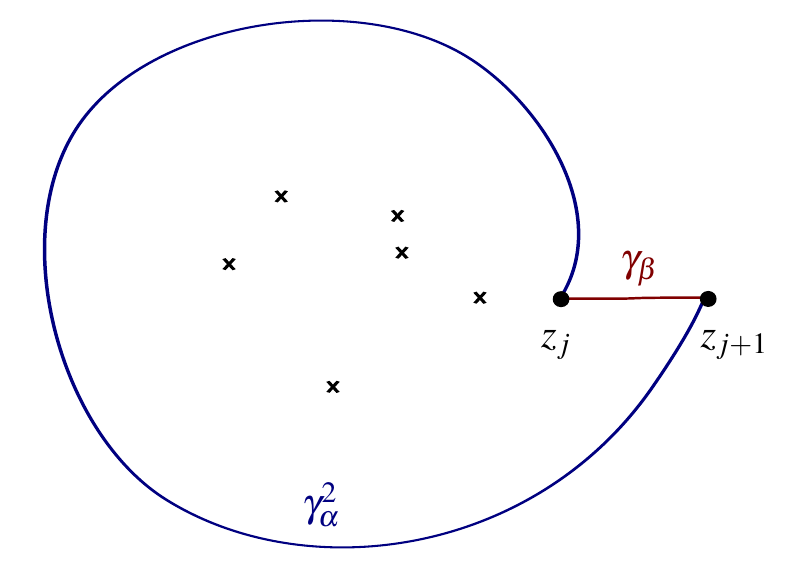}
	\caption{The projection in $\C$ of a pseudoholomorphic strip connecting $\xi_j$ to $\xi_{j+1}$.}\label{projectionstrip}
\end{figure}

\begin{lemma}\label{degreeshift}
	The degree shift is given by\begin{align}
	s_{\varphi}:=\deg(\xi_{z_{j+1}, i})- \deg(\xi_{z_j, i})=d-2.
	\end{align} 
\end{lemma}
\proof 
We compute the index of a pseudoholomorphic strip $u\colon \R\times [0,1]\to E$ with boundary conditions \begin{align*}
\forall s \in \R: \ u(s,0)\in \Delta_{\alpha}^k, \ u(s,1)\in \Delta_{\beta}, \\
\lim\limits_{s\to -\infty}u(s,t)=\xi_{z_j,i}, \ \lim\limits_{s\to +\infty} u(s,t)=\xi_{z_{j+1},i}.
\end{align*}

This will give the shift in grading between a fibre and the next, $s_{\varphi}$, which is independent of the homotopy class of $u$.

After fixing a trivialisation of $u^*|_{TM}$, the maps \begin{align*}
\delta_{\alpha}:=u^*|_{\R\times \{ 0 \}}T\Delta_{\alpha}^k, \ \delta_{\beta}:=u^*|_{\R\times \{ 1 \}}T\Delta_{\beta}
\end{align*}
can be viewed as paths $\delta_{\alpha,\beta}(s) \colon [-\infty, +\infty] \to LGr(n)$ connecting $T_{\xi_{z_j,i}}\Delta_{\alpha}^k$ to $T_{\xi_{z_{j+1},i}}\Delta_{\alpha}^k$ and $T_{\xi_{z_j,i}}\Delta_{\beta}$ to $T_{\xi_{z_{j+1},i}}\Delta_{\beta}^k$ respectively.

The Maslov index of $u$ can be computed as (\cite[Section 3]{maslov}, see also \cite[1.3]{introfuk}) the number of times (counted with multiplicities) at which $\delta_{\alpha}(s)$ and $\delta_{\beta}(s)$ are not transverse to each other. 
Since for every $x \notin \Crit(\pi)$, the tangent space $T_xE$ admits symplectic splitting \eqref{splitting} $T_xE\cong T_x^vE \oplus T^h_xE=ker(D_x\pi)\oplus \pi^*(T_x\C)$, and the Maslov index is additive under such decomposition (\cite[Theorem 2.3]{maslov}), we can break the computation of the index into horizontal (base) and vertical (fibre) direction.

In the base direction, it suffices to look at the paths $\gamma_{\alpha}^k$ and $\gamma_{\beta}$ to see that the contribution is $+2$ (see Figure \ref{projectionstrip}). In the vertical direction, the winding of $\delta_{\alpha}$ relative to $\delta_{\beta}$ is given by the change in phase of the relative quadratic form $\eta_X/dz^2$ under a $2\pi$-rotation, which, as in Lemma \ref{lemmamonodromy}, is given by $d-4$.

We therefore obtain a total degree shift of
\begin{align}\label{grshift}
s_{\varphi}:=(d-4)+2=d-2.
\end{align}
\qed

This value indicates the shift in the grading of the vanishing cycles from one intersection fibre to the next. Namely, after parallel transporting $(V_{z_j,\alpha}, V_{z_j,\beta})$ and $(V_{z_{j+1},\alpha}, V_{z_{j+1}, \beta})$ into the smooth fibre $M=\pi^{-1}(z_*)$ (to be able to compare them), we can write $(V_{z_{j+1}, \alpha}, V_{z_{j+1}, \beta})=(V_{{z_j}, \alpha}[-s_{\varphi}], V_{z_j, \beta})$ and therefore, for $j= 1, \dots , k-1$
\begin{align}
(V_{z_{j}, \alpha}, V_{z_j,\beta}) = (V_{z_0,\alpha}[-j \cdot s_{\varphi}], V_{z_j, \beta})= (V_{z_0,\alpha}[j (2-d)], V_{z_j, \beta}).
\end{align}

The formula \eqref{grshift} indicates how this shift affects the Floer cohomology of the cycles; for $j=1, \dots k-1$, we have
\begin{align}\label{cycleshift}
\HF^*(V_{z_{j},\alpha}, V_{z_{j},\beta}; \Zmod)\cong \HF^{*}(V_{z_0,\alpha}[j \cdot (2-d)], V_{z_0,\beta};\Zmod) \cong \HF^{*-j\cdot (2-d) }(V_{z_0,\alpha}, V_{z_0,\beta};\Zmod).
\end{align}

Now recall the subsets $I_p=\{ z \in I, \ \ind(z)\geq 2p \}\subset I$ used in the definition of the spectral sequence of Proposition \ref{spectral}, $-(k-1)\leq p \leq 0$. Using \eqref{cycleshift}, we can see that the nontrivial $(p,q)$-entries of the first page of the spectral sequence \eqref{spectralsequence} are given by (recall $-(k-1)\leq p \leq 0$)
\begin{equation}\label{spsequence}
E_1^{pq}=\HF^{p+q}(V_{z_{(-p)},\alpha}, V_{z_{(-p)},\beta})\cong   \HF^{p+q} \left( V_{z_0,\alpha}\left[ -p(2-d)\right], V_{z_0,\beta}\right) \cong  
\HF^{p+q-p(d-2)} ( V_{z_0,\alpha}, V_{z_0,\beta}).
\end{equation}

\section{Model projective twists}\label{chaptermodels2}

As showed in Section \ref{modelfibred}, projective twists can be identified with $S^1$-fibred twists, in a local model in which the coisotropic submanifolds (defining these twists) are given by the unit cotangent bundles $ST^*\APn$, $\A\in \{ \R, \C\}$. 
In this section, we use the constructions of Section \ref{chaptermodels} to introduce an alternative local model for real and complex projective twists in dimension two. 

Consider the exact Lefschetz fibrations $\pi\colon E_{\APtwo} \to \C$ of Sections \ref{pencil1} ($\A=\R$) and \ref{lfcp2} ($\A=\C$) respectively.
The generalised lantern relations \eqref{lantern} and \eqref{lantern2} force the total monodromy of these fibrations to preserve every vanishing cycle. This property ensures that the conditions necessary to build a compactly symplectomorphism on the total space as the one of Section \ref{chaptermodels} are met. For $\A \in \{ \R, \C \}$, that yields a compactly supported symplectomorphisms $\varphi \in \Symp_{ct}(T^*\APtwo)$ that we aim to compare with the real and complex projective twists.

In the real case, we can combine these computations with the knowledge of the mapping class group $\pi_0(\Symp_{ct}(T^*\RPtwo))$ (computed in \cite{evans}) to obtain:

\begin{theorem}
	The symplectomorphism $\varphi \in \Symp_{ct}(T^*\RPtwo)$ is isotopic to a power of the projective twist $\tau_{\RPtwo}^k$, $k\in \Z^*$.
\end{theorem}

However, in the complex case, very little is known about $\pi_0(\Symp_{ct}(T^*\CPtwo))$. As we have shown in \cite{bct}, for $n\geq 19$, there are cases in which a non-standard framing of the projective twist produces a symplectomorphism that is not Hamiltonian isotopic to the standard $\tau_{\CP}\in\Symp_{ct}(T^*\CP)$, but that doesn't say anything about $n=2$. As a consequence, we only obtain a partial result that is based on our Floer cohomological computations.

\begin{theorem}
	The symplectomorphism $\varphi \in \Symp_{ct}(T^*\CPtwo)$ is of (symplectic) infinite order.
\end{theorem}
We can only conjecture that $\varphi$ is isotopic to the complex planar projective twist. However, there is strong evidence that our constructions correspond to the projective twists. Namely, for the two examples, the current literature seems to confirm these symplectomorphisms are indeed projective twists (this is shown in Sections \ref{expectedrp2} and \ref{expectedcp2}).

\subsection{$\RPtwo$ twist}\label{rp2}
We define a non-trivial, compactly supported symplectomorphism $\varphi \in \Symp_{ct}(T^*\RPtwo)$ on the total space of the Lefschetz fibration $\pi\colon E_{\RPtwo} \to \C$ built in Section \ref{pencil1}, by applying the construction of Section \ref{construction} to this fibration.

Using the Floer theoretical computations for Lefschetz thimbles of Section \ref{hfthimbles}, we prove that $\varphi$ is isotopic to a power of the projective twist $\tau_{\RPtwo}^k \in \pi_0(\Symp_{ct}(T^*\RPtwo))$, $k\in \Z$.

\subsubsection{The real projective twist}\label{hfrp2}
Let $\pi\colon E_{\RPtwo} \to \C$ be the Lefschetz fibration of Section \ref{pencil1}, with base point $z_*\in \C$, and smooth fibre $\pi^{-1}(z_*)$ exact symplectomorphic to a $2$-sphere with four boundary components. Consider the vanishing paths $(\gamma_0, \gamma_1, \gamma_2)$, Lefschetz thimbles $(\Delta_{0}, \Delta_1, \Delta_2)$ and vanishing cycles $(V_0,V_1,V_2)$ as in Section \ref{monodromy1}.

According to the conventions of Section \ref{hfthimbles}, for $j, k \in \{ 0,1,2\}$, $j\neq k$ any two thimbles $\Delta_j, \Delta_k$ are disjoint and satisfy $ h(\Delta_j)>h(\Delta_k)$ if $j>k$. Any element in the basis of vanishing cycles satisfies Condition \eqref{seidelcondition} by Lemma \ref{condition1}, so Section \ref{construction} yields a well-defined compactly supported symplectomorphism $\varphi \in \Symp_{ct}(T^*\RPtwo)$ on the total space of the Lefschetz fibration $\pi\colon E_{\RPtwo}\to \C$. 

Let $\gamma_2^k=b_{2\pi}^k(\gamma_2)$ and $\Delta_{2}^k:=\Delta_{\gamma_2^k} \cong \varphi^k(\Delta_2)$ and $I:= \gamma_2^k(\R^+) \cap \gamma_0(\R^+)  =
\{ z_0, \dots ,z_{k-1} \}$. Over each intersection point $z \in I$, the thimbles $\Delta_2^k$ and $\Delta_0$ intersect in their vanishing cycles $V_{z,2}$ and $V_{z,0}$ respectively, which are Lagrangian circles meeting at two points (see Figure \ref{fibrecircles}) \begin{align}
V_{z,0}\cap V_{z,2}:=\{ \xi_z^-, \xi_z^+ \} \subset \pi^{-1}(z).
\end{align}
For each $z \in I$, the Lagrangians $(V_{z,0}, V_{z,2})$ trace four punctured strips between $\xi_z^-$ and $\xi_z^+$, so there cannot be any non-trivial differentials and hence the pairs $(\xi_z^-, \xi_z^+)$ generate the Floer cohomology $\HF^*(V_{z,2}, V_{z,0};\Zmod) \cong \Zmod \oplus \Zmod$. Together with Proposition \ref{firstpowerhf}, we obtain:
\begin{prop}\label{proprp2}
	The symplectomorphism $\varphi \in \pi_0(\Symp_{ct}(T^*\RPtwo))$ is a non-trivial element of the symplectic mapping class group. \qed
\end{prop}


Endow the Lagrangian vanishing cycles in the smooth fibres and the Lagrangian thimbles with a $\Z$-grading as in Section \ref{vcycles}.
Set $\ind(z_0)=0$ and $\ind(z_{j})=-2j$ for $j=1, \dots , k-1$ as in \ref{arcs}.

Let $(\xi_0^-, \xi_0^+)$ be the pair of generators of $\CF^*(V_{z_0,0}, V_{z_0, 2})$. 

\begin{lemma}\label{lemmadegreegens}
	For every $z\in I$, the degree of both generators $(\xi_z^-,\xi_z^+)$ is the same, and $(\xi_0^-,\xi_0^+)$ can be set in degree 0.
	
	Therefore, we fix \begin{align}
	\HF^*(V_{z_0,2}, V_{z_0,0}; \Zmod) \cong \Zmod \oplus \Zmod \text{ for }  *=0.
	\end{align}
\end{lemma}
\proof 
The relation $\mathcal{K}_{\CPtwo}^{\otimes 2}\cong \mathcal{O}(2)^{\otimes -3}$ for the canonical bundle $\mathcal{K}_{\CPtwo}\cong \mathcal{O}(-3)$
induces a quadratic volume form $\eta_X$ with four poles at the base locus points, no zeroes. As discussed in Section \ref{gradingshift}, in turn this yields a relative quadratic form on the smooth fibres, that is holomorphic except at the four base points where it has poles of single order (this can be seen with a computation in local coordinates). 

Since $\mathcal{K}_{\CPone}^{\otimes 2}\cong \O_{\CPone}(-4)$, $\mathcal{K}_{\CPone}^{\otimes2}\otimes \O(D)=\O$ for a divisor $D\subset \CPone$ of four points, the space of sections of $\mathcal{K}^2_{\CPone}$ with four single poles is the space of sections of the trivial bundle. Therefore, on a smooth fibre there is a unique (up to scale) quadratic form as above. Fix one such form, which will be of the shape $\eta_M:=dz^2\cdot r(z)$, for a rational function $r(z)$ with four simple poles.

Consider now the double cover of the sphere, branched at four points\begin{equation}\label{doublecover}
\begin{split}
p\colon  T^2&\lra S^2 
\end{split}
\end{equation} 
obtained by quotienting by the covering involution $w\mapsto -w$ (where $w$ is a coordinate on $\C/\Z^2$).

\begin{figure}[htb]
	\centering
\def\svgwidth{200pt}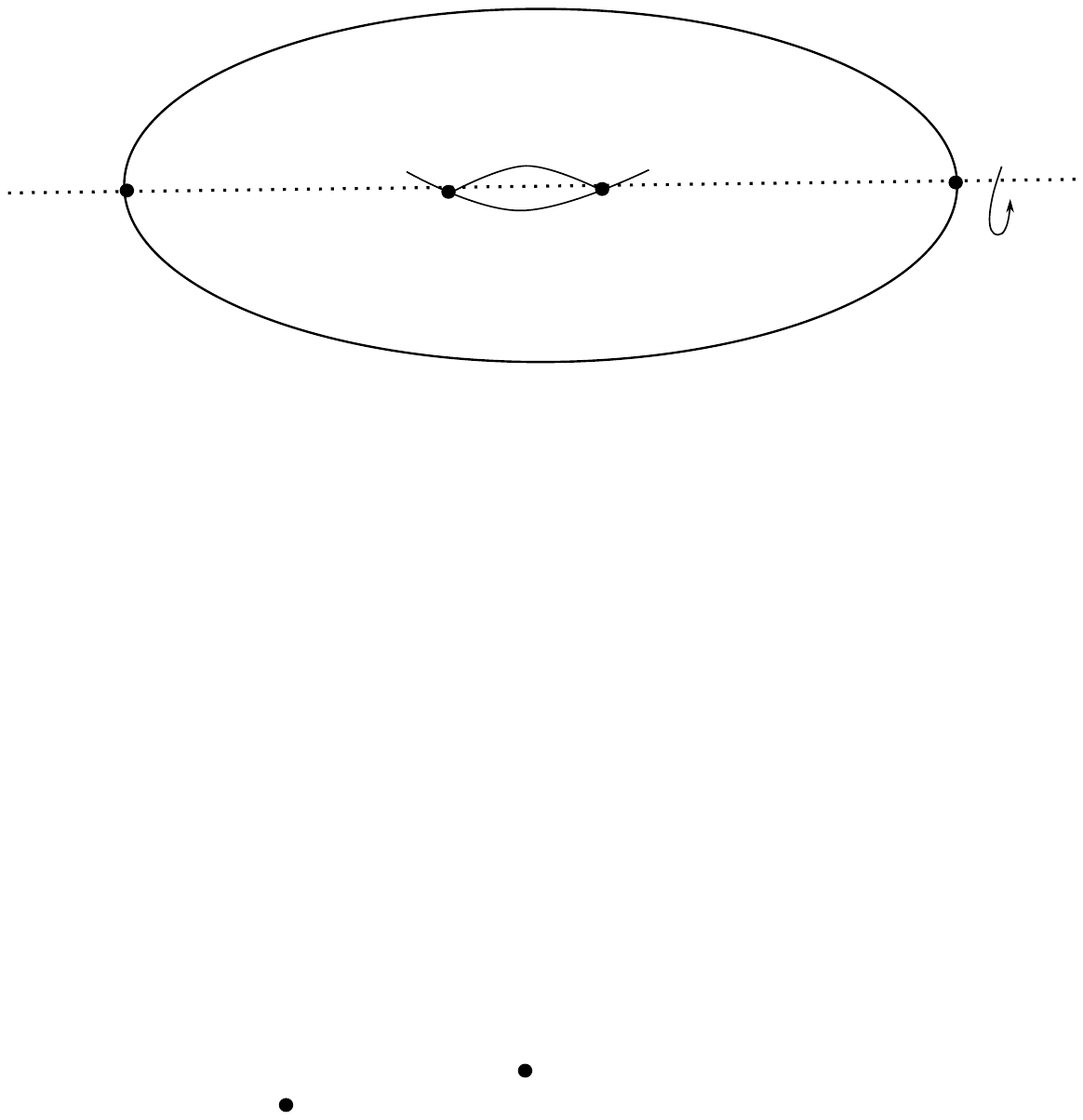
	\caption{The branched double cover of the sphere.}
	\label{doublecoverfigure}
\end{figure}

The pull-back $\eta_M$ under \eqref{doublecover} is a holomorphic form (the pullback of $dz^2/z$ is holomorphic), so it must be a multiple of the square of the standard form on the torus, $dw^2$. This form is translation invariant.

The vanishing cycles in the sphere lift to the double cover, so the degree of their intersection points can be computed by computing the degree of the the intersections of their lifts on the torus, using $dw^2$. The lift of any of the vanishing cycles has two components (two disjoint curves, two longitudes and two meridians, see Figure \ref{doublecoverfigure}). Namely, every vanishing cycle encircles two of the base points in the sphere, and since $p$ is a double cover, the monodromy of $p$ around two of the branched points is the identity. On the torus, the lifts of the vanishing cycles intersect in four points pairwise related by the covering involution $w\mapsto -w$ . Within one of these pairs, the two points are related by translation. As $w$ is translation invariant, all intersection points must have the same degree. 

\qed

With the gradings of \ref{lemmadegreegens}, the non-trivial entries of the spectral sequence \eqref{spsequence}, generated by the Floer cohomology groups of the vanishing cycles in the intersection fibres, are given by (see \eqref{spsequence})

\begin{align}\label{spsequencerp2}
E_1^{pq}=  \HF^{p+q}\left(V_{z_0,2}[ p(d-2)], V_{z_0,0}\right) \cong  \HF^{p+q}\left(V_{z_0,2} [ p ], V_{z_0,0}\right) \cong \HF^{q} (V_{z_0,2},V_{z_0,0}),	\end{align}
where $d=3$ comes from \eqref{pencilrp2}.

\begin{eqfloat}
	\begin{equation*}
	\begin{tikzpicture}
	\matrix (m) [matrix of math nodes,
	nodes in empty cells,nodes={minimum width=5ex,
		minimum height=4ex,outer sep=-5pt},
	column sep=0ex,row sep=0ex]{
		&&&&&&&q \\ 
		&0& 0&0 &0 &0 &  &\vdots\\
		& 0 &  0  & 0& 0 &0 & & 1\\
		&\Zmod^2 &\Zmod^2  & \Zmod^2 & \Zmod^2& \Zmod^2 &&0 \\
		&&&&&&& \\
		p   &	-4&-3&-2&-1&0&&\\
	};
	
	\draw[thick] (m-1-7.east) -- (m-6-7.east) ;
	\draw[thick] (m-5-1.north) -- (m-5-8.north);
	\end{tikzpicture}
	\end{equation*}
	\caption{The shape of the spectral sequence for $k=5$.}
\end{eqfloat}

Note that there is room for a nontrivial differential $d_1 \colon E^{(p,q)}\to E^{(p+1,q)}$, so it is not possible to directly compute the Floer cohomology groups $\HF^*(\varphi^k(\Delta_{2}), \Delta_0)$. However, using \cite{evans}, we can infer the following.

\begin{cor}\label{corollaryreal}
	The symplectomorphism $\varphi \in \pi_0(\Symp_{ct}(T^*\RPtwo))$ is isotopic to a power of the projective twists $\tau^k_{\RPtwo}$, $k\in \Z^*$.
\end{cor}
\proof By Proposition \ref{proprp2}, $\varphi$ is a non-trivial compactly supported symplectomorphism of $T^*\RPtwo$.
The symplectic mapping class group $\Symp_{ct}(T^*\RPtwo)$ is known to be generated by $\tau_{\RPtwo}$ (\cite{evans}) so $\varphi$ is isotopic to a power $\tau_{\RPtwo}^k$, $k \in \Z$.
\endproof 

We expect the rank of $\HF(\varphi^k(\Delta_2), \Delta_0)$ to increase linearly with $k$ (see below); if that was indeed the case, $\varphi$ would indeed be isotopic to the projective twist, or its inverse. 

\subsubsection{Expected results}\label{expectedrp2}

Let $\N_{k}:=\{ 0, 1 \dots , k-1 \}$. The Floer cohomology group are expected to be of the shape

\begin{equation}\label{expectedhfrp2}
\HF^*(\varphi^k(\Delta_2),\Delta_{0}; \Zmod)=\left\lbrace
\begin{array}{ll}
\Zmod\oplus \Zmod	&*=- \N_k\\
0 & \text{otherwise}.
\end{array} \right.
\end{equation}

Note that if all the differentials of the spectral sequence \eqref{spsequencerp2} vanish, \eqref{expectedhfrp2} is the limit to which the sequence converges to. The nature of the above prediction rests in a series of claims, some of which are still in conjectural state in the literature. We sketch the reasoning below. 

Consider the symplectomorphism of the total space $\Phi_{2\pi} \in \Symp(T^*\RPtwo)$ we defined in Section \ref{construction}.

Assume that $\Delta_{0}$, $\Delta_2$ have been suitably deformed into Lagrangians with Legendrian boundary at infinity, to become objects of the wrapped Fukaya category $\mathcal{W}(T^*\RPtwo;\Zmod)$. 

Then the following is expected (see \cite[Remark 3.1]{maysei}, \cite[Appendix]{bee}, \cite[Lemma 3.37]{gps}, \cite{aglgmodels}) \begin{align}\label{conjecturehw}
\varinjlim_k \HF^*(\Phi_{2\pi}^k(\Delta_2), \Delta_0;\Zmod)\cong  \text{HW}^*(\Delta_2, \Delta_0;\Zmod).
\end{align}

Moreover, if we identify the thimbles with cotangent fibres $T_{q_1}^*, T_{q_2}^* \in T^*\CPtwo$, we obtain 

\begin{align}\label{wrapped1}
\varinjlim_k\HF^*(\Phi_{2\pi}^k(\Delta_2), \Delta_0)=\text{HW}^*(\Delta_{2}, \Delta_0) \cong \text{HW}^*(T_{q_1}^*, T_{q_2}^*)\cong  \text{HW}^*(T_{q_1}^*, T_{q_1}^*).
\end{align}
where the last isomorphism is induced by invariance under Hamiltonian isotopies of wrapped Floer cohomology (\cite[Section 3.1]{absei}). 

Let $\Omega\RPtwo$ be the based loop space of $\RPtwo$. By \cite{absch}, there is an isomorphism $\text{HW}^*(T_{q_1}^*, T_{q_1}^*;\Zmod)\cong H_{-*}(\Omega\RPtwo)$; which, combined with \eqref{wrapped1} yields the conjectural relation
\begin{align}
\varinjlim_k \HF^*(\Phi_{2\pi}^k(\Delta_2), \Delta_0) \cong H_{-*}(\Omega\RPtwo;\Zmod).
\end{align}

\begin{lemma}
	The homology of the based loop space $\Omega \RPtwo$ is given by \begin{align}
	\forall *\in  \N, \ \	H_{*}(\Omega \RPtwo;\Zmod) \cong H_*(\Omega S^2)\oplus H_*(\Omega S^2) \cong		\Zmod\oplus \Zmod.
	\end{align}

\end{lemma}
\proof Consider the fibration $\Zmod \to S^2 \to \RPtwo$, and let $P\RPtwo$ the space of based paths on $\RPtwo$.
Pulling back the path-loop fibration $\Omega \RPtwo \to P\RPtwo \to \RP$ to $S^2\to \RPtwo$ yield another fibration $\Omega \RPtwo \to \Zmod \to S^2$, see for example \cite[Section 5]{atcourse} (alternatively, this can be deduced from the homotopy lifting property). Iterating this process we obtain $\Omega S^2 \to \Omega \RPtwo \to \Zmod$. 

In this situation, despite having a disconnected base $\Zmod$ (and $\Omega \RPtwo$ has two components), the fibration does have homotopy equivalent fibres, an explicit homotopy equivalence being the following. 
Let $\gamma$ be a fixed based loop in $\RPtwo$ representing the non-trivial homology class. Then, concatenation with $\gamma$ defines a map from a component of $\Omega \RPtwo$ to the other; concatenating twice gives a map from each component to itself.
Therefore, $H_*(\Omega \RPtwo) \cong H_*(\Omega S^2)\oplus H_*(\Omega S^2)$, and $H_{*}(\Omega S^2)\cong k[x]$, for $x$ in degree $1$. (The homology of based loop spaces of spheres is known (\cite[p.235]{htpytheory}), and can be obtained by applying the Serre spectral sequence). \endproof

\begin{rmk}
	The expected result \eqref{expectedhfrp2} matches the computation of \cite[Theorem 1.3]{frauenschlenk}.
\end{rmk}

\subsection{$\CPtwo$-twist}\label{cp2twist}

Following the instructions of Section \ref{construction}, this time for the Lefschetz fibration $\pi\colon E_{\CPtwo}\to \C$ built in Section \ref{lfcp2}, we define a symplectomorphism $\varphi \in \Symp_{ct}(T^*\CPtwo)$.

Using the techniques of Section \ref{spectral}, we then compute the Floer cohomology of two vanishing thimbles, one of which acted upon by a power $\varphi^k$. We find that the rank of the Floer cohomology is given by $2k$, which implies that $\varphi$ has infinite order in $\Symp_{ct}(T^*\CPtwo)$. Comparing this outcome with results of \cite{frauenschlenk} and predictions from wrapped Floer cohomology computations, we conjecture the following.

\begin{conj}
	The symplectomorphism $\varphi$ is isotopic to the $\CPtwo$-twist.
\end{conj}

\subsubsection{A symplectomorphism $\varphi \in \Symp_{ct}(T^*\CPtwo)$ of infinite order.}	\label{symplecto}

Let $\pi \colon E_{\CPtwo} \to \C$ be the Lefschetz fibration of Section \ref{lfcp2}. Recall that the smooth fibres are Stein domains whose symplectic completion is exact symplectomorphic to the plumbing of spheres $T^*S^3\#_{S^1}T^*S^3$. Let $z_*\in \C$ be the base point and consider the vanishing paths $(\gamma_0, \gamma_2, \gamma_4)$, Lefschetz thimbles $(\Delta_{0}, \Delta_2, \Delta_4)$ and vanishing cycles $(V_0,V_2,V_4)$ as in Section \ref{lanterncp2}. By Lemma \ref{monodromycp2}, the ungraded monodromy preserves the vanishing cycles, and the graded version shifts their grading as $\phi(V_i) = V_i[-2]$. Then Corollary \ref{commutingrel} grants the well-definedness of the construction of Section \ref{construction}, that we use to build a compactly supported symplectomorphism $\varphi \in \Symp_{ct}(T^*\CPtwo)$ of the total space of the Lefschetz fibration $\pi\colon E_{\CPtwo} \to \C$. 

To compute the Floer cohomology group $\HF(\varphi^k(\Delta_4), \Delta_0; \Zmod)$, we utilise the spectral sequence of Proposition \ref{spectral}, and much of the discussion for the $\RPtwo$ case (Section \ref{hfrp2}) applies.

Over each intersection point $z \in I$, the intersection of the thimbles is determined by the intersection of their associated vanishing 3-spheres $(V_{z,4}, V_{z,0})=(\Delta_4\cap \pi^{-1}(z), \Delta_0\cap \pi^{-1}(z))$. The latter intersect cleanly in a circle $V_{z,4}\cap V_{z,0} \cong S^1$. 
In a neighbourhood of the intersection locus, one can apply a perturbation by a Morse function with two critical points (one minimum and one maximum), so that after perturbing, the two vanishing cycles intersect in two points $\{ \xi_z^-, \xi_z^+ \} \subset \pi^{-1}(z)$, and the latter generate the Floer complex $\CF(V_{z,4}, V_{z,0})$ in each fibre $\pi^{-1}(z)$, $z\in I$.
The generators $\{ \xi_z^-,\xi_z^+ \}$ represent non-trivial cocycles in $\CF(V_{z,4}, V_{z,0})$, since the clean intersection $V_{z,4} \cap V_{z,0}$ has one component, and by exactness, the results of \cite[Theorem 3.4.11]{pozniak} apply (see also \cite[Theorem 3.1]{seidelknotted}). It follows that the Floer cohomology of the vanishing cycles is the standard cohomology (up to a grading shift $l \in \Z$ discussed below)  \begin{equation}\label{hffibre}
\forall z \in I, \ \HF^{*+l}(V_{z,4}, V_{z,0};\Zmod ) \cong H^*(S^1;\Zmod) \cong \Zmod\oplus\Zmod.
\end{equation}

By Lemma \ref{firstpowerhf}, we know that $\HF^*(\varphi(\Delta_4), \Delta_0;\Zmod)\cong \Zmod \oplus \Zmod$ generated by $ \{ \xi_0^-, \xi_0^+\}$. On the other hand, $\HF(\Delta_{4}, \Delta_0;\Zmod)=0$ so we obtain:
\begin{prop}
	The symplectomorphism $\varphi\in \Symp_{ct}(T^*\CPtwo)$ is a non-trivial element of the symplectic mapping class group.\qed 
\end{prop}

Let $\gamma_4^k=b_{2\pi}^k(\gamma_4)$ and $\Delta_{4}^k:=\Delta_{\gamma_4^k} \cong \varphi^k(\Delta_4)$ and $I:= \gamma_4^k(\R^+) \cap \gamma_0(\R^+)  =
\{ z_0, \dots ,z_{k-1} \}$. Equip the Lagrangian vanishing cycles in the smooth fibres and the Lagrangian thimbles and with a $\Z$-grading as in Section \ref{vcycles}.

Set $\ind(z_0)=0$ and $\ind(z_{j})=-2j$ for $z_{j}\in I$, as in Section \ref{arcs}. The elements $\{ \xi_j^-, \xi_j^+ \}  \subset \pi^{-1}(z)$ of any generating pair of the complex $\CF^*(V_{z, 4}, V_{z,0}; \Zmod)$ must be shifted in their degrees as $\deg(\xi_j^+)= \deg(\xi_j^-)\pm 1$. Namely, 
this shift is established by the Morse perturbation applied above (see \cite[(4.4)]{seidelknotted}). Set $\deg(\xi_0^-)=-1$, $\deg(\xi_0^+)=0$. This fixes \begin{align*}
\HF^*(V_{z_0,4}, V_{z_0,0}; \Zmod) \cong \Zmod, \ *=-1,0
\end{align*}

\noindent This choice determines the grading of the Floer cohomology groups $\HF^*(V_{z, 4}, V_{z,0}; \Zmod)$, for every $z \in I$. Inserting $d=6$ (from \eqref{pencilcp}) in the formula \eqref{spsequence}, we obtain, for $-(k-1) \leq p \leq 0$,
\begin{align}\label{spsequencecp2}
E_1^{pq}=\HF^{p+q}(V_{z_{(-p)},4}, V_{z_{(-p)},0})\cong \HF^{p+q-p(d-2)}(V_{z_0,4}, V_{z_0,0}) \cong  \HF^{q-3p}(V_{z_0,4}, V_{z_0,0}).	\end{align}

\begin{eqfloat}
	\begin{equation}\label{spseq2}
	\begin{tikzpicture}

	\matrix (m) [matrix of math nodes,
	nodes in empty cells,nodes={minimum width=5ex,
		minimum height=5ex,outer sep=-5pt},
	column sep=1ex,row sep=1ex]{
		p   &-3 & -2   &-1 &0 &  &   \\
		&&&&&& \\
		&0&0 &0 &\Zmod  &&0\\
		&0 &0  & 0&\Zmod&&-1\\
		&0&0& 0& 0& &-2\\
		&0&0&\Zmod & 0&&-3\\
		&0&  0  &\Zmod&0  &&-4 \\
		&0& 0 &0&0& &-5\\
		&0&\Zmod&0&0&   &-6 \\
		&0&\Zmod&0&0&& -7\\
		&&&&&&q\\
		\\};
	\draw[thick] (m-1-6.east) -- (m-11-6.east) ;
	\draw[thick] (m-2-1.north) -- (m-2-7.north) ;

	\end{tikzpicture}
	\end{equation}
	\caption*{The shape of the spectral sequence for $k=3$.}
\end{eqfloat}

The difference in total degree between generators of different fibres is always greater than one: by Lemma \ref{degreeshift}, we know it is at least $s_{\varphi}=d-2=6-2=4$. Therefore, there is no non-trivial ``horizontal'' differential (counting pseudoholomorphic curves between different fibres). This can also be noticed by the position of the generators in the first page of the sequence, in Diagram \eqref{spseq2} (the total degree of a generator is obtained by adding $p$ and $q$).
With \eqref{hffibre} it follows that all differentials are trivial, and the sequence collapses at the first page, and we have the ungraded result $\HF(\Delta_{4}^k, \Delta_0)\cong \oplus_{j=0}^{k-1} \HF(V_{z_j,4}, V_{z_j,0})$. Together with the grading computations above, this implies the following.

\begin{thm}\label{propcp2}Let $\N_{k}:=\{ 0, 1 \dots , k-1 \}$.
	The graded Floer cohomology groups are given by
	\begin{equation}\label{hfcp2}
	\HF^*(\varphi^k(\Delta_4),\Delta_{0}; \Zmod)=\left\lbrace
	\begin{array}{ll}
	\Zmod	&*=-4 \N_k, -4\N_k-1\\
	0 &.
	\end{array} \right.
	\end{equation}
	In particular, omitting the gradings, $\HF(\varphi^k(\Delta_4),\Delta_{0}; \Zmod) \cong (\Zmod\oplus \Zmod)^k$. \end{thm}\endproof

\begin{cor}\label{corollaryorder}
	The symplectomorphism $\varphi$ has infinite order in $\pi_0(\Symp_{ct}(T^*\CPtwo))$.
\end{cor}

\begin{rmk}
	The symplectic mapping class group $\pi_0(\Symp_{ct}(T^*\CPtwo))$ has not yet been computed. In particular, it is not known whether it is solely generated by the standard projective twist $\tau_{\CPtwo}$.
	For example, there could be other potential generators such as $\CPtwo$-twists which are not Hamiltonianly isotopic to $\tau_{\CPtwo}$. However, as shown in \cite{bct}, such phenomena have only been observed in dimensions $n\geq18$, so we cannot draw further conclusions.
\end{rmk}

\subsubsection{Expected results}\label{expectedcp2}

As for the real projective twists, we compare our results with the literature. First note that if we identify the Lefschetz thimble with cotangent fibres of the cotangent bundle $T^*\CPtwo$, Theorem \ref{propcp2} matches with the computations of \cite{frauenschlenk}.

The Floer cohomology of a cotangent fibre $T^*_{q_1}$ of $T^*\CPtwo$ twisted by a projective twist and another (untwisted) cotangent fibre $T^*_{q_2} \subset T^*\CPtwo$ was computed in \cite[Theorem 2.13]{frauenschlenk}. Let $q_1,q_2 \in \CPtwo$ be two distinct points and $T^*_{q_1},T^*_{q_2} \in T^*\CPtwo$ their associated cotangent fibres. Then, according to \cite{frauenschlenk}, \begin{align}\label{frschlmapcp}
\HF(\tau_{\CPtwo}^k(T^*_{q_1}),T^*_{q_2}; \Zmod) \cong (\Zmod\oplus \Zmod)^k.
\end{align}

We can also apply the observations of Section \ref{expectedhfrp2} to get

\begin{equation}
\begin{split}
\varinjlim_k\HF^*(\Phi_{2\pi}^k(\Delta_4), \Delta_0;\Zmod) &\cong \text{HW}^*(\Delta_{4}, \Delta_0;\Zmod) \cong \label{expectedhwcp2} \\
\cong \text{HW}^*(T_{q_1}^*,T_{q_2}^*;\Zmod)&\cong  \text{HW}^*(T_{q_1}^*, T_{q_1}^*;\Zmod).
\end{split}
\end{equation}

In this case, $\CPtwo$ is not a Spin manifold, but since we work over a field of coefficient two, the isomorphism $\text{HW}^*(T_{q_1}^*, T_{q_1}^*;\Zmod) \cong H_{-*}(\Omega \CPtwo; \Zmod)$ from \cite{absch} still holds. We therefore obtain the expected isomorphism
$$\varinjlim_k\HF^*(\Phi_{2\pi}^k(\Delta_4), \Delta_0;\Zmod)\cong H_{-*}(\Omega \CPtwo;\Zmod).$$

The homology of the loop space $\Omega \CPtwo$ can be computed by considering the fibration $\Omega S^5 \to \Omega \CPtwo \to S^1$ obtained by ``looping twice'' the standard fibration $S^1 \to S^5 \to \CPtwo$ and (as we did in Section \ref{expectedrp2}) and applying the Serre spectral sequence (see \cite[Chapter 16 (a)]{htpytheory}).

\begin{lemma}
	For a ground field $k$ of characteristic $2$, $H_{*}(\Omega \CPtwo)\cong 
	H_*(S^1)\otimes H_*(\Omega S^5)$ with generators $x$ in degree $+4$ and $y$ in degree $+1$, and $y^2=0$.
	
	\qed 
\end{lemma}

	\section{Implementations}\label{hp2}
In this last section, we discuss a set of further potential developments (some more conjectural than others) that can be catalised by our investigations.

\subsection{The $\HPtwo$-twist}
As a first corollary, we explain how the main ideas and arguments can be applied (with some gaps) to illustrate a plausible local model for the $\HPtwo$-twist. The model should be built from an appropriate Lefschetz fibration on the cotangent bundle of the quaternionic projective plane.

By the works of Johns, there exists an exact Lefschetz fibration with three singular fibres, and whose total space is exact symplectomorphic to a disc cotangent bundle of $\HPtwo$ (\cite[Theorem A]{ johns1, johns}). The generic smooth fibre of such a fibration is a clean plumbing of spheres $T^*S^7\#_{S^3}T^*S^7$ along $S^3\subset S^7$, such that the three vanishing cycles of the fibration are the two $7$-spheres of the core of the plumbing and their Morse--Bott surgery (\cite[Section 4.i]{johns1}). 

In the corresponding algebro-geometric picture, the fibration should be induced by a Lefschetz pencil of hyperplane sections on $Gr_{\C}(6,2)$, the Grassmannian of $2$-dimensional linear subspaces of a complex $6$-dimensional vector space ($Gr_{\C}(6,2)$ is a projective variety with a Plücker embedding into $\C\PP^{14}$).

This is inferred by the following observations relative to the real and complex projective planes. In both cases, the projective variety $X$ they embed in as Lagrangians and on which the Lefschetz pencil is defined, can be considered their `complexification' (as a complex projective algebraic variety). This means that the variety $X$ admits an anti-holomorphic involution fixing $\APtwo$ so that the latter embeds in $X$ as the real locus. For $\RPtwo\subset \CPtwo$ the involution is given by conjugation, whereas in the case $\CPtwo \subset \CProduct$ the involution is given by conjugation on the second factor.
The same happens for the quaternionic projective plane $\HPtwo$, which embeds as a Lagrangian into $Gr_{\C}(6,2)$ (whose real dimension is indeed $16$), via an an anti-holomorphic involution on $Gr_{\C}(6,2)$ preserving it (see \cite[Section 6]{atiyahberndt}).

On the other hand, the the smooth fibre of such a pencil should coincide with the space of geodesics of $\HPtwo$ (according to Proposition \ref{audinprop}).

Assuming the existence of such a fibration, we could build an element $\varphi \in \Symp_{ct}(T^*\HPtwo)$ from the construction of Section \ref{construction}, and run the same Floer theoretical computations with two Lefschetz thimbles $\Delta_0, \Delta_{8} \subset T^*\HPtwo$. Assuming all the differentials would be forced to vanish due to degree reasons, the Floer cohomology groups would be given by\begin{align}
\HF(\varphi(\Delta_{8}), \Delta_{0};\Zmod)\cong \bigoplus_{i=1}^k \HF(V_{p_i,8}, V_{p_i,0};\Zmod) \cong \bigoplus_{i=1}^k H^*(S^3; \Zmod) \cong (\Zmod\oplus \Zmod)^k.
\end{align}

The second isomorphism is justified as in \eqref{hffibre}, since the two $7$-spheres are exact and meet cleanly in a $3$-sphere.

This outcome would match the existing literature results concerning the local projective twist $\tau_{\HPtwo} \in \pi_0(\Symp_{ct}(T^*\HPtwo))$.

\begin{lemma}[{{\cite{frauenschlenk}}}] Let $x,y \in \HPtwo$ be distinct points and $T_x,T_y \in T^*\HPtwo$ their associated cotangent fibres. Then
	\begin{align}
	\HF(\tau_{\HPtwo}^k(T_x),T_y; \Zmod) \cong (\Zmod\oplus \Zmod)^k. 
	\end{align}
\end{lemma}

\subsection{The $\RPtre$- twist}\label{rp3}

It should also be possible to apply the construction of Section \ref{chaptermodels} to provide an adequate local model for the $\RPtre$-twist.

Consider a pencil of quadrics in $\CPtre$, i.e a pencil generated by hyperplane sections of the line bundle $\mathcal{L}:=\mathcal{O}_{\CPtre}(2)$. Any smooth quadric $\Sigma\subset \CPtre$ is isomorphic to the product $\Sigma\cong \CPone \times \CPone$, via the Segre embedding, and two quadrics in $\CPtre$ intersect generically in a $2$-torus. Namely, since we have $\mathcal{O}_{\CPtre}(2)|_{\CPone\times \CPone}=\mathcal{O}_{\CPone\times \CPone}(2,2)$, the base locus has bidegree $(2,2)$, and the genus formula $g=(d_1-1)(d_2-1)$ yields $g=1$. 

Hence, a pencil described above has base locus $B$ a divisor of bidegree $(2,2)$, and $4$ singular fibres. Recall that the number of singular fibres follows from 
$$\chi(\CPtre)=2\chi(\CPone \times \CPone)- \chi(B)-r.$$
and $\chi(B)=0$. 

By Example \ref{examplerp}, there is a symplectomorphism $\CPtre \setminus (\CPone\times \CPone) \cong \mathring{D}_{\eps}T^*\RPtre$, for an open disc bundle $\mathring{D}_{\eps}T^*\RPtre$, $\eps>0$. The fibration $\CPtre \setminus (\CPone\times \CPone) \to \CPone$ derived from the pencil therefore induces an exact Lefschetz fibration \begin{align}\label{lfrp3}
\pi\colon T^*\RPtre \to \C.
\end{align}
The smooth fibre of $\pi$ is exact symplectomorphic to the completion of $M:= \CPone\times \CPone \setminus U_{B/\Sigma}$, where $U_{B/\Sigma}$ is an open neighbourhood of $B\subset \CPone \times \CPone$ (which is modelled on an open subbundle of the normal bundle to $B\subset \Sigma$).

Such a fibration admits four singular fibres with four associated Lagrangian vanishing spheres, that we denote by $V_0, V_1, V_2, V_3 \subset M$ (relative to a choice of vanishing paths and thimbles).

The global monodromy of the fibration $\pi$ is a product of the twists along the four vanishing $2$-spheres and satisfies a generalised lantern relation
\begin{align}\label{monosubsrp3}
\tau_{V_0}\tau_{V_1}\tau_{V_2}\tau_{V_3} \simeq \tau_{V}
\end{align}
where $\tau_{V}$ is the fibred twist along $\partial M\to B$.

The construction of Section \ref{construction} can be applied to the fibration $\pi\colon T^*\RPtre\to \C$ above, to obtain a compactly supported symplectomorphism $\varphi \in \Symp_{ct}(T^*\RPtre)$. We conjecture that $\varphi$ is isotopic to the projective twist $\tau_{\RPtre}\in \pi_0(T^*\RPtre)$.

Note that in this case we can make no direct Floer theoretical conclusions as we do not know the intersection pattern of the vanishing spheres.

\subsection{Monodromy substitution and symplectic fillings of unit cotangent bundles}\label{monodromysubst}\label{substit}\label{surgery}

Let $\pi\colon E \to D$ be an exact Lefschetz fibration over the disc (in the sense of \cite[(15d)]{seidelbook}), with smooth fibre a Liouville domain $(M, \omega=d\lambda)$ and total monodromy $\phi \in \Symp_{ct}(M)$. After rounding off the corners of the total space, its boundary $\partial E$ admits a contact structure, which (by construction) is compatible with an open book decomposition $(M,\lambda, \phi)$ with page $M$ and monodromy $\phi$.

Assume that as before, $\pi$ is induced by a Lefschetz pencil on a complex projective variety $X$ seen as a Kähler manifold $(X, \omega_X)$, with closed fibre a Kähler submanifold $(\Sigma, \omega_{\Sigma})$ and base locus $(B, \omega_B) \subset (\Sigma, \omega_{\Sigma})$. Let $\mathcal{L}\to X$ be the line bundle defining the pencil, satisfying $c_1(\mathcal{L})=[\Sigma]$.

Then, as explained in Section \ref{monodromy}, the total space $(E,\Omega_E)$ is exact symplectomorphic to the complement $X\setminus U_{\Sigma/X}$ of an open neighbourhood $U_{\Sigma/X}$ of the smooth fibre $\Sigma \subset X$.
Similarly, the smooth fibre $(M, \omega)$ is exact symplectomorphic to the complement $\Sigma\setminus U_{B/\Sigma}$ of an open neighbourhood $U_{B/\Sigma}$ of the base locus $B\subset \Sigma$. Recall that the neighbourhoods $U_{\Sigma/X}$ and $U_{B/\Sigma}$ are disc subbundles of the symplectic normal bundles to $\Sigma \subset X$ and to $B\subset \Sigma$, and hence modelled on the restrictions $\L|_{\Sigma}$ and $\L|_{B}$ respectively.

In this section we discuss an operation (called \emph{monodromy substitution}) which utilises the generalised lantern relation \eqref{generalantern} to derive a Morse--Bott--Lefschetz (MBL) fibration $\pi'\colon E' \to D$ with fibre symplectomorphic to $(M,\omega)$, boundary contactomorphic to $\partial E$ but whose total space $E'$ is topologically distinguished from $E$.

The monodromy $\phi\in \Symp_{ct}(M)$ of the fibration $\pi$ is a product of positive powers of Dehn twists, which, by the generalised lantern relation, is isotopic to a fibred Dehn twist $\phi \simeq \tau_{V}$ along the coisotropic $V=\partial M\to B$ (see Remark \ref{coisotropicboundary} and Lemma \ref{lantern1}). On the other hand, by the general results summarised in Section \ref{mbldiscussion}, the data $(M, \omega)$, $V \subset M$ also define a MBL fibration $\pi'\colon E'\to D$ with fibre (exact) symplectomorphic to $(M, \omega)$ and monodromy the fibred twist $\tau_V \in \Symp_{ct}(M)$.

The total space $E'$ of the resulting MBL fibration can be read off by using Biran's decomposition of Theorem \ref{birandecompo}: it is symplectomorphic to an open disc subbundle $U_{\Sigma/X}^*\subset N_{\Sigma/X}^*$ of the dual symplectic normal bundle to $\Sigma\subset X$ (this is modelled on the restriction $\mathcal{L}^*|_{\Sigma}$ with $c_1(\mathcal{L}^*|_{\Sigma})=-[B]$), see also \cite{oba}. 

Since the boundary of a MBL fibration $\pi'\colon E' \to D$ with smooth fibre $(M, \omega=d\lambda)$ and monodromy $\phi$ admits a contact structure compatible to the open book $(M,\phi, \lambda)$ (\cite[A.2]{oba}), by construction the contact boundaries $\partial E$ and $\partial E'$ are both contactomorphic to the same open book $(M, \lambda, \phi)$.

Therefore, the contact manifold defined by the open book $(M,\lambda, \phi)$ has at least two fillings: one given by the total space of the exact Lefschetz fibration $E$, and the other by the total space of the MBL fibration $E'$. The former is an exact filling, while the latter is a strong filling (\cite[Proposition A.3]{oba}) which in general contains rational curves.

\begin{thm}{\cite[Corollary A.4]{oba}}\label{strongfilling} Let $(M, \omega=d\lambda)$ be a Liouville domain, $\phi \in \Symp_{ct}(M)$. Suppose there is a collection of spherically fibred coisotropic submanifolds $(V_1, \dots , V_k) \subset M$ such that $\phi \simeq \tau_{V_1}\dots \tau_{V_k}$ in $\pi_0(\Symp_{ct}(M)$. Then, a contact structure compatible with the open book $(M, \lambda, \phi)$ is strongly fillable.
	\qed 
\end{thm}

So the choice of identification of the monodromy, either to a product of twists, or to a fibred Dehn twist via the generalised lantern relation, indicates the presence of two (possibly distinguished) symplectic fillings, and passing from one identification to the other is called a \emph{monodromy substitution} (see \cite{engu, emvhm}).

\begin{definition}\cite[4.3]{oba}\label{defsubstitution}
	Let $\pi_i \colon E_i \to \C$, $i=1,2$ be two (Morse--Bott--) Lefschetz fibrations with fibre a symplectic manifold $(M, \omega)$ and let $(V_{i,1}, \dots , V_{i,k_{i}})$ the collection of (coisotropic) vanishing cycles of $\pi_i$ for some distinguished basis. The two fibrations are related by a monodromy substitution if $\exists 1 \leq j_0 \leq \min\{ k_1,k_2 \}$ such that $V_{1,j}=V_{2,j}$ for any $j <j_0$ and there is an isotopy in $\pi_0(\Symp_{ct}(M))$, $\tau_{V_{1,j_0}} \circ \cdots \circ  \tau_{V_{1,k_1}} \simeq \tau_{V_{2,j_0}} \circ \cdots \circ  \tau_{V_{2,k_2}}$ (then the global monodromies are also isotopic).
\end{definition}

A monodromy substitution can be performed in a general context, and when applicable the resulting fillings are potentially non-equivalent (see \cite[Example 3.13, Proposition 4.2]{oba}). We will see applications in Sections \ref{fillingrp2}, \ref{fillingrp3}, \ref{fillingcp2}.

To examine the topology of the new filling, we need to consider the decomposition (given by Theorem \ref{birandecompo}) of the Kähler manifold $(X,\omega_X)$ that generates the Lefschetz pencil/fibration. In every example under study we are in a fortuitous setting in which the decomposition has the shape\begin{align}
X\cong \Xi \cup U_{\Sigma/X},
\end{align}
for a smooth Lagrangian $\Xi \subset X$ endowed with a Riemannian metric with periodic geodesic flow and an open disc subbundle of the symplectic normal bundle $U_{\Sigma/X} \subset N_{\Sigma/X}$ to the smooth fibre of the pencil.
On one hand, it yields the total space of an exact symplectic fibration induced by the pencil, which (by Proposition \ref{audinprop}) is exact symplectomorphic to a disc cotangent bundle $(D_{\eps}^*\Xi, d\lambda_{\Xi})$, whose boundary is contactomorphic to the sphere cotangent bundle $(ST^*\Xi, \xi_{std})$, with the canonical contact structure obtained from the Liouville form $\xi_{std}:=\ker(\lambda_{\Xi}|_{ST^*\Xi})$ (see Section \ref{examplescotangent}).

On the other hand, the monodromy substitution replaces the standard filling of $(ST^*\Xi, \xi_{can})$ with a subbundle $U_{\Sigma/X}^*\subset N^*_{\Sigma/X}$ of the dual symplectic normal bundle to the hypersurface $\Sigma \subset X$, yielding the total space of a MBL fibration with smooth fibre $(M, \omega)$ and $S^1$-fibred coisotropic $V=\partial M \to B$.

\subsubsection{Known example: two Stein fillings of $ST^*\RPtwo$}\label{fillingrp2}

Consider the Lefschetz fibration $\pi\colon T^*\RPtwo \to \C$ of Section \ref{pencil1} (recall it is generated by a Lefschetz pencil of conics on $\CPtwo$, i.e sections of $\O_{\CPtwo}(2)$). 
After smoothing the corners, the boundary of the total space is the unit cotangent bundle $(ST^*\RPtwo, \xi_{std})$, which, equipped with the standard contact structure, is contactomorphic to the lens space $(L(4,1), \overline{\xi}_{S^3})$, where $\overline{\xi}_{S^3}$ is the structure induced by the standard contact structure on $S^3$. By Theorem \ref{strongfilling}, the lantern relation \begin{align}\label{lanternfilling}
\tau_{V_0} \tau_{V_1} \tau_{V_2} \simeq  \tau_{d_1} \tau_{d_2}\tau_{d_3}\tau_{d_4}
\end{align}
indicates that the lens space $(ST^*\RPtwo, \xi_{std})$ admits two fillings. On one hand, the total space of the Lefschetz fibration $\pi$ yields the standard filling, i.e a disc cotangent bundle of $\RPtwo$.
The other filling is presented as the total space of a MBL (which in this case is also Lefschetz) fibration, with fibre isomorphic to a 4-punctured sphere, and monodromy given by the product of Dehn twists around the boundary circles of the fibre (the RHS of \eqref{lanternfilling}). By the decomposition \eqref{decomporpn}, the total space of the latter is given by the dual disc normal bundle to the quadric $\Sigma \subset \CPtwo$, a subbundle $U_{\Sigma/\CPtwo}^* \subset N^*_{\Sigma/\CPtwo}$.  In this case, it is easy to see that $N^*_{\Sigma/\CPtwo} \cong \mathcal{O}_{\CPone}(-4)$.

Note that $H_*(D^*\RPtwo;\Q)\cong H_*(B^4;\Q)$, where $B^4$ is the $4$-ball. Therefore, as shown in \cite{engu}, a monodromy substitution via the isotopy \eqref{lanternfilling} can be geometrically interpreted as a rational blowdown, in which a neighbourhood of a sphere with self intersection $-4$ is cut out and replaced by a rational homology ball.

These two distinct fillings coincide with the two minimal Stein fillings of $(ST^*\RPtwo, \xi_{std})\simeq (L(4,1), \overline{\xi}_{S^3})$ found in \cite{mcduff}.

\subsubsection{Fillings of $ST^*\RPtre$}\label{fillingrp3}
We can apply the principles of Section \ref{monodromysubst} using the Lefschetz fibration $\pi\colon T^*\RPtre \to \C$ constructed in Section \ref{rp3} to obtain an alternative symplectic filling of $ST^*\RPtre$.

Let $\pi\colon T^*\RPtre \to \C$ be the Lefschetz fibration \eqref{lfrp3}, with smooth fibre $M \cong (\CPone\times \CPone) \setminus U_{B/\Sigma}$, monodromy $\phi:=\tau_{V_1}\tau_{V_2}\tau_{V_3}\tau_{V_4} \in \Symp_{ct}(M)$ (see Section \ref{rp3}).

\begin{thm}\label{fillingrp3thm}
The monodromy substitution performed using the relation \eqref{monosubsrp3} yields a strong symplectic filling of the contact manifold $(ST^*\RPtre, \xi_{std})$, that is not a Stein filling.
\end{thm}
\proof As in the other examples, we use the lantern relation to perform a monodromy substitution: between the Lefschetz fibration (restricted over the disc) $\pi|_D\colon T^*\RPtre|_D \to D$ with 
fibre $M$, and monodromy $\phi \in \Symp_{ct}(M)$ 

and a MBL fibration $\pi'\colon E' \to D$ with fibre $M$ and monodromy the fibred tiwst along $\partial M\to B$ (isotopic to $\phi$)

The total space $E'$ is a strong filling of $ST^*\RPtre$ obtained as the dual normal bundle to $ \CPone \times \CPone \subset \CPtre$,
\begin{align}\label{fillinge'}
E'\cong \mathcal{O}_{\CPone\times \CPone}(-2,-2).
\end{align}
The topology of $E'$ reveals that this filling is distinguished from the standard cotangent bundle: $H^4(E';\Q)$ is of rank $2$, while $H^4(D^*\RPtre;\Q)$ has rank zero. Moreover, the filling \eqref{fillinge'} is not exact as the zero section is embedded as a symplectic submanifold.

\endproof

\subsubsection{Fillings of $ST^*\CPtwo$}\label{fillingcp2}
In the spirit of the previous two examples, we can use the generalised lantern relation \eqref{lantern2} to obtain the following.

\begin{prop}\label{fillingcp2thm}
	Consider the Lefschetz fibration $\pi\colon T^*\CPtwo \to \C$ as in \eqref{lfadjusted2}.
	The monodromy substitution performed using the relation \eqref{lantern2} yields a strong symplectic filling of the contact manifold $(ST^*\CPtwo, \xi_{std})$ which is not exact.
\end{prop}
\proof 
Consider the restriction $\pi\colon T^*\CPtwo|_D \to D$ of the Lefschetz fibration \eqref{lfadjusted2} with fibre exact symplectomorphic to $M\cong DT^*S^3\#DT^*S^3$, $B=\CPtwo\#3\overline{\CPtwo}$ and monodromy $\phi$. After rounding off the corners, $\pi$ has contact boundary contactomorphic to $(ST^*\CPtwo, \xi_{std})$. The lantern relation \eqref{lantern2} can be used to operate a monodromy substitution (Section \ref{substit}) to get a MBL fibration $\pi'\colon E' \to D$ with fibre $M$ and monodromy $\tau_V$, for $V=\partial M\to B$. The total space $E'$ is a strong filling (Theorem \ref{strongfilling}) given by the dual bundle 
\begin{align}\label{cp2filling}
E' =N^*_{Fl_3/\CProduct} \cong \mathcal{O}_{\CProduct}(-1,-1)|_{Fl_3}
\end{align}
in that replaces the Lagrangian $\CPtwo$ with $\mathcal{O}_{\CProduct}(-1,-1)$, the normal bundle over the Flag threefold $Fl_3\subset \CProduct$ with reverse orientation. 
This symplectic filling is not exact (the Flag $3$-fold is embedded as a symplectic submanifold).
\endproof

\bibliographystyle{alpha}




\end{document}